\documentclass[12pt]{amsart}
\usepackage[active]{srcltx}
\usepackage{atveryend}
\makeatletter
\let\origcitation\citation
\AtEndDocument{\def\mycites{}%
  \def\citation#1{\g@addto@macro\mycites{#1^^J}\origcitation{#1}}}
\AtVeryEndDocument{\newwrite\citeout\immediate\openout\citeout=\jobname.cit
  \immediate\write\citeout{\mycites}\immediate\closeout\citeout}
\makeatother
\usepackage{enumerate}
\usepackage{cancel}
\usepackage{mathdots}
\usepackage{epsfig}
\usepackage{hyperref}
\usepackage{amsmath}
\usepackage{amssymb}
\usepackage{graphics} 

\usepackage{yfonts}
\newtheorem{Pantachie}{Pantachie\ }
\newtheorem{Hahn Field}{Hahn Field}

\newtheorem{Surreal}{Surreal}
\newtheorem*{Definition *}{Definition}
\newtheorem*{Saturation 1 *}{Saturation for $V(\mathbb{^*R})$}
\newtheorem*{Saturation 2 *}{Saturation Axiom for $\mathbb{^*R}$}

\title[Infinitesimalist Theories of Continua]{Contemporary Infinitesimalist Theories of Continua and their late 19th- and early 20th-century forerunners}

\author{Philip Ehrlich}\thanks{{\bf{Acknowledgement}}. The idea of writing this historical overview grew out of our (2007). We are grateful to the editors for providing us this opportunity to do so, and especially for tolerating our going well beyond the designated page limits. During the course of writing the paper we have profited from correspondence with C. Ward Henson and H. J. Keisler about nonstandard analysis, and with Marta Bunge, Paolo Giordano, Wolfgang Bertram and Mikhail Katz about SDG/SDT, IDG, TDC and NSDG respectively, and we are grateful to them all for their comments and helpful responses to our queries. We are especially grateful to Marta Bunge for permitting us to quote from our correspondence.}

\begin{document}
\maketitle

\section{introduction}

 In the decades bracketing the turn of the twentieth century the real number system was dubbed the arithmetic continuum because it was held that this number system is completely adequate for the analytic representation of all types of continuous phenomena. In accordance with this view, the geometric linear continuum is taken to be isomorphic with the arithmetic continuum, the axioms of geometry being so selected to insure this would be the case. In honor of Georg Cantor and Richard Dedekind, who first proposed this mathematico-philosophical thesis, the presumed correspondence between the two structures is sometimes called the \emph{Cantor-Dedekind axiom}. Since that time, the Cantor-Dedekind philosophy has emerged as a pillar of standard mathematical philosophy that underlies the standard formulation of analysis, the standard analytic and synthetic theories of the geometrical linear continuum, and the standard axiomatic theories of continuous magnitude more generally. 

	Since its inception, however, there has never been a time at which the Cantor-Dedekind philosophy has either met with universal acceptance or has been without competitors. The period that has transpired since its emergence as the standard philosophy has been especially fruitful in this regard having witnessed the rise of a variety of constructivist and predicativist theories of real numbers and corresponding theories of analysis as well as the emergence of a number of alternative theories that make use of infinitesimals. Whereas the constructivist and predicativist theories have their roots in the early twentieth-century debates on the foundations of mathematics and were born from critiques of the Cantor-Dedekind theory, the infinitesimalist theories are intended to either provide intuitively satisfying (and, in some cases, historically rooted) alternatives to the Cantor-Dedekind conception that have the power to meet the needs of geometry, analysis or portions of differential geometry, or to situate the Cantor-Dedekind system of real numbers in a grander conception of an arithmetic continuum. 

	The purpose of this paper is to provide a historical overview of some of the contemporary infinitesimalist alternatives to the Cantor-Dedekind theory of continua.  Among the theories we will consider are those that emerge from nonstandard analysis, nilpotent infinitesimalist approaches to portions of differential geometry and the theory of surreal numbers. Since these theories have roots in the algebraic, geometric and analytic infinitesimalist theories of the late nineteenth and early twentieth centuries, we will also provide overviews of the latter theories and some of their relations to the contemporary ones.\footnote{Among the infinitesimalist conceptions we will not consider is the one arising from Alain Connes's  \emph{noncommutative (differential) geometry}. One of Connes's original motivations for developing noncommutative geometry was to apply geometric ideas and concepts to spaces that are intractable when considered from the usual set-theoretic framework of Riemannian geometry. The principal feature of noncommutative geometry is a novel conception of geometric space, encoded within an algebra, in general noncommutative, instead of a set of points, the geometry in classical spaces being recovered when the algebra is commutative. In Connes's  framework the ordinary differential and integral calculus is replaced by a new calculus of infinitesimals--\emph{the quantized calculus}--based on the operator formalism of quantum mechanics, and therein for each positive real $\alpha$, infinitesimals of order $\alpha$ (forming a two-sided ideal) are introduced in connection with compact operators on a Hilbert space. According to Connes, the framework that emerges from this perspective is ``a geometric space which is neither a continuum nor a discrete space but a mixture of both" (Connes 1994, p. 30; also see, Connes 1998, 2006).}  We will find that the contemporary theories, while offering novel and possible alternative visions of continua, need not be (and in many cases are not) regarded as replacements for the Cantor-Dedekind theory and its corresponding theories of analysis and differential geometry, but rather as companion theories or methods to be situated in the toolkits of mathematicians, or, as in the case of the surreal numbers, as a panorama from within which one may view the classical arithmetic continuum as the Archimedean portion of a transfinite recursively unfolding maximal arithmetico-tree-theoretic continuum whose elements are individually definable in terms of sets of standard set theory and whose first-order theory is that of the reals.

With the exception of \S 9 and a portion of \S11, all of the contemporary theories discussed below can be formalized in \emph{Zermelo-Fraenkel set theory with Choice} (ZFC) occasionally supplemented with one or another familiar assumption about infinite cardinals. In \S 9, where sets and proper classes are required, the underlying set theory is \emph{von Neumann-Bernays-G\"odel set theory with Global Choice} (NBG), where all proper classes have the ``cardinality" of the proper class $On$ of ordinals, and a sentence in the language of sets is true if and only if it is true in ZFC.\footnote{For NBG and its relation to ZFC, see (Smullyan and Fitting 2010).}

Since most of the earliest infinitesimalist theories of continua are based on \emph{non-Archimedean} (totally) ordered fields, rings and abelian groups, we will begin with these and their historical roots. Following standard practice, an ordered abelian group, an ordered ring or an ordered field $A$ is said to be \emph{Archimedean} if for all $x,y \in A$, where $0<x<y$, there is an $n \in \mathbb{N}$ such that $nx>y$.\footnote{The concepts of ordered class, divisible ordered abelian group, ordered commutative ring with identity, ordered field and positive cone of an ordered abelian group play prominent roles in portions of the text. For the reader's convenience, the definitions of these and some related concepts are collected below, where ``ordered" is understood to mean ``totally ordered".

An \emph{ordered class} is a structure  $\langle A,\leq\rangle$, where  $A$ is a class (a set or proper class) and  $\leq$ is a binary relation on $A$ that satisfies the conditions: $\forall xy(x\leq y \lor y \leq x)$; $\forall xyz((x\leq y \wedge y\leq z)\rightarrow x \leq z)$; $\forall xy((x \leq y \land y\leq x) \rightarrow x=y)$. If $x\leq y$ and $x \neq y$, we write $x<y$ or $y>x$. An \emph{ordered abelian group} (written additively) is a structure $\langle A,\leq ,+,0 \rangle$, where $\langle A,\leq \rangle$ is an ordered class and $+$ is commutative, associative binary operation on $A$ satisfying the conditions: $\forall x(x+0=x)$; $\forall x \exists y(x+y=0)$; $\forall xyz(x\leq y \rightarrow (x+z \leq y+z))$. An ordered abelian group $A$ is  \emph{divisible}, if for each $x \in A$ and each positive integer $n$, there is an $a \in A$ such that $na=x$. An ordered abelian group $\langle A,\leq ,\cdot,1 \rangle$ (written multiplicatively) is defined similarly using $\cdot$ and $1$ in place of $+$ and $0$, and divisibility is defined similarly by the condition: for each $x \in A$ and each positive integer $n$, there is an $a \in A$ such that $a^n=x$. An \emph{ordered commutative ring with identity} (or with \emph{unity}) is a structure $\langle A,\leq ,+, \cdot,0,1 \rangle$, where $\langle A,\leq ,+,0 \rangle$ is an ordered abelian group, $0 \neq1$, and $\cdot$ is commutative, associative binary operation on $A$ that satisfies the conditions: $\forall x(x\cdot 1=x)$; $\forall xyz[ (x \cdot (y+z)=(x \cdot y) +(x\cdot z)]$; $\forall xy((0\leq x \wedge 0\leq y) \rightarrow 0\leq xy )$. An \emph{ordered field} is an ordered commutative ring with identity that satisfies the further condition: $\forall x[x \neq0 \rightarrow \exists y (x \cdot y=1)]$. The ordered additive group of every ordered field is divisible. The \emph{positive cone} of an ordered additive abelian group $A$ is $\{x \in A: x > 0 \}$.

A nonzero element $x$ of an ordered commutative ring $A$ with identity is said to be a \emph{proper zero-divisor} if for some nonzero $y \in A$, $xy=0$. An ordered commutative ring with identity is said to be an \emph{ordered integral domain} if it has no proper zero-divisor. In the case of an ordered integral domain the condition $\forall xy((0\leq x \wedge 0\leq y) \rightarrow 0\leq xy )$ can be replaced by $\forall xy((0< x \wedge 0< y) \rightarrow 0< xy)$. An ordered field is an ordered integral domain.}

\section{the emergence of non-archimedean systems of magnitudes}

Even before Cantor (1872) and Dedekind (1872) had published the modern theories of real numbers that would be employed to all but banish infinitesimals from late 19th- and pre-Robinsonian, 20th-century analysis, Johannes Thomae (1870) and Paul du Bois-Reymond (1870-71) were beginning the process that would in the years bracketing the turn of the century not only establish consistent and relatively sophisticated theories of infinitesimals in mainstream mathematics but make them the focal point of great interest and a mathematically profound and philosophically significant research program. One theory grew out of the pioneering investigations of \emph{non-Archimedean geometry} of Giuseppi Veronese (1889, 1891, 1894), Tullio Levi-Civita (1892/93, 1898) and David Hilbert (1899), and led to the celebrated algebraico-set-theoretic work of Hans Hahn (1907). And another emerged from a parallel development of du Bois-Reymond's (1870-71, 1875, 1877, 1882) groundbreaking work on the \emph{rates of growth of real functions} and led in the same period to the famous works of G. H. Hardy (1910, 1912) and some less well known though important work of Felix Hausdorff (1907, 1909). Each of these research programs, which collectively gave rise to modern-day non-Archimedean mathematics, led to nonstandard theories of continua that have come to play important roles in the mathematics of our day, albeit not always in the guise of theories of continua.

Before turning to these research programs and the nonstandard theories of continua they gave rise to, we will first provide some historical background on the Archimedean axiom.

\section{the archimedean axiom}
	
	In his historically important paper \emph{Zur Geometrie der Alten, insbesondere \"uber ein Axiom des Archimedes}  (On the Geometry of the Greeks, in Particular, on the Axiom of Archimedes) Otto Stolz observed that:

\begin{quote}
            It has often been noted that Euclid implicitly used the principle: \emph{a magnitude can be so often multiplied that it exceeds any other of the same kind}....  Archimedes employed this principle as an explicit axiom in some of his works .... For brevity, we will therefore henceforth refer to this principle as the \emph{Axiom of Archimedes}.  To investigate whether or not this is a necessary proposition, requires us first to have agreement on a characterization of the concept of ``magnitude" (1883, p. 504). 
 \end{quote}           
  
 \noindent         
Such agreement was required for, as Stolz emphasized, the term ``magnitude" occurs in Euclid's \emph{Elements}, but he nowhere explains the concept. In response to his query, Stolz provided an axiomatization for the type of ordered additive systems of line segments occurring in the \emph{Elements}--an ordered, additive system of equivalence classes (of congruent line segments) constituting what we today call \emph{the positive cone of a divisible (additively written) ordered abelian group}, and therewith established the following result: \emph{whereas every such system of magnitudes that is continuous in the sense of Dedekind is Archimedean, there are such systems of magnitudes that are non-Archimedean}. To establish the existence of a latter such system he made use of an algebraic development of a fragment of du Bois-Reymond system of \emph{orders of infinity} that emerged from the latter's aforementioned work on the rates of growth of real functions (see \S 6).\footnote{See (Ehrlich 2012, pp. 24-27) for Stolz's construction.}  It was with these and related discoveries that Stolz (1883, 1884, 1885), Rodolfo Bettazzi (1890), Veronese (1889) and Otto H\"older (1901) laid the groundwork for the modern theory of magnitudes, the branch of late
19th- and early 20th-century mathematical philosophy that would, in the decades that
followed, evolve into the theories of Archimedean and non-Archimedean ordered algebraic systems, theories in which talk of systems of magnitudes would be replaced by talk of ordered abelian groups, ordered fields, the positive cones of such structures, and so on. 

The non-Archimedean systems of magnitudes studied in the just-mentioned works are additive structures that sometimes have modest multiplicative structures as well. Unlike the system of real numbers, however, none of them is an ordered field. Just as ordered fields of real numbers arose in conjunction with the study of Euclidean geometry, it was from the study of non-Archimedean geometry that non-Archimedean ordered fields emerged. It was also from non-Archimedean geometry that the first well-developed non-Archimedean theory of continua emerged.

\section{veronese's theory of continua}

Following Wallis's and Newton's incorporation of directed segments into
Cartesian geometry, it became loosely understood that given a unit segment
$AB$ of a line
$L$ of a classical Euclidean space, the collection of directed segments of
$L$ emanating
from $A$, including the degenerate segment $AA$ itself, constitutes an Archimedean ordered
field with $AA$ and $AB$ the additive and multiplicative identities of the field and addition
and multiplication of segments suitably defined. This idea was made precise by
Veronese (1891, 1894) and Hilbert (1899) in their pioneering works on the
foundations of geometry. Also emerging from these works, and inspired in part by the aforementioned work of Stolz, was the idea that it is
possible to construct an axiomatization for the central theorems of
Euclidean geometry that is \emph{independent of the Archimedean axiom} and
for which the aforementioned system of line segments in a model of
the geometry continues to be an ordered field; however, in those models
of the geometry in which the Archimedean axiom fails, the ordered fields
in question are non-Archimedean ordered fields.

Following Veronese,\footnote{For the remainder of this section, except when referring to an ordered field of segments, we employ the terms ``segment" and ``subsegment" without the modifiers ``directed" or ``degenerate" to refer to a nondegenerate line segment without a direction. } two segments $s$ and $s'$ are said to be \emph{finite relative to one another} if there are positive integers $m$ and $n$ such that $s' \prec ms$ and $s \prec ns'$, where $a \prec b$ indicates that $a$ is congruent to a proper subsegment of $b$, i.e. a subsegment of $b$ that is not identical to $b$. In accordance with this terminology, a system $S$ of segments may be said to be Archimedean if every pair of such segments is finite relative one another. If $S$ is non-Archimedean, there are segments $s$ and $s'$ that are not finite relative to one another, in which case $s$ is said to be \emph{infinitesimal relative to} $s'$ and $s'$ is said to be \emph{infinite relative to} $s$, if $s\prec s'$. Thus, $s$ is infinitesimal relative to $s'$ and $s'$ is infinite relative to $s$ just in case $ns \prec s'$ for all positive integers $n$. If $S$ is the positive cone of a non-Archimedean ordered field, then for each segment $s$ in $S$ there are segments $s'$ and $s''$ in $S$ such that $s$ is infinite relative $s'$ and infinitesimal relative to $s''$.

Unlike the analytic constructions of Hilbert, Veronese's construction of a non-Archimedean
ordered field of line segments is clumsy and quite complicated, though
to some extent the complexity is a by-product of what he is attempting
to achieve. After all, Veronese is not merely attempting to construct a
non-Archimedean ordered field of segments that is appropriate to the
geometry in question, but moreover an ordered field of such segments that models his novel theory of non-Archimedean
continua, and one that is synthetically constructed to boot. 

Veronese had a distinctive, well-developed philosophy of geometry that underlay his philosophy of the continuum. For the sake of space, we simply note that in addition to holding that our conception of geometric continua should be grounded in our conception of magnitude rather than in number or in points, Veronese held that just as it is compatible with geometrical intuition that the parallel postulate fails, it is compatible with our geometrical intuition that continua need not be Archimedean. It was for these reasons that, unlike Cantor and Dedekind, Veronese sought a segment-based, synthetic theory of geometric continua that is independent of the Archimedean condition.\footnote{While Veronese has a segment-based, rather than a point-set based, theory of continua, the notion of a point is a primitive of his geometric system.}

\begin{sloppypar} 
Despite the lack of elegance in its presentation and elements of
obscurity in its formulation, the theory of rectilinear continua developed
in Veronese's \emph{Fondamenti di Geometria} (1891) is a profound and relatively sophisticated scheme,
several of whose central concepts and ideas permeate the 20th-century
theory of ordered algebraic systems and through it nonstandard analysis, the theory of the rates of growth of functions and the theory of surreal numbers.
For our purpose here, however, we limit our attention to its two continuity conditions, each of which, unlike the
Dedekind continuity condition, is satisfiable by Archimedean as well as
non-Archimedean ordered abelian groups and ordered fields. In fact, as Veronese
 (1889, 1891) and Levi-Civita (1898) collectively demonstrate, each of the conditions is equivalent to the Dedekind continuity 
condition if and only if the Archimedean axiom is assumed.
\end{sloppypar} 

Veronese formulates his two continuity conditions as follows.

\smallskip
\emph{Relative Continuity Condition}. Every segment $XX'$ whose ends vary in opposite directions and becomes \emph{indefinitely small} contains an element outside the domains of variability of its ends (Veronese 1891, p. 128).
\smallskip

\emph{Absolute Continuity Condition}. Every segment $XX'$ whose ends vary in opposite directions and becomes \emph{indefinitely small in the absolute sense}  contains an element outside the domains of variability of its ends (Veronese 1891, p. 150).
\smallskip

Veronese unpacks his continuity conditions in terms of a variable segment $XX'$ that is the difference $AX'-AX$ of a pair of subsegments of a segment $AB$,  where $AX$ is a proper subsegment of $AX'$, henceforth written $AX \prec^{*} AX'$. While keeping $A$ fixed and preserving the condition that $AX \prec^{*}AX'$, $X$ is envisioned to increase in a strict monotonic fashion (without a greatest member) as $X'$ decreases in a strict monotonic fashion (without a least member), $$A..........\overrightarrow{X}..........\overleftarrow{X'}..........B $$ subject to the following conditions that flesh out Veronese's conceptions of becoming ``indefinitely small" and ``indefinitely small in the absolute sense", respectively. 

\smallskip

\emph{Indefinitely Small}. For each segment $s$ \emph{that is finite relative to an arbitrarily given unit segment which could be taken to be $AB$},  $X$ and $X'$ take on values $X_s$ and $X'_s$, respectively, where $X_sX'_s  \prec s$.

\smallskip

\emph{Indefinitely Small in the Absolute Sense}. For each segment $s$,  $X$ and $X'$ take on values $X_s$ and $X'_s$, respectively, where $X_sX'_s \prec s$.

\smallskip
\noindent
Given the satisfaction of these respective conditions, Veronese's continuity conditions assert that there is a segment $AY$ such that $$AX_s \prec^{*} AY \prec^{*} AX'_s$$ for all such values $X_s$ of $X$ and all such values $X'_s$ of $X'.$

Veronese refers to his first continuity condition as a relative continuity condition since it is concerned with families of segments that grow arbitrarily small subject to the proviso that they remain finite relative to a given unit segment. The absolute continuity condition, by contrast, is concerned with families of segments that grow arbitrarily small subject to the limits of the geometric space itself. 

Veronese's relative continuity condition ensures that if one limits oneself to the segments that are finite relative to an arbitrarily selected segment $s$ and if one collects together into equivalence classes all such segments that differ from one another by amounts that are infinitesimal relative to $s$, the resulting system of equivalence classes with order defined in the expected manner is isomorphic to the standard continuum. Moreover, if (as in Veronese's geometry) the system of directed segments on a line emanating from a point is an ordered field, then if one takes the equivalence class containing $s$ as the unit and defines addition and multiplication of the equivalence classes in the familiar geometrical fashion, the resulting system is isomorphic to the positive cone $\mathbb{R}^+$ of the ordered field of real numbers. That is, \emph{Veronese's non-Archimedean continuum is indistinguishable from the classical continuum when infinite and infinitesimal differences are ignored}.\footnote{In addition to a non-Archimedean Euclidean space, which is what we are considering above, Veronese considers a non-Archimedean elliptic space, where the the system of directed segments emanating from a point has a somewhat different structure. However, in the elliptic case, the import of his relativity continuity condition is much the same.} 

\sloppypar
Unlike the segment $AY$ in the relative continuity condition, which is unique if and only if the Archimedean axiom holds, the segment $AY$ in the absolute continuity hypothesis is invariably unique, as Veronese was well aware. For this reason, Veronese's absolute continuity condition may also be stated in the following algebraic form that more clearly highlights its relation to the continuity condition of Dedekind.

\smallskip
\emph{Absolute Continuity: Algebraic Formulation}. Let $G$ be an ordered abelian group (or the positive cone thereof). If $(A, B)$ is a Dedekind cut of $G$ and if for each positive $\epsilon \in G$ there are elements $a$ of $A$ and $b$ of $B$ for which $b-a < \epsilon$, then either $A$ has a greatest member or $B$ has a least member, but not both.
\smallskip

It is a simple matter to show that in the Archimedean case, and only in
the Archimedean case, Veronese's metrical condition on cuts is invariably satisfied. Thus, for Veronese, unlike
for Dedekind, continuous systems of magnitudes need not be completely
devoid of Dedekind gaps, though they must be devoid of those Dedekind gaps that satisfy
the metrical condition satisfied in the classical case. Veronese maintained that insofar as the intuitive conception of a continuum does not require the Archimedean axiom, it is his absolute continuity condition, rather than Dedekind's continuity condition, that is intuitively more justifiable.

The algebraic and geometric formulations of Veronese's absolute continuity condition were widely discussed during the first decade of the twentieth century by authors such as H\"older, A. Schoenflies, L. E. J. Brouwer, K. Vahlen, G. Vitali, F. Enriques and Hahn\footnote{For references, see (Ehrlich 2006, pp. 66-70).} before sinking into relative obscurity as the Cantor-Dedekind conception of the continuum solidified its status as the standard conception. However, it was resurrected (without reference to Veronese) by a number of authors including R. Baer (1929, 1970),  L. W. Cohen and C. Goffman (1949), K. Hauschild (1966)  and Dana Scott (1969a), who (along with others) carefully studied it as a \emph{completeness} condition. Among the results that emerged from these investigations is that every nontrivial densely ordered abelian group (ordered field) admits an extension, unique up to isomorphsim, to a least ordered abelian group (ordered field) that satisfies Veronese's absolute continuity condition. Motivated by the work of Cohen and Goffman and especially  Scott, and unaware of the Veronesean roots of the condition, Zakon (1969, p. 226) asked if nonstandard models of analysis are complete in the absolute sense of Veronese. In response it was found that some are (e.g. Keisler 1974; Keisler and Schmerl 1991; Jin and Keisler 1993) and others are not (e.g. Kamo 1981, 1981a;  Keisler and Schmerl 1991; Ozawa 1995), something we will return to in \S 8. Those that are are usually said to be \emph{Scott Complete}.\footnote{Veronese's absolute continuity condition has emerged as a standard concept in the theory of ordered algebraic systems, albeit usually under a variety of other names. For further references, see (Ehrlich 1997, p. 224).} 

On the other hand, since (as can be readily shown) an ordered field is continuous in the relative sense of Veronese if and only if contains an isomorphic copy of $\mathbb{R}$, every nonstandard model of analysis satisfies Veronese's relative continuity condition. The same is also true of the ordered fields due to Hahn discussed in the following section, the first group of which drew inspiration from the Cantor-Dedekind continuum as well as the non-Archimedean continua of Veronese.

\section{hahn's non-archimedean generalizations of the archimedean arithmetic continuum}

Though Veronese's construction of his non-Archimedean
geometric continuum is synthetic, he represents the line segments
that emerge from his construction using a loosely defined, complicated system of
numbers consisting of finite and transfinite series of the form
\[
\infty_1^{y_1}r_1+\infty_1^{y_2}r_2+\infty_1^{y_3}r_3+\dotsb
\]
where $r_1,r_2,r_3,\dots$ are real numbers, and
$\infty_1^{y_1},\infty_1^{y_2},\infty_1^{y_3},\dots$ is a sequence of
\textit{units}, each of
which is infinitesimal relative to the preceding units, $\infty_1$ being
the number (of ``infinite order~1'' (1891, p. 101)) introduced by Veronese to represent the infinitely large line
segment whose existence is postulated by his ``hypothesis on the existence
of bounded infinitely large segments'' (1891, p. 84). Veronese's number system was provided an 
analytic foundation by Levi-Civita (1892--1893, 1898), who therewith
provided the first analytic constructions of non-Archimedean ordered fields. Building on the work of Levi-Civita, Hahn (1907) constructed 
non-Archimedean ordered fields (and ordered abelian groups more generally) having properties that generalize the familiar
continuity properties of Dedekind and Hilbert (Ehrlich 1995, 1997, 1997a), and
he demonstrated (vis-\`{a}-vis his celebrated \emph{Hahn Embedding Theorem} for ordered abelian 
groups\footnote{For a historical account of Hahn's embedding theorem and his momentous contribution to non-Archimedean mathematics more generally, see (Ehrlich 1995) which includes an extensive discussion with references of all the material in this section.})
that his number systems provide a panorama of the finite, infinite and infinitesimal numbers that
can enter into a non-Archimedean theory of continua based on the concept of an ordered
field (Ehrlich 1995, 1997, 1997a). This idea was later brought into sharper focus when it
was demonstrated that every ordered field can be embedded in a perspicuous fashion in a suitable \emph{Hahn field},\footnote{This theorem, which extends Hahn's embedding theorem for ordered abelian groups to ordered fields, has a complicated history that makes it difficult to attribute it to any single author. However, by 
the early 1950s, as a result of the work of Kaplansky (1942), it appears to have assumed 
the status of a ``folk theorem'' among knowledgeable field theorists, with numerous proofs published thereafter. For references and details, see (Ehrlich 1995).} the general construction of the latter of which is given as follows.

\begin{Hahn Field}[Hahn 1907]
Let $\mathbb{R}$ be the ordered field of real numbers and $G$ be a nontrivial
ordered abelian group. The collection, $$\mathbb{R}((t^{G}))$$ of all series
\[
\sum_{\alpha<\beta} t^{g_{\alpha}}\cdot r_{\alpha},
\] 
where $\beta$ is an ordinal, $\{g_{\alpha}\colon \alpha<\beta\}$ is a (possibly empty) strictly decreasing sequence of members of
$G$ and $\{r_{\alpha}\colon \alpha<\beta\}$ is a sequence of members
of $\mathbb{R}-\left \{ 0 \right \}$ is a non-Archimedean ordered
field when the order is defined lexicographically and sums and products are defined 
termwise, it being understand that $t^{\gamma}\cdot
t^{\kappa}=t^{\gamma+\kappa}$. 
\end{Hahn Field}

Hahn further established 
\begin{Hahn Field}[Hahn 1907]
If in the above construction one restricts the members of the universe to series indexed over all ordinals having cardinality less that an uncountable cardinal $\aleph_\alpha$, the resulting structure, which we will denote $$\mathbb{R}((t^{G}))_\alpha$$ is likewise an ordered field.
\end{Hahn Field}

While Hahn assumed $G$ to be a set, Hahn Field 2 continues to hold in NBG when $G$ is a proper class and the series are indexed over all ordinals $\alpha<\beta$, where $\beta$ is an arbitrary ordinal. Henceforth, we will assume that Hahn Field 2 is so extended to include such cases, the latter of which will be denoted  $$\mathbb{R}((t^{G}))_{On}.$$

\smallskip
The theories of Hahn groups and Hahn fields play critical roles in the theory of ordered algebraic systems and in the corresponding model theory thereof. In \S 9 we shall see how particular $\aleph_\alpha$\emph{-restricted Hahn fields} from Hahn Field 2 shed light on some of the most distinguished non-Archimedean continua of our day. In the remainder of this section we will introduce a number of related concepts that will be employed in the subsequent discussion and which shed further light on Hahn's celebrated constructs.

Using absolute values, Veronese's comparative notions of segments that are finite, infinite, and infinitesimal relative to one another were extended to the members of nontrivial ordered abelian groups by Levi-Civita (1898), and through the work of Levi-Civita and Hahn have become fixtures in post-nineteenth century mathematics, albeit under a variety of different rubrics. For example, following B. H. Neumann (1949, p. 205), instead of saying that the nonzero elements are finite relative to one another, they are now commonly said to be \emph{Archimedean equivalent to} or \emph{Archimedean relative to} to one another. Archimedean equivalence partitions the nonzero elements of an (additively written) ordered abelian group into disjoint
classes called \emph{Archimedean classes}.

Formally speaking, if $a$ and $b$ are nonzero members of an ordered abelian group $G$, then  $a$ is said to be Archimedean equivalent  to $b$   if there are positive integers $m$  and $n$   such that  $m|b|>a$  and $n|a|>b$; if $a$ and $b$ are not Archimedean equivalent, then $a$  is said to be \emph{infinitesimal (in absolute value) relative to}  $b$ and $b$  is said to be \emph{infinite (in absolute value) relative to} $a$, if $|a|<|b|$. In accordance with these conventions, $0$ is infinitesimal (in absolute value) relative to every other member of  $G$. Moreover, if $G$ is the additive group of an ordered field or of an ordered ring with a unit or identity more generally, the elements are simply said to be \emph{infinite (in absolute value)} and \emph{infinitesimal (in absolute value)}, respectively, if they are infinite (in absolute value) relative to and infinitesimal (in absolute value) relative to the identity. 

In virtue of the lexicographical ordering, every nonzero member of a Hahn field $\mathbb{R}((t^{G}))$ is Archimedean equivalent to exactly one element of the form $t^g$ ($g\in G$), the latter of which may be regarded as a canonical ``unit" element of the Archimedean class. Two such elements $x$ and $y$ are Archimedean equivalent if and only if their zeroth exponents ($g_0$ in the statement of Hahn Field 1) are equal, and $x$ is infinitesimal (in absolute value) relative to $y$ if and only if the zeroth exponent of $x$ is less than the zeroth exponent of $y$. Thus, like the numbers representing Veronese's non-Archimedean continuum, Hahn's numbers are formal sums of terms, each being infinitesimal (in absolute value) relative to the preceding terms, where each term is a nonzero real multiple of the canonical unit of its Archimedean class. Furthermore, every nonzero $x \in \mathbb{R}((t^{G}))$ is the sum of three components: the \emph{purely infinite part} of $x$, whose terms have positive exponents; the \emph{real part} of $x$, whose sole term has exponent $0$; and the \emph{infinitesimal part} of $x$, whose terms have negative exponents. The appellation ``real part" is motivated by the fact that $\{rt^0: r \in \mathbb{R}\}$ is a canonical copy of the ordered field of reals in $\mathbb{R}((t^{G}))$.

\section{the pantachies of du bois-reymond and hausdorff }

Although interest in the rates of growth of real functions is already found
in Euler's {\itshape De infinities infinitis gradibus tam infinite magnorum
quam infinite parvorum} (On the infinite degrees
of infinity of the infinitely large and infinitely small) (1778), their systematic study was
first undertaken by Paul du~Bois-Reymond  (cf. 1870-71, 1875, 1877, 1882), under the rubric \textit{Infinit\"{a}rcalc\"{u}l} (infinitary calculus).\footnote{For a complete list of
du~Bois-Reymond's writings on his \textit{Infinit\"{a}rcalc\"{u}l} and a 
survey of the contents thereof, see (Fisher 1981).} And while du-Bois-Reymond never attempted to employ this work to develop a non-Archimedean theory of continua, he and others including Poincar\'e (1893/1952, pp. 28-29) believed it provided an intimation of the possibility of such a theory. Moreover, as we shall see in \S 8, there are intriguing historical and conceptual relations between his theory and one of the most important such theories of our day. 

Though du Bois-Reymond's contribution to the \textit{Infinit\"{a}rcalc\"{u}l} is concerned solely with functions and is analytic by nature, it is intimately related to his ideas on quantity, in general, and the geometric linear continuum, in particular. Unlike Veronese, who admitted the possibility of an Archimedean geometric linear continuum, for du Bois-Reymond the geometric linear continuum is necessarily non-Archimedean. Indeed, on the basis of a misguided argument that has been aptly characterized as ``breathtaking"  (Fisher 1981 p. 114), du Bois-Reymond maintained that the infinite divisibility of the line implies that ``the unit segment decomposes into infinitely many subsegments of which none is finite" (1882, p. 72). Moreover, as in Veronese's non-Archimedean continuum, for every segment of du Bois-Reymond's continuum, there are segments that are infinitesimal relative to it. However, unlike Veronese's non-Archimedean continuum, du Bois-Reymond's continuum of segments is not the positive cone of an ordered field (or even of an ordered abelian group) since, according to du Bois-Reymond, ``two finite segments are equal when there is no finite
difference between them. [That is,] a finite quantity does not change if an infinitely
small quantity is added to it or taken away from it" (1882, pp. 73-74), even when the infinitely
small quantity is nondegenerate.

It was the comparison of quantities (of the same type) having different orders of magnitude that du Bois-Reymond took to be the primary object of his infinitary calculus (1882, pp. 75-66), but as was mentioned above, he only developed the theory for functions. In particular, du Bois-Reymond erects his \textit{Infinit\"{a}rcalc\"{u}l} primarily on families of increasing
functions from $\mathbb{R}^+=\{x\in \mathbb{R}\colon x>0\}$ to $\mathbb{R}^+$ such that for
each function $f$ of a given family, 
$\lim\limits_{x\to \infty} f(x)=+\infty$, and for each pair of functions $f$ and
$g$ of the family,
$0\leq \lim\limits_{x\to \infty}f(x)/g(x)\leq +\infty$. He assigns to each such
function $f$ a so-called \textit{infinity}, and
defines an ordering on the infinities of such functions by stipulating that for each pair of
such functions $f$ and~$g$:
\begin{align*}
&\text{$f(x)$ has an infinity \textit{greater than} that of $g(x)$, if
$\lim\limits_{x\to
\infty} f(x)/g(x)=\infty$;}\\
&\text{$f(x)$ has an infinity \textit{equal to} that of $g(x)$, if
$\lim\limits_{x\to \infty}f(x)/g(x)=a\in \mathbb{R}^+$;}\\
&\text{$f(x)$ has an infinity \textit{less than} that of $g(x)$, if
$\lim\limits_{x\to \infty}f(x)/g(x)=0$.} 
\end{align*}

\begin{figure}[h!]
\centering\includegraphics[width=0.8\textwidth]{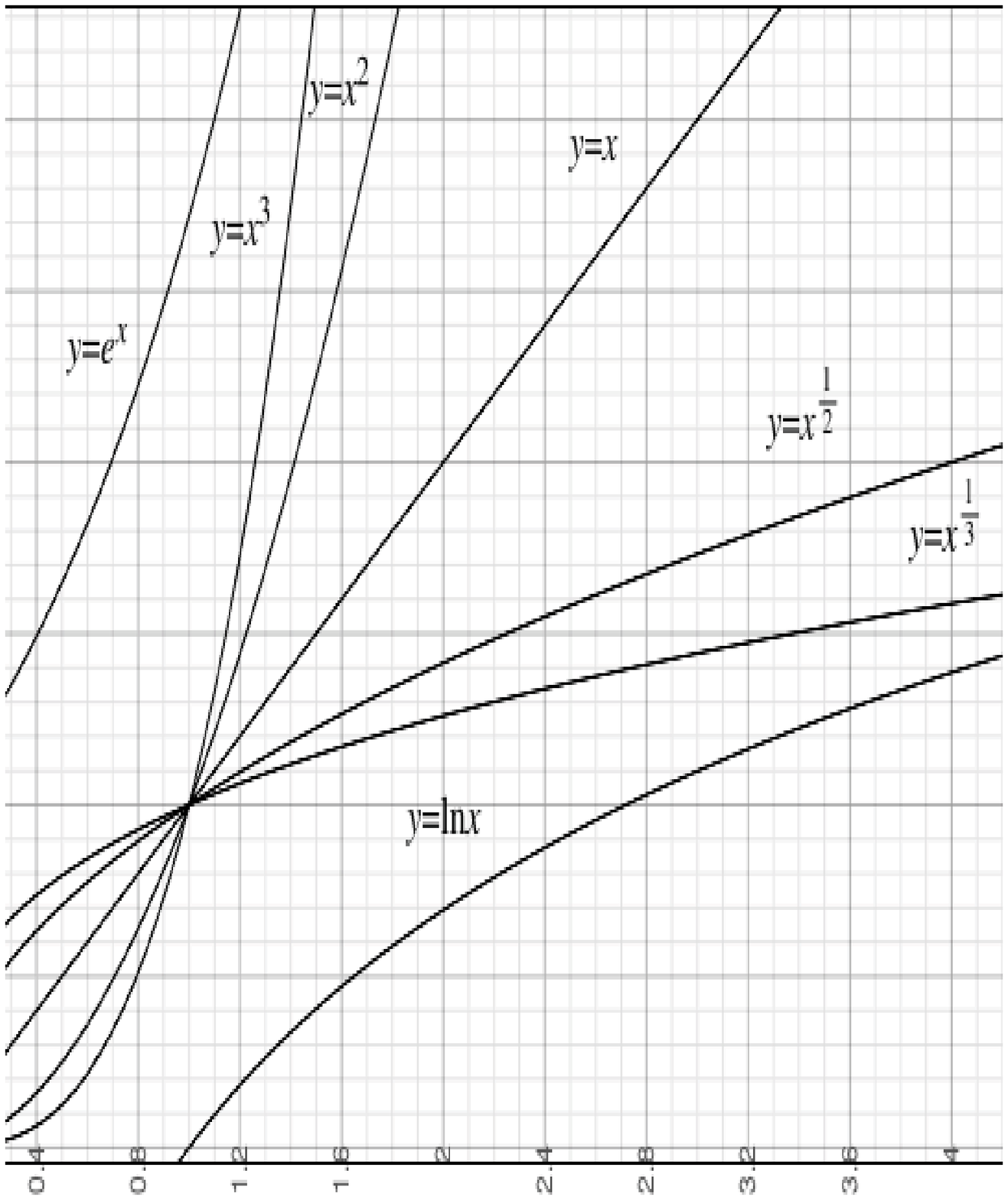}
\vspace{-1.1 in}
\caption{\ }\label{Ehrlichfig2}
\end{figure}

In accordance with this scheme, the infinities of the following functions \[\dots,\ln(\ln x),\ln x,\dots, x^{1/n},\dots,x^{1/3},x^{1/2},x, x^2,x^3\dots,x^n,\dots,e^x,e^{e^x},
\dots
\] increase as we move from left to right. Moreover, as the comparative graphs of several of
these functions illustrate (see Figure~\ref{Ehrlichfig2}), given any two
functions $f$ and $g$ having different infinities from a family of
the just-said kind, $f(x)$ has a greater infinity than $g(x)$ if and only if
$f(x)>g(x)$ for all $x>{}$some $x_0$. Furthermore, since, for example, $x^2$ has a greater infinity than $x$ and $$\lim\limits_{x\to \infty}(x^2 + x)/x^2 =1,$$ $x^2 + x$ has the same infinity as $x^2$, which illustrates for functions du Bois-Reymond's corresponding idea for segments, that if a segment $s$ is infinitesimal relative to a segment $s'$, then the segment resulting from adding $s$ to $s'$ is equal in length to $s'$.\footnote{Stolz, who was a student of du Bois-Reymond (as well as Weierstrass), developed another number system (Stolz 1884) that modeled this absorptive aspect of du Bois-Reymond's conception. For a discussion of this little-known system of ``moments", see (Ehrlich 2006, \S 3).}  

While du~Bois-Reymond usually restricted his investigations to families of functions of the
above said kind, on occasion he mistakenly assumed that each pair of increasing
functions $f$ and $g$ from $\mathbb{R}^+$ to $\mathbb{R}^+$ for which $\lim\limits_{x\to
\infty}f(x)=+\infty$ and $\lim\limits_{x\to \infty}g(x)=+\infty$ could be
compared in the manner described above. This led him to postulate the
existence of an all-inclusive ordering of the infinities of such real
functions---an \textit{infinitary pantachie}, as he called 
it (1882, p. 220).\footnote{Du Bois-Reymond explains that his adjective ``pantachie'' derives from the Greek 
words for ``everywhere''.}
Such a pantachie, according to du~Bois-Reymond, 
would provide a conception of a numerical linear continuum ``denser" than that of
Cantor and Dedekind. 

Georg Cantor, however, who was an ardent opponent of infinitesimals and non-Cantoriean infinities (cf. Ehrlich 2006) would have no part of this. Indeed, having demonstrated as Stolz (1879, p. 232) and Pincherle (1884, p. 742)
had before him that there could be no such all-inclusive ordering of the infinities of
such functions, Cantor proclaimed: ``the `\textit{infinitary pantachie},' of 
du~Bois-Reymond, belongs in the wastebasket \textit{as nothing but paper
numbers!}'' (1895, p. 107).
Felix Hausdorff, on the other hand, suggested, ``[t]here is no reason to reject the entire theory because
of the possibility of incomparable functions as G. Cantor has done'' (1907, p. 107), and in
its place undertook the study of {\itshape maximally inclusive sets of
pairwise comparable real functions}, each of which, retaining du~Bois-Reymond's term, he calls an \textit{infinitary pantachie} or a
\textit{pantachie} for short. This led him to his well-known investigation
of $\eta_{\alpha}$-\emph{orderings} (1907), and to the following less 
well-known theorem.

\begin{Pantachie}[Hausdorff 1907, 1909]
Infinitary pantachies exist. If $P$ is an
infinitary pantachie, then $P$ is an $\eta_1$-ordering of power
$2^{\aleph_0}$; in fact, $P$ is (up to
isomorphism) the unique $\eta_1$-ordering of power $\aleph_1$, assuming
(the Continuum Hypothesis)~$\text{CH}$.
\end{Pantachie}

An ordered set $L$ is said to be an $\eta_{\alpha}$-ordering, if for all subsets $A$ and $B$ of $L$ of power $<\aleph_\alpha$, where $A<B$ (i.e. every member of $A$ precedes every member of $B$), there is a $y \in L$ such that $A <\left \{ y \right \} <B$. The ordered field $\mathbb{R}$ of real numbers is an $\eta_0$-ordering, though not an $\eta_1$-ordering, and as such Hausdorff's pantachies are in a precise sense more dense than $\mathbb{R}$, thereby lending precision to du Bois-Reymond's intimation. 

As was noted above, du Bois-Reymond further believed that pantachies exhibit numerical aspects, but he never undertook an arithmetization of them. Motivated by different considerations, however, Hausdorff did. Moreover, like his theory of pantachies, more generally, he did so in the following broader setting.

In addition to modifying du~Bois-Reymond's conception of
an infinitary pantachie, Hausdorff redirected du~Bois-Reymond's investigation
by investigating \textit{numerical sequences} rather than continuous
functions, though he showed that the desired results about the latter can be obtained as corollaries from results about the former. He also deleted the monotonicity assumption and
replaced the infinitary rank ordering with the \textit{final rank
ordering} (which was illustrated above). That is, Hausdorff redirected du~Bois-Reymond's investigation to the study of subsets of the set $\mathfrak{B}$ of all numerical sequences $A=(a_1,a_2,a_3,\dots,a_n,\dots)$
in which the $a_n$ are real numbers, and he defines the ``final
ordering'' on $\mathfrak{B}$ (and subsets
 thereof) by the conditions $A<B$ if \textit{eventually} $a_n<b_n$, $A=B$
 if \textit{eventually} $a_n=b_n$, $A>B$
if \textit{eventually} $a_n>b_n$, and $A\parallel B$ (i.e., $A$ is
\textit{incomparable} with $B$) in all other cases,
where ``eventually'' means for all values of $n$ with the exception of a finite
number, thus
for all $n\geq{}$some~$n_0$ [1909 in (Plotkin 2005, p.~276)].\footnote{Strictly
  speaking, in his (1907), unlike his (1909), 
Hausdorff only considers numerical sequences in which the $a_n$ are positive real numbers. 
However, the proof of the above theorem carries over to the more general case. 
Moreover, in his (1907), unlike his later works, he uses the term ``finally'' instead of 
``eventually'' in the definitions of $<$, $>$, $=$ and $\parallel$.}
Hausdorff, who bases his theory on 
representative elements of the equivalences classes of eventually equal
numerical sequences rather than on the equivalence classes themselves,
calls a subset $\mathfrak{B}'$ of $\mathfrak{B}$ totally ordered by the final order a
\textit{pantachie} if it is not properly contained in another
subset of $\mathfrak{B}$ totally ordered by the final order.

In his aforementioned investigation of 1907, Hausdorff first raised the question of the existence of
a pantachie that is algebraically a field, but he only made partial headway in providing
an answer. However, in 1909 he returned to the problem and provided a stunning positive
answer. Indeed, beginning with the ordered set of numerical sequences of the form
$(r,r,r,\dots,r,\dots)$ where $r$ is a rational number, and utilizing what appears to be the very first
algebraic application of his maximal principle, Hausdorff proves the following 
little-known, remarkable result.

\begin{Pantachie}[Hausdorff 1909]
There is a pantachie of numerical sequences that is an ordered field, whose field operations are given by
\begin{align*}
A+B&=(a_1+b_1,a_2+b_2,\dots, a_n+b_n,\dots),\\
A-B&=(a_1-b_1,a_2-b_2,\dots, a_n-b_n,\dots),\\
AB&=(a_1b_1,a_2b_2,\dots,a_nb_n,\dots),\\
A/B&=(a_1/b_1,a_2/b_2,\dots, a_n/b_n,\dots),
\end{align*}
where $A+B$, $A-B$, $AB$ and $A/B$ are defined up to final equality. 
Any ordered field that is a pantachie is, in fact, real-closed.\footnote{This result had been all but forgotten until it was resurrected in (Plotkin 2005, pp. 271-301) and (Ehrlich 2012, p. 29).}
\end{Pantachie}

Writing before Artin and Schreier (1926), Hausdorff does not employ the term ``real-closed", though as Ehrlich (2012, p. 29) observed, Hausdorff proves enough to show that an ordered field that is a pantachie is real-closed. However, to fully appreciate the significance of Hausdorff's result as well as the extent of its intimate relation to certain contemporary non-Archimedean continua, we require some background about real-closed fields.\footnote{By a classical result of R. Baer (1927) and W. Krull (1932), the collection of Archimedean classes of an ordered field $A$ constitutes an ordered abelian group induced by the order and multiplication in $A$ (see, Ehrlich 1995, p. 186). Though we will not further pursue the matter here, we note that the system of orders of infinity associated with the members of a Hausdorff pantachie $\mathbb{H}_p$ is isomorphic to the positive cone of the ordered abelian group of Archimedean classes of $\mathbb{H}_p$, the latter being (up to isomorphism) the unique divisible ordered abelian group that is an $\eta_1$-ordering of power $\aleph_1$, assuming $\text{CH}$ (e.g., (Ehrlich 1988)).}

\section{elementary continua}	

In their groundbreaking work on ``real algebra" Emil Artin and Otto Schreier (1926) sought to characterize the algebraic content of a kind of ordered field of which the arithmetic continuum is paradigmatic. This led to their theory of \emph{real-closed} fields, which is a jewel of twentieth-century mathematics. Following Artin and Schreier, an ordered field $K$ may be said to be \emph{real-closed}
if it admits no extension to a more inclusive
ordered field that results from supplementing $K$ with
solutions to polynomial equations with coefficients in $K$. Intuitively speaking, real-closed ordered
fields are precisely those ordered fields having no ``holes"
that can be filled by algebraic means alone. Among a number of important properties of the arithmetic continuum they showed are satisfied by real-closed ordered fields more generally is the \emph{intermediate value theorem for polynomials of a single variable}. In fact, as soon became clear, for an ordered field $K$, the satisfaction of the intermediate value theorem for polynomials of a single variable over $K$ is equivalent to $K$ being real-closed (cf. Tarski 1939; Warner 1965, pp. 492-493). Accordingly, the idea of a real-closed ordered field is equivalent for the case of polynomials of a single variable to satisfying what has traditionally been regarded as one of the quintessential properties of a continuous function--the formalization of the intuitive idea that a continuous curve connecting points on the opposite sides of a straightline intersects the given line. A similar equivalence was later established for the \emph {extreme value theorem for polynomials in a single variable} (Gamboa 1987), another prototypical classical property of continuous functions. Shedding still more light on the relation between real-closed ordered fields and $\mathbb{R}$, Tarski (1948, 1959) demonstrated
that real-closed ordered fields are precisely the
ordered fields that are first-order indistinguishable from
$\mathbb{R}$, or, to put this another way, they are precisely the
ordered fields that satisfy the \emph{elementary} (first-order)
content of the Dedekind continuity axiom. For this reason they are sometimes called \emph{elementary
continua}. When Tarski's result is combined with a classical result of Artin and Schreier (1926) on real-closed ordered fields, one obtains the critical fact that every ordered field is contained in a (to within isomorphism) smallest elementary continuum, the latter having the same cardinality of the given field.

While $\mathbb{R}$ is the best known elementary continuum,
as is evident from the above it is not the only one.
Some elementary continua, like $\mathbb{R}$, are Archimedean,
though most are non-Archimedean, and among the latter
many are extensions of $\mathbb{R}$. In the following section, we will discuss a class of real-closed extensions of $\mathbb{R}$ that are among the foremost non-Archimedean continua of our day, and in the subsequent section we will draw attention to some of the interesting connections between them, Hausdorff's pantachies, the surreals and some of the elementary continua that emerge from the constructions of Hahn. Unlike Veronese's non-Archimedean continuum, which was motivated by generalizing the geometrical properties of the classical geometric continuum, those to which we now turn are motivated by the idea of providing a non-Archimedean treatment of the analytic properties of $\mathbb{R}$.

\section{nonstandard (robinsonian) continua}	

\begin{sloppypar}
In the early 1960s Abraham Robinson (1961, 1966) 
made the momentous discovery that among the real-closed extensions of the reals
 there are number systems that can provide the basis for a consistent and entirely
 satisfactory nonstandard approach to analysis based on infinitesimals, a possibility that had been called into question by many since the latter decades of the nineteenth century.\footnote{There were some earlier successes in this direction, including the modest contributions of Levi-Civita (1898),  Neder (1941, 1943) and Gesztelyi (1958). The most successful of these is the theory of Schmieden and Laugwitz (1958; also see, Laugwitz 1983, 1992, 2001). Unlike Robinson, Schmieden and Laugwitz make use of a partially-ordered number system containing zero-divisors, with the result that many of the classical results need to be reformulated. As Laugwitz aptly acknowledged, Robinson's system ``was much more powerful when it came to applications in contemporary research" (2001, p. 128).} By analogy with
 Thoralf Skolem's (1934) \textit{nonstandard model of arithmetic}, a number system from which
 Robinson drew inspiration, Robinson called his number systems
 \textit{nonstandard models of analysis}. These number systems, which are now more often called
\textit{hyperreal number systems} (Keisler 1976, 1994), may be characterized as follows: let
$\langle\mathbb{R},S\colon S\in \frak F\rangle$ be a relational structure where
$\frak F$ is the set of all finitary relations defined
on $\mathbb{R}$ (including all functions). Furthermore, let $\mathbb{^*R}$ be a proper extension of $\mathbb{R}$ and for
each $n$-ary relation $S\in \frak F$ let ${^*\!S}$ be an $n$-ary relation on
$\mathbb{^*R}$ that is an extension of $S$. The
structure $\langle\mathbb{^*R},\mathbb{R},{^*\!S}\colon S\in\frak F\rangle$ is said to be a hyperreal number system if it satisfies the
\textit{transfer principle}: every $n$-tuple of real numbers satisfies the same first-order formulas in
$\langle\mathbb{R},S\colon S\in \frak F\rangle$ as it satisfies in
$\langle\mathbb{^*R},\mathbb{R},{^*\!S}\colon S\in\frak F\rangle$.
\end{sloppypar}

The existence of hyperreal number systems is a consequence of the
compactness
theorem of first-order logic and there are a number of algebraic techniques that can be
employed to construct such systems. The earliest and still one of the most commonly employed such techniques is the
\emph{ultrapower} construction (e.g. Keisler 1976, pp. 48-57; Goldblatt 1998, ch. 3), a construction that was introduced in full generality by \L o\'s (1955), identified as a source of hyperreal number systems by Robinson (1961, p. 3) and popularized as such by Luxemburg (1962) and Stroyan and Luxemburg (1976).\footnote{Hewitt (1948) had already employed the ultrapower construction to obtain a non-Archimedean real-closed extension of the reals, but unlike \L o\'s he did not establish the transfer principle that accrues from the construction, the latter being crucial to hyperreal number systems in the sense used above.} Not all hyperreal number systems can be obtained this way, however. By results of H. J. Keisler (1963,
1976, pp. 58-59), on the other hand, every hyperreal number system is isomorphic to
a \emph{limit ultrapower}.

Using the transfer principle, one can develop satisfactory nonstandard conceptions and treatments of \emph{all} of the concepts and theorems of the calculus (e.g. Keisler 1976; Goldblatt 1998; Loeb 2000). For example, it follows from the transfer principle that:

\begin{quote}
 A real-valued function $f$  is continuous at $a \in  \mathbb{R}$ (in the standard sense) if and only if $^*\!f(x)$ is infinitely close to  $^*\!f(a)$ whenever $x$  is infinitely close to $a$, for all  $x \in  \mathbb{^*R};$ 
 
 \end{quote}

 \noindent
and on the basis of this one may prove the familiar classical results concerning the continuity of real-valued functions including the intermediate and extreme value theorems (e.g. Goldblatt 1998, pp. 79-80). 

Of course modern analysis goes well beyond the traditional province of the calculus, dealing with arbitrary sets of reals, sets of sets of reals, sets of functions from sets of reals to sets of reals, and the like. For example, it is commonplace for analysts to prove theorems about the set of all continuous functions on the reals or about the set of all open subsets of the reals, sets to which the just-said transfer principle does not apply. To obtain nonstandard treatments of these aspects of analysis a more general transfer principle is required. For this purpose Robinson (1966) originally employed a type-theoretical version of higher order logic, but it proved to be unpopular. Since then, following Robinson and Zakon (1969), it has become most common to employ a transfer principle associated with a structure $\langle V(\mathbb{^*R}),V(\mathbb{R}),*\rangle$ that generalizes the corresponding hyperreal number system $\langle\mathbb{^*R},\mathbb{R},{^*\!S}\colon S\in\frak F\rangle$. In this setup the transfer principle emerges from an elementary embedding $*$ that relates the members of the superstructure $V(\mathbb{R})$ over $\mathbb{R}$ with those of the superstructure $V(\mathbb{^*R})$ over $\mathbb{^*R}$, where for any set of individuals $X$, the \emph{superstructure $V(X)$ over $X$} is defined by $V(X) = \bigcup_{n<\omega}V_n(X)$, where $V_0(X) = X$ and $V_{n+1}(X) = V_{n}(X) \cup \mathcal{P}(V_{n}(X))$, $\mathcal{P}$ being the power set operator.\footnote{While the superstructure approach is the most popular approach, there are a variety of alternative approaches to nonstandard analysis, each with its own virtues, that cannot be touched upon here. For references, overviews and insightful discussions of a number of such approaches including the nonstandard set-theoretic approaches of Nelson, Hrbacek, and Kanovei and Reeken, see (Kanovei and Reeken 2004; Benci \emph{et al} 2006; Fletcher \emph{et al} 2017). Also see (Sanders 2018) for constructivist approaches to, and aspects of, nonstandard analysis.} More specifically, the \emph{generalized transfer principle} asserts that any property expressible in the language $\{=,\in\}$ with only bounded quantifier formulas holds
in $V(\mathbb{R})$ for an entity $S$ if and only if it holds in $V(\mathbb{^*R})$ for the corresponding entity $^*\!S$. 
 
Robinson and Zakon's superstructure paper was presented at a conference on applications of model theory to algebra, analysis and probability held at Cal Tech in May of 1967.  Luxemburg, who organized the conference, also presented a paper (1969) that has had a profound impact on the development of nonstandard analysis. In his paper, Luxemburg incorporated the conception of a $\kappa$\emph{-saturated model} (developed by Morley and Vaught 1962 and Keisler 1964: see Chang and Keisler 1990) into the subject and applied it to the theory of topological spaces. Soon thereafter, Peter Loeb (1975; also see, Hurd and Loeb 1985, ch. IV), making critical use of Luxemburg's ideas on saturation, introduced the \emph{Loeb measure} and therewith (along with Anderson (1976), Henson (1979) and others) developed the nonstandard treatments of the classical theories of measure and integration due to Lebesgue. Since that time, the use of saturation assumptions in nonstandard analysis has become standard. 

As the term ``saturation" suggests, models satisfying saturation assumptions are rich in elements. The saturation assumptions most commonly employed make use of the critical notion of an \emph{internal set} or an \emph{internal relation}, where an entity $b \in V(\mathbb{^*R})$ is said to be internal if $b \in {^*\!a}$ for some $a \in V(\mathbb{R})-\mathbb{R}$. Roughly speaking, the internal sets are those sets of hyperreals that inherit the first-order properties of their real counterparts.

As the following formulation makes clear, saturation assumptions come in varying strengths. 
\begin{Saturation 1 *}
$V(\mathbb{^*R})$ is said to be $\kappa$-\emph{saturated} if every subset $X$ of $V(\mathbb{^*R})$ consisting of $<\kappa$ internal sets has a nonempty intersection whenever every finite subset of $X$ has a nonempty intersection. If $\kappa$ is the cardinality of ${\mathbb{^*R}}$, $V(\mathbb{^*R})$ is said to be \emph{saturated}.
\end{Saturation 1 *}

For many purposes, $\omega_1$-\emph{saturation} is all the saturation of $V(\mathbb{^*R})$ a nonstandard analyst requires. However, as Henson and Keisler (1986) demonstrated, even a reliance on $\omega_{1}$-saturation is of metamathematical significance. In particular, if, as mathematical practice seems to suggest, standard analysis can be formalized in a conservative extension of second-order arithmetic, then ``nonstandard analysis (i.e. second order nonstandard arithmetic) with the $\omega_1$-saturation axiom scheme has the same strength as third order arithmetic" (1986, p. 377) and is therefore stronger than standard analysis as is typically practiced. Indeed, as Henson and Keisler put it:

\begin{quote}
This shows that in principle there are theorems which can be proved with nonstandard analysis but cannot be proved by the usual standard methods. The problem of finding a specific and mathematically natural example of such a theorem remains open. However, there are several results, particularly in probability theory, whose only known proofs are nonstandard arguments which depend on saturation principles; see, for example, the monograph (Keisler 1984). Experience suggests that it is easier to work with nonstandard objects at a lower level than with sets at a higher level. This underlies the success of nonstandard methods in discovering new results. To sum up, nonstandard analysis still takes place within ZFC, but in practice it uses a larger portion of full ZFC than is used in standard mathematical proofs. (1986, pp. 377-378)\footnote{The result of Henson and Keisler stands in contrast with those of Harvey Friedman and Kreisel (Kreisel 1969), who established conservation results when no saturation assumption is assumed.}
\end{quote}

While a theorem of nonstandard analysis expressible in standard mathematics that cannot be established by standard means has yet to be identified, the recognition of the value of nonstandard techniques for discovering new results, a value emphasized by Robinson (1974, p. IX) and G\"{o}del (Robinson 1974, p. X), and reiterated by Henson, Keisler and others, appears to be growing. Noteworthy in this regard is the high-profile championing in both words and deeds of the use of nonstandard techniques for discovering new results by Terence Tao (e.g. Tao  2008, 2013, 2014), one of the most celebrated mathematicians of our day.

	Unlike $\mathbb{R}$, the structures that may play the role of $\mathbb{^*R}$ in a hyperreal number system are far from being unique up to isomorphism. In addition to there being such systems satisfying  varying types and strengths of saturation conditions, in an attempt to find hyperreal number systems that look more like $\mathbb{R}$ some authors have sought and established the existence of hyperreal number systems that are continuous in the sense of Veronese (see \S 4) or that satisfy generalizations of the Bolzano-Weierstrass property (e.g. Keisler and Schmerl 1991), properties that are not universally satisfied by hyperreal number systems. From a purely mathematical point of view this absence of uniqueness causes no difficulty, and it has been argued that from the standpoint of varying applications can even be advantageous (e.g. Keisler 1994, p. 229; Jin 1992, 1997). On the other hand, if one takes $\mathbb{^*R}$ to be a model of \emph{the} continuous straight line of geometry--\emph{something practitioners of nonstandard analysis tend not to do}--the absence of uniqueness is a bit disconcerting.

	Still, as several nonstandard analysts including Robinson (1973, p. 130), Lindstr\o m (1988, p. 82) and Keisler (1994, p. 229) emphasize, even $\mathbb{R}$ is not as unique as one would like to think since its uniqueness up to isomorphism is relative to the underlying set theory. Indeed, as Dana Scott aptly observed:
	
\begin{quote}
Maybe we have to face the fact that there are many distinct theories of the continuum.
(Scott 1969, p. 87)
\end{quote}
	
\noindent	
In particular, by retaining the construction of $\mathbb{^*R}$ and supplementing the set theory with additional axioms, one can change the second-order theory of the real line. This led Keisler (1994) to suggest that not only is ZFC not the appropriate underlying set theory for the hyperreal number system but also that set theory might have developed differently had it been developed with the hyperreal numbers rather than the real numbers in mind. According to Keisler, an appropriate set theory Òshould have the power set operation to insure the unique existence of the real number system, and another operation which insures the unique existence of the pair consisting of the real and hyperreal number systemsÓ (1994, p. 230).

Motivated by the above, Keisler (1976, pp. 57-60) introduced the following saturation assumption for hyperreal number systems to secure categoricity. Unlike the saturation condition for $*$, the saturation assumption for $\mathbb{^*R}$ does not appeal to the notion of an internal set.  

\begin{Saturation 2 *}
Let $S$ be a set of equations and inequalities involving real functions, hyperreal constants, and variables, such that $S$ has
a smaller cardinality than $\mathbb{^*R}$. If every finite subset of $S$ has a hyperreal solution, then $S$ has a hyperreal solution. 
\end{Saturation 2 *}

Saturated hyperreal number systems cannot be proved to exist in ZFC. However, in virtue of classical results from the theory of saturated models there is (up to isomorphism) a unique saturated hyperreal number system of power $\kappa$ whenever either $\kappa$ is (strongly) inaccessible and uncountable or the generalized continuum hypothesis holds at $\kappa$ (i.e., $\kappa=\aleph_{\alpha +1}=2^{\aleph_\alpha}$ for some $\alpha$). So, for example, by supplementing ZFC with the assumption of the existence of an uncountable inaccessible cardinal, one can obtain uniqueness (up to isomorphism) by limiting attention to saturated hyperreal lines having the least such power (Keisler 1976, p. 60). Keisler (1976, p. 59) and others (e.g. Kanovei and Reeken 2004; Ehrlich 2012; Borovik, Jin and Katz 2012) have also noted that it is possible to obtain uniqueness (up to isomorphism) by employing a saturated model that is a proper class, something we will say more about in the following section. In addition, following Henson (1974) and Ross (1990), it is possible to introduce axioms which, while falling short of full saturation, imply varying degrees of saturation and guarantee uniqueness up to isomorphism in certain powers that provably exist in ZFC.\footnote{The aforementioned ultrapower construction of hyperreal number systems is another source of nonuniqueness since the construction depends on an arbitrary choice of a nonprincipal ultrafilter. Hyperreal number systems that arise via iterated ultrapower constructions are analogous sources of nonuniqueness. However, Kanovei and Shelah (2004) have shown that in certain important cases these sources of nonuniqueness are eliminable via the existence of definable hyperreal number systems.} 

In the Preface to the Second Edition of his monograph \emph{Non-Standard Analysis}, Robinson expressed the view that:

\begin{quote}
the application of non-standard analysis to a particular mathematical discipline is a matter of choice and that it is natural for the actual decision of an individual to depend on his early training. (Robinson 1974, ix)
\end{quote}

Following a talk given by Robinson at the Institute for Advanced Study in Princeton in March of 1973, Kurt G\"odel went further when he maintained: 
\begin{quote}
There are good reasons to believe that non-standard analysis, in some version or other, will be the analysis of the future. (Robinson 1974, p. x) 
\end{quote}
\noindent
Among the reasons offered by G\"odel was that ``non-standard analysis frequently simplifies substantially the proofs, not only of elementary theorems, but also of deep results" with the consequence that it can ``facilitate discovery" as witnessed by ``the proof of the of invariant subspaces for compact operators". (Robinson 1974, p. x)
\noindent

``An even more convincing reason", said G\"odel, is that:

\begin{quote}
Arithmetic starts with the integers and succeeds by successively enlarging the number system by rational and negative numbers, irrational numbers, etc. But the next quite natural step after the reals, namely the introduction of infinitesimals, has simply been omitted. (Robinson 1973, p. x)
 \end{quote}

Roughly two decades later, motivated by considerations of simplicity and facility of discovery, H. J. Keisler more guardedly observed: 

\begin{quote}
At the present time, the hyperreal number system is regarded as somewhat of a novelty. But because of its broad potential, it may eventually become part of the basic toolkit of mathematicians. This process will probably take a very long time, perhaps 50 to 100 years. (1994, p. 235) 
\end{quote}

As we mentioned above, nonstandard analysis has indeed made inroads among standard analysts and mathematicians more generally and has been employed to prove a wide range of new results. Given its growing visibility and the growing recognition of its heuristic virtues, there is some reason to believe that nonstandard analysis, together with its non-Archimedean continuum, may indeed become part of the basic toolkit of mathematicians  in the timeframe Keisler suggested.

\section{the absolute arithmetic continuum: conway's system of surreal numbers}	

In his monograph \emph{On numbers and Games} (1976), J. H. Conway introduced a real-closed ordered field of \emph{surreal numbers} containing the reals and the ordinals as well as a great many less familiar numbers including $ - \omega $, ${\omega  \mathord{\left/
 {\vphantom {\omega  2}} \right.
 \kern-\nulldelimiterspace} 2}$, ${1 \mathord{\left/
 {\vphantom {1 \omega }} \right.
 \kern-\nulldelimiterspace} \omega }$, $\sqrt \omega  $, ${e^\omega }$, $\log \omega $ and $\sin \left( {{1 \mathord{\left/
 {\vphantom {1 \omega }} \right.
 \kern-\nulldelimiterspace} \omega }} \right)$ to name only a few, where $\omega $ is the least infinite ordinal. This particular real-closed field, which Conway calls ${\mathbf{No}}$, is so remarkably inclusive that, subject to the proviso that numbers--construed here as members of ordered fields--be individually definable in terms of \emph{sets} of NBG,  it may be said to contain, ``All Numbers Great and Small" (Conway 1976, p. 3). In this regard, ${\mathbf{No}}$ bears much the same relation to ordered fields that ${\mathbb{R}}$ bears to Archimedean ordered fields. 
 
The relation between the inclusiveness of ${\mathbb{R}}$ and ${\mathbf{No}}$ may be made precise by the following result where an ordered field (Archimedean ordered field) $A$ is 
said to be \textit{homogeneous universal} if it is
\textit{universal}--every ordered field (Archimedean ordered field)
 whose universe is a set or a proper class of NBG can be embedded in $A$--and it is 
\textit{homogeneous}--every isomorphism between subfields of $A$ whose universes are sets can 
be extended to an automorphism of $A$.\footnote{The just-said concepts are straightforward adaptations for Archimedean structures and proper classes  of conceptions from the classical theory of \emph{homogeneous universal structures} due to J\'onsson (1960) and Morley and Vaught (1962).}

 \begin{Surreal}[Ehrlich 1988, 1989, 1992] Whereas  $\mathbb{R}$ is (up to isomorphism) the unique homogeneous universal Archimedean ordered field, $\mathbf{No}$ is (up to isomorphism) the unique homogeneous universal ordered field.
  \end{Surreal}

In the case of Archimedean ordered fields, universality and homogeneous universality are equivalent since there is one and only one embedding of an Archimedean ordered field into $\mathbb{R}$. In the non-Archimedean case, on the other hand, this is not so since there are ordered fields that are universal, but not homogeneous. The homogeneity of \textbf{No} ensures that whenever subfields of \textbf{No}, whose universes are sets, are structurally indistinguishable when considered in isolation, they along with their surroundings in \textbf{No} are structurally indistinguishable as well.

	Since there is a multitude of real-closed ordered fields, it is natural to inquire if, like  $\mathbb{R}$, it is possible to distinguish $\mathbf{No}$ (to within isomorphism) from the remaining real-closed ordered fields by appealing solely to its order. The following definition, which is a straightforward extension for proper classes of Hausdorff's conception of an $\eta_\alpha$-ordering (see \S6), enables one to do just this. 
	
 \begin{Definition *}[Ehrlich 1987, 1988] 
 An ordered class $\langle A,<\rangle$ is said to be an \emph{absolute linear continuum} if for all subsets $L$ and $R$ of $A$  where $L<R$  there is a $y \in A$ such that $L<\{y\}<R$.
\end{Definition *}

	An absolute linear continuum $\langle A,<\rangle$ is both \emph{absolutely dense} in the sense that for each pair of nonempty subsets $L$ and $R$ of $A$ where $L<R$, there is a $y \in A$ such that $L<\{y\}<R$, and \emph{absolutely extensive} in the sense that given any (possibly empty) subset $X$ of $A$ there are members $a$ and $b$  of $A$ that are respectively smaller than and greater than every member of $X$. In fact, an ordered class is an absolute linear continuum just in case it has both of the just-stated properties. Accordingly, since every element of an ordered class must either lie between two of its nonempty subclasses or lie to the right or to the left of some (possibly empty) subclass, absolute linear continua are ordered classes having no order-theoretic limitations that are definable in terms of sets of standard set theory.

The following is the analog for absolute linear continua of Cantor's (1895) classical order-theoretic characterization of $\mathbb{R}$.

 \begin{Surreal}[Ehrlich 1988]
$\mathbf{No}$ is (up to isomorphism) the unique absolute linear continuum. 
 \end{Surreal}
 
	Unlike $\mathbb{R}$, however, $\mathbf{No}$ is not distinguished (up to isomorphism) from other ordered fields by its structure as an ordered class. Indeed, there are infinitely many pairwise nonisomorphic ordered fields that are absolute linear continua. What one can prove, however, is   

 \begin{Surreal}[Ehrlich 1988, 1992] $\mathbf{No}$ is (up to isomorphism) the unique real-closed ordered field that is an absolute linear continuum.  
 \end{Surreal}

\begin{sloppypar}	
	Thus, in virtue of Surreal 3, $\mathbf{No}$ is not only devoid of set-theoretically defined order-theoretic limitations, it is devoid of algebraic limitations as well; moreover, to within isomorphism, it is the unique ordered field that is devoid of both types of limitations or ``holes", as they might more colloquially be called. That is, 
$\mathbf{No}$ not only exhibits all possible algebraic and set-theoretically defined order-theoretic gradations consistent with its structure as an ordered field, it is to within isomorphism the unique such structure that does.\footnote{A historically important infinitesimalist theory of continua not touched on in this paper is that of C. S. Peirce (circa 1897, 1898, 1900). For a reconstruction (with arithmetization) of aspects of Peirce's conception based on variations of Surreal 1, Surreal 2 and Surreal 3 above, see (Ehrlich 2010). For a discussion of Peirce's theory more generally, see Matthew Moore's contribution to this volume.} It is ultimately this together with a number of related results (Ehrlich 1989, 1992) that underlies Ehrlich's contention that whereas $\mathbb{R}$ should be regarded as the arithmetic continuum (modulo the Archimedean axiom), $\mathbf{No}$ may be regarded as an \emph{absolute arithmetic continuum} (modulo NBG). 
\end{sloppypar}
	
However, to fully appreciate the nature and significance of this absolute continuum, according to Ehrlich, one must appeal to its rich \emph{algebraico-tree-theoretic} structure that emerges from the recursive clauses in terms of which it is defined. From the standpoint of Conway's construction, this \emph{simplicity hierarchical} or \emph{s-hierarchical} structure, as Ehrlich calls it (1994, 2001, 2011, 2012), depends upon ${\mathbf{No}}$'s implicit structure as a \emph{lexicographically ordered full binary tree}\footnote{A \emph{tree}  $\left\langle {A, < _s } \right\rangle $ is a partially ordered class such that for each  $x \in A$, the class  $\left\{ {y \in A:y < _s x} \right\}$ of \emph{predecessors} of $x$, written  `$pr_A \left( x \right)$', is a set well ordered by  $ < _s $. A maximal subclass of $A$ well ordered by  $ < _s $ is called a \emph{branch} of the tree. Two elements  $x$ and $y$ of $A$ are said to be \emph{incomparable} if  $x \ne y$,  $x\not  < _s y$ and  $y\not  < _s x$. The \emph{tree-rank} of  $x \in A$, written  `$\rho _A (x)$', is the ordinal corresponding to the well-ordered set $\left\langle {pr_A \left( x \right), <_s}\right\rangle $; the  $\alpha $th \emph{level} of  $A$ is $\left\langle {x \in A:\rho _A (x) = \alpha } \right\rangle$; and a root of $A$ is a member of the zeroth level. If  $x,y \in A$, then  $y$ is said to be an \emph{immediate successor} of  $x$ if  $x < _s y$ and $\rho _A (y) = \rho _A (x) + 1$; and if  $(x_\alpha )_{\alpha  < \beta }$ is a chain in $A$ (i.e., a subclass of  $A$ totally ordered by  $ < _s $), then $y$ is said to be an \emph{immediate successor of the chain} if  $x_\alpha   < _s y$ for all $\alpha  < \beta$ and  $\rho _A (y)$ is the least ordinal greater than the tree-ranks of the members of the chain. The \emph{length} of a chain  $(x_\alpha  )_{\alpha  < \beta }$ in  $A$ is the ordinal  $\beta $. A tree $\left\langle {A, < _s } \right\rangle $ is said to be \emph{binary} if each member of  $A$ has at most two immediate successors and every chain in $A$ of limit length (including the empty chain) has at most one immediate successor. If every member of  $A$ has two immediate successors and every chain in  $A$ of limit length (including the empty chain) has an immediate successor, then the binary tree is said to be \emph{full}. Since a full binary tree has a level for each ordinal, the universe of a full binary tree is a proper class. A binary tree  $\left\langle {A, < _s } \right\rangle $ together with a total ordering  $ < $ defined on  $A$ may be said to be \emph{lexicographically ordered} if for all $x,y \in A$,  $x$ is incomparable with  $y$ if and only if  $x$ and  $y$ have a common predecessor lying between them (Ehrlich 2001, p. 1234). A subtree $B$ of a tree $\left\langle {A, < _s } \right\rangle $ is said to be \emph{initial} if for each $x \in B$, $pr_B \left( x \right)= pr_A \left( x \right)$.} and arises from the fact that the sums and products of any two members of the tree are the \emph{simplest possible} elements of the tree consistent with ${\mathbf{No}}$'s structure as an ordered group and an ordered field, respectively, it being understood that $x$ is \emph{simpler than} $y$ (written $x<_s y$) just in case $x$ is a predecessor of $y$ in the tree. In (Ehrlich 1994, 2001, 2002, 2012) the just-described simplicity hierarchy was brought to the fore and made a part of an algebraico-tree-theoretic definition of ${\mathbf{No}}$. Intuitively speaking, in accordance with this approach, $x$ is simpler than $y$ if one cannot construct $y$ until one has already constructed $x$.\footnote{In Conway's development of ${\mathbf{No}}$, and most expositions thereof (e.g. Gonshor 1986; Alling 1987; Siegel 2013, ch. VIII), a surreal number $x$ is said to be \emph{simpler than} a surreal number $y$ if $x$ is constructed at an earlier stage of recursion than $y$ (i.e. the \emph{birthday} of $x$ is less than the birthday of $y$). In Ehrlich's approach, which is growing more standard in research papers on the surreals, the tree-rank of a surreal plays the role of the birthday. While tree-ranks (or birthdays) remain of importance, it is the simpler than relation defined in terms of the predecessor relation that has emerged as the theoretically central notion. The results designated Surreal 4 and Surreal 6 below only begin to show the power of the simpler than relation so defined. 

In the investigation of Conway's class of \emph{games}, which subsumes the surreal numbers, but where the corresponding tree structure is lacking, the birthday structure remains of central importance.}

Among the striking s-hierarchical features of ${\mathbf{No}}$ is that much as the surreal numbers emerge from the empty set of surreal numbers by means of a transfinite recursion that provides an unfolding of the entire spectrum of numbers great and small (modulo the aforementioned provisos), the recursive process of defining  ${\mathbf{No}}$'s arithmetic in turn provides an unfolding of the entire spectrum of ordered fields in such a way that an isomorphic copy of each such system either emerges as an \emph{initial subfield}--a subfield that is an initial subtree of ${\mathbf{No}}$ (see Note 25)--or is contained in a theoretically distinguished instance of such a system that does. In particular:

 \begin{Surreal} [Ehrlich 2001, p. 1253]  
 Every real-closed ordered field is isomorphic to an initial subfield of ${\mathbf{No}}$.
 
 \end{Surreal}
 
 Another striking s--hierarchical feature of  ${\mathbf{No}}$ is the following generalization of the Cantor Normal Form Theorem.

 \begin{Surreal}[Conway 1979, pp. 31-33]
Every surreal number can
be assigned a canonical ``proper name'' (or normal form) that is a reflection of its characteristic
s--hierarchical properties. These \emph{Conway names}, as they are sometimes called, are expressed as
formal sums of the form  $$\sum_{\alpha<\beta}  {\omega ^{y_\alpha  } .r_\alpha} $$
where  $\beta $ is an ordinal,  $\left( {y_\alpha} \right)_{\alpha  < \beta } $ is a strictly decreasing
sequence of surreals, and  $\left( {r_\alpha  } \right)_{\alpha  < \beta } $ is a sequence of nonzero
real numbers. Every such expression is in fact the Conway name of some surreal number, the Conway name of an ordinal being just its \emph{Cantor normal form} i.e. the unique Cantorian sum $$\sum_{\alpha< m}  {\omega ^{y_\alpha  } .n_\alpha} $$ equal to the given ordinal where $m$ is a natural number, $\left( {y_\alpha} \right)_{\alpha  < m } $ is a strictly decreasing series of ordinals, and $\left( {n_\alpha  } \right)_{\alpha  < m } $  is a series of nonzero natural numbers. 
 
 \end{Surreal}

\begin{figure}[h!]
\centering\includegraphics[width=0.95\textwidth]{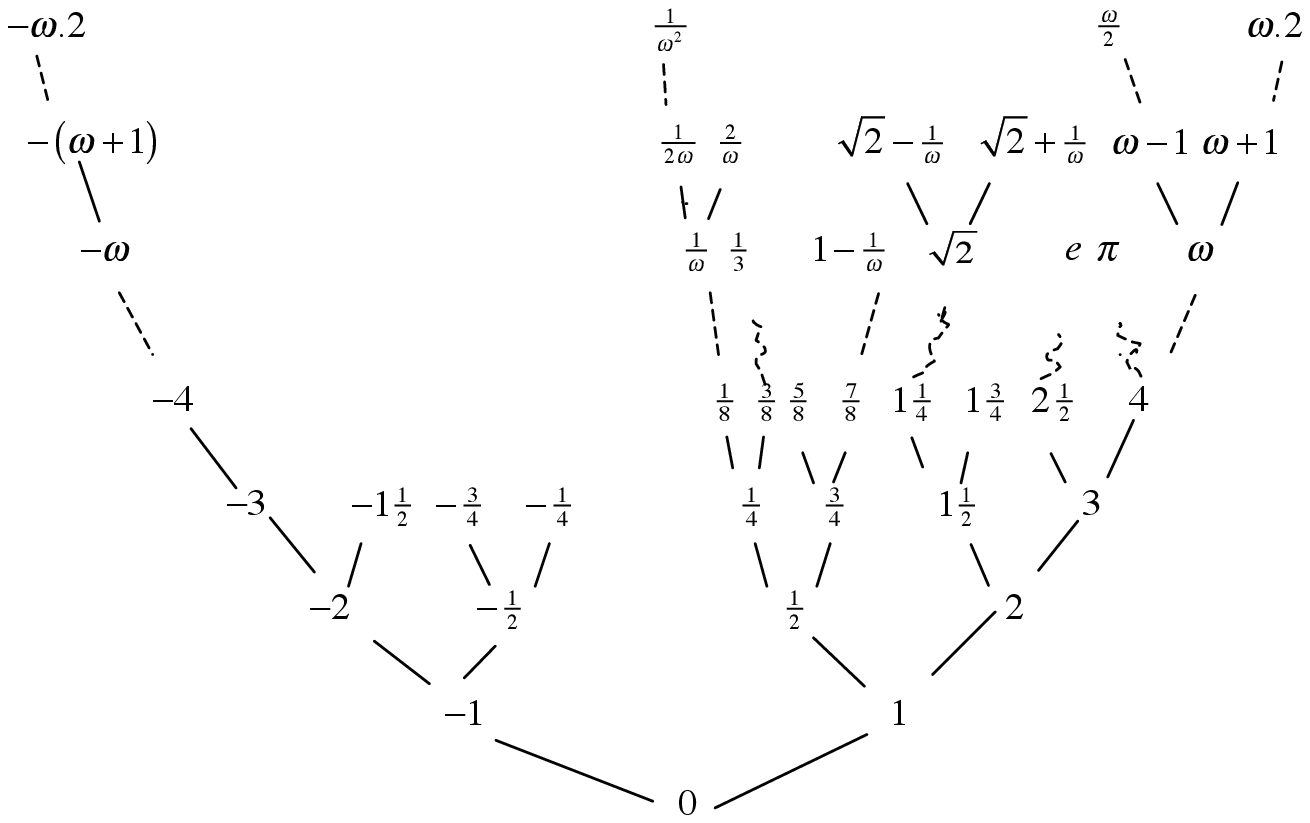}
\caption{\ }\label{Ehrlichfig3}
\end{figure}

	Figure \ref{Ehrlichfig3} below offers a glimpse of the some of the early stages of the recursive unfolding of this s-hierarchical absolute continuum, where $\omega $ is the least infinite ordinal as well as the simplest positive infinite number, $ - \omega $ is the additive inverse of $\omega $ as well as the simplest negative infinite number, $1/\omega $ is the multiplicative inverse of $\omega $ as well as the simplest positive infinitesimal number, and so on.\footnote{For a complete characterization of the recursive unfolding of $\left\langle {{\mathbf{No}},{ < _s}} \right\rangle $ in terms of Conway names of surreal numbers, see (Ehrlich 2011).}

At present there are four basic ways of constructing the surreals: Conway's inductive cut construction, which is a synthesis of the Dedekind cut construction and the von Neumann ordinal construction (Conway 1976); Conway's sign-expansion construction (Conway 1976, 2001, p. 65), which was made popular by Gonshor (1986); Ehrlich's inductive cut construction using Cuesta Dutari cuts (Ehrlich 1988, Alling and Ehrlich 1986, 1987), which is a generalization of the Dedekind cut construction; and Ehrlich's inductive generalization of the von Neumann ordinal construction (2002, 2012)\footnote{
Ehrlich (1994, p. 242) has also provided an explicit definition of the surreals as they arise in the latter construction. In addition, working in the HoTT framework, M. Shulman (2013) has introduced an interesting variation on Conway's original construction.}. In the latter construction each surreal number $x$ emerges as a canonical partition $\left( {L,R} \right)$ of the set of all surreal numbers that are simpler than $x$ in the full binary surreal number tree $\left\langle {{\mathbf{No}},{ < _s}} \right\rangle $. Although $L$ and $R$ are defined independently of the total (lexicographic) ordering $ < $ that is defined on $\left\langle {{\mathbf{No}},{ < _s}} \right\rangle $, they are found to coincide with the sets of all surreal numbers that are simpler than $x$ and less than $x$ and simpler than $x$ and greater than $x$, respectively, in the lexicographically ordered full binary surreal number tree $\left\langle {{\mathbf{No}}, < ,{ < _s}} \right\rangle $. ${\mathbf{No}}$'s ordinals (whose arithmetic in ${\mathbf{No}}$ is different than the familiar Cantorian arithmetic) are identified with the members of the ``rightmost" branch of $\left\langle {{\mathbf{No}}, < ,{ < _s}} \right\rangle $, i.e. the unique initial subtree of ${\mathbf{No}}$ that is a proper class satisfying the condition $x<y$ whenever $x<_s y$. ${\mathbf{No}}$'s system of reals is identified with the unique initial subtree of ${\mathbf{No}}$ that is a Dedekind complete ordered field, and, as is evident from Surreal 4 above, many of the the non-Archimedean number systems that have been the focus of our discussion likewise emerge as initial substructures of ${\mathbf{No}}$--substructures of ${\mathbf{No}}$ that are initial subtrees of ${\mathbf{No}}$. Indeed, in some cases, as we will see in Surreal 7, they emerge as highly distinguished initial substructures of ${\mathbf{No}}$. 

It is this algebraico-tree-theoretic s-hierarchical structure of ${\mathbf{No}}$ together with its remarkably inclusive nature that underwrites the theorem that is central to Ehrlich's portrayal of ${\mathbf{No}}$ as a \emph{unifying, s-hierarchical absolute arithmetic continuum}. To motivate the result, recall that besides being (up to isomorphism) the unique Dedekind complete ordered field,
${\mathbb{R}}$ is (up to isomorphism) the unique universal and the unique maximal (or
nonextensible) Archimedean ordered field, the condition of maximality or of 
nonextensibility being Hilbert's (1900) classical continuity condition. Analogs of these
results also hold for ${\mathbb{R}}$ considered as an $s$-hierarchical ordered field. More importantly,
however, as Surreal 6 below indicates, $s$-hierarchical analogs of these classical characterization
theorems also hold for ${\mathbf{No}}$.

Following Ehrlich (2001, p. 1236), $\langle A,+, \cdot,0,1,<,<_s \rangle$ is said to be an \emph{s-hierarchical ordered field} if $\langle A,+, \cdot,0,1,< \rangle$ is an ordered field, $\langle A,<,<_s \rangle$ is a lexicographically ordered binary tree and $+$ and $\cdot$ are defined in an analogous matter that Conway defines the field operations in ${\mathbf{No}}$, thereby ensuring the sums and products of the members of $A$ are the simplest possible members of $A$ consistent with $A$'s structure as an ordered field (Ehrlich 1994, pp. 251-252, 2001, p. 1236, 2012, pp. 15-16). If $A$ is an s-hierarchical ordered field, then $A$ is said to be \emph{universal} if every s-hierarchical ordered field is isomorphic to an initial subfield of $A$; $A$ is said to be \emph{maximal} (or \emph{nonextensible}) if there is no s-hierarchical ordered field that properly contains $A$ as an initial subfield; and $A$ is be said to be \emph{complete} if whenever $L$ and $R$ are subsets of $A$ for which $L<R$, there is a simplest member of $A$ lying between the members of $L$ and the members of $R$.

\begin{sloppypar} 
Mirroring the classical equivalences regarding Dedekind's and Hilbert's characterizations of the classical arithmetic continuum, we have
\end{sloppypar} 

\begin{Surreal}[Ehrlich 2001, p. 1239]
Let $A$ be an $s$-hierarchical ordered field. $A$ is complete if and only
if $A$ is
universal if and only if $A$ is maximal if and only if $A$ is isomorphic to ${\mathbf{No}}$.
\end{Surreal}

Whereas Surreal 1 and Surreal 3 may be said to characterize $\mathbf{No}$ as an absolute arithmetic continuum (modulo NBG), Surreal 6 may be said to characterize $\mathbf{No}$ as an \emph{absolute s-hierarchical arithmetic continuum} (modulo NBG). Moreover, since every initial subfield of ${\mathbf{No}}$ is an $s$-hierarchical ordered field (Ehrlich 2001, p. 1236), it is evident from Surreal 5 that the universal component of this characterization has real teeth, as one would expect of a purported absolute continuum. Moreover, as was intimated above, not only does every elementary continuum emerge as an initial subfield of $\mathbf{No}$, in some cases they emerge as distinguished initial substructures of ${\mathbf{No}}$. As apt illustrations, we note

\begin{Surreal}[Ehrlich 1988, 2012]
Let $\mathbf{No}(\lambda)$ be the class of surreal numbers having tree-rank $< \lambda$ and let $\alpha$ be an ordinal.  

\smallskip	
{\rm{I}}. Assuming {\rm\textbf{CH}}, $\mathbf{No}(\omega_1)$  (with $+,\cdot$, and $<$ defined \'a la Conway) is isomorphic to the underlying ordered field of the (unique up to isomorphism) saturated hyperreal number system of power $\aleph_1$.

\smallskip
{\rm{II}}. Assuming the existence of an inaccessible cardinal, $\aleph_\alpha$  being the least, $\mathbf{No}(\omega_\alpha)$ (with $+,\cdot$, and $<$ defined \'a la Conway) is isomorphic to the underlying ordered field of the (unique up to isomorphism) saturated hyperreal number system of power $\aleph_\alpha$. 

\smallskip
{\rm{III}}. ${\mathbf{No}}$ is isomorphic to the underlying ordered field of the maximal hyperreal number system that exists in {\rm{NBG}}.\footnote{There is a general result relating the underlying ordered fields of saturated hyperreal number systems to particular $\mathbf{No}(\lambda)$ (Ehrlich 1988), but since the ones commonly employed in the literature are those referred to in I and II, we have focused attention on them.}

\end{Surreal}

It is worth noting that while nonstandard analysis can be carried out in $\mathbf{No}$ (and in suitable subfields thereof), $\mathbf{No}$ was neither motivated by, nor is presently well suited for, that purpose. As Conway observed, ``$\mathbf{No}$ is really irrelevant to non-standard analysis" (1976, p. 44). Central to nonstandard analysis is the transfer theorem, and at present the only way to obtain the fruit of transfer in $\mathbf{No}$ would be by inducing it vis-\'a-vis an ordered field embedding of a given hyperreal number system. Hence, while cross-fertilization between the surreal and hyperreal number systems might in time lead to employing subfields of $\mathbf{No}$ as a vehicle for doing nonstandard analysis, at present there is little point in doing so.

As was mentioned above, assuming \textbf{CH}, an ordered field that is an infinitary pantachie in the sense of Hausdorff is isomorphic to the underlying ordered field of one of the most familiar non-Archimedean continua of our day. Indeed, this follows immediately from Part I of Surreal 7 and 

\begin{Surreal}[Ehrlich 2012]
Let $\mathbb{H}_p$ be an ordered field that is an infinitary pantachie in the sense of Hausdorff. Assuming {\rm\textbf{CH}}, $\mathbf{No}(\omega_1)$  (with $+,\cdot$, and $<$ defined \'a la Conway) is isomorphic to $\mathbb{H}_p$.\footnote{Recently, Aschenbrenner, van den Dries and van der Hoeven (2019) have shown that assuming {\rm\textbf{CH}}, the underlying ordered field of every maximal Hardy field is likewise isomorphic to $\mathbf{No}(\omega_1)$  (with $+,\cdot$, and $<$ defined \'a la Conway). At present, it is not known if without \textbf{CH} one can either show that $\mathbf{No}(\omega_1)$ is isomorphic to the underlying ordered field of a nonstandard model of analysis, to an infinitary pantachie in the sense of Hausdorff, or to a maximal Hardy field. In each case, however, it is known that such a structure is isomorphic to an initial subfield of $\mathbf{No}$ containing $\mathbf{No}(\omega_1)$ (Ehrlich 2014, p. 29).}

\end{Surreal}

The most common way of constructing a saturated hyperreal number system of power $\aleph_1$ is as the reduced power of $\mathbb{R}$ over the index set $\mathbb{N}$ modulo a nonprincipal ultrafilter on $\mathbb{N}$ (e.g. Goldblatt 1998, ch. 3). Hausdorff's construction, while similar to the just-said \emph{ultrapower} construction, differs in that he uses the reduced power of $\mathbb{R}$ over the index set $\mathbb{N}$ modulo the filter of cofinite subsets of $\mathbb{N}$ (Plotkin 2005, p. 269). Not being a maximal filter, the filter of cofinite subsets of $\mathbb{N}$ is \emph{not} an ultrafilter on $\mathbb{N}$ (e.g. Goldblatt 1998, p. 18).  Despite this difference, it is truly remarkable that the byproduct of Hausdorff's perfection of an idea rooted in du Bois-Reymond's belief that an infinitesimally enriched alternative to the Cantor-Dedekind conception of the continuum could be obtained is isomorphic (modulo \textbf{CH}) to the underlying ordered field of what is perhaps the most familiar nonstandard model of analysis employed today.

In \S 5 it was observed that, by appealing to ideas of Hahn one may shed light on the structure of $\mathbf{No}$ and on some of the most important non-Archimedean continua of our day more generally.  Although there is much to be said about this, we will simply illustrate this with respect to the non-Archimedean continua referred to in Surreal 7. While the result embraces the subfields of $\mathbf{No}$ relevant to Parts I and II of Surreal 7, unlike those results, which assume {\rm\textbf{CH}} or the existence of an inaccessible cardinal, it is provable in NBG.

\begin{Surreal}[Ehrlich 1988, 2012]

{\rm{1}}. If $\aleph_\alpha$ is a regular cardinal, then  $\mathbf{No}(\omega_{\alpha})$  (with $+,\cdot$, and $<$ defined \'a la Conway) is isomorphic to $$\mathbb{R}\left (\left( t^{\mathbf{No\left ( \omega _{\alpha} \right )}} \right )\right)_{\omega _{\alpha}.}$$

{\rm{II}}. ${\mathbf{No}}$ is isomorphic to $$\mathbb{R}\left (\left( t^{\mathbf{No}} \right )\right)_{On}.$$

\end{Surreal}

The surreal unification of systems of numbers great and small has been extended beyond the range of number systems referred to above through the work of van den Dries and Ehrlich (2001, forthcoming), Berarducci and Mantova (2018, forthcoming), Fornasiero (2016), Ehrlich and Kaplan (2018) and Aschenbrenner, van den Dries and van der Hoeven (2019, forthcoming). With the exception of (Ehrlich and Kaplan 2018), which is concerned with ordered abelian groups and ordered domains, this work deals with ordered exponential fields and ordered differential fields, making use of the exponential function on $\mathbf{No}$  introduced by Martin Kruskal and developed by Harry Gonshor (1987) and a derivation $\partial_{\rm{BM}}$ on $\mathbf{No}$ due to Berarducci and Mantova (2018).
Rudiments of analysis on the surreals have also been developed by Alling (1987), Rubinstein-Salzedo and Swaminathan (2014), and Costin, Ehrlich and Friedman (24 Aug 2015; also see, Costin and Ehrlich 2016). Costin, Ehrlich and Friedman, in particular, have developed a theory of integration (and differentiation) that extends the range of analysis from the reals to the surreals for a large subclass of \emph{analyzable functions}. The analyzable functions, which constitutes a broad generalization of the analytic functions and includes most standard functions that arise in applied analysis, was isolated by Jean \'Ecalle (1981-1985,1992,1993) in connection with his proof of Dulac's Conjecture, which is the best currently known result on Hilbert's 16th problem. Unlike nonstandard analysis, which provides an infinitesimalist approach to integration on the extended real number system $\mathbb{R}\cup \{\pm \infty\}$, the theory of surreal integration deals with integrals whose bounds and values need not be extended reals at all. For example, in the surreal theory
\[\int_{0}^{\omega}e^{x}dx=e^{\omega}-1=\omega^{\omega}-1,\] where $\omega$ is the least infinite ordinal and $e$ is the Kruskal-Gonshor exponential function on $\mathbf{No}$ that extends the familiar exponential function on $\mathbb{R}$. 

The just-mentioned works of Berarducci and Mantova, Aschenbrenner, van den Dries and van der Hoeven, and Costin, Ehrlich and Friedman all make use of the system $\mathbb{T}$ of \emph{transseries} or of a surreal isomorphic copy thereof.\footnote{For discussions of transerries, see (van der Hoeven 2006; Costin 2009; Edgar 2010; Aschenbrenner, van den Dries, and van der Hoeven  2017).}  A typical transseries is a formal sum (of finite or countable infinite length) of the form

$$e^{e^x}-3e^{x^2}+5x^{1/2}-\log x+1+x^{-1}+x^{-2}+x^{-3}+\cdot \cdot \cdot+e^{-x}+x^{-1}e^{-x},$$

\noindent
where $x$ is an infinitely large variable. Transseries naturally arise as asymptotic solutions of differential or more general functional equations. The system $\mathbb{T}$, which has proven to be of considerable interest to analysts, model theorists, differential algebrists, computer algebrists, surrealists and theoretical physicists,\footnote{With respect to model theory, we note that Aschenbrenner, van den Dries, and van der Hoeven have been awarded the 2018 Karp Prize for their work (2017) on the model theory of transseries and asymptotic differential algebra.} was introduced by Dahn and G\"oring (1987), as an ordered exponential field, in connection with their work on Tarski's problem of the decidability of the theory of real numbers with the exponential function (Tarski 1948, p. 45), and independently by \'Ecalle (1992), as an ordered differential field, in connection with his work on analyzable functions.\footnote{An \emph{ordered exponential field} is an ordered field $K$ together with an isomorphism $\exp_K$ from the additive group of $K$ onto the multiplicative group of $K$, and an \emph{ordered differential field} is an ordered field $K$ together with a map $\partial : K \rightarrow K$ such that $\partial(a+b)=\partial(a)+\partial(b)$ and $\partial(ab)=\partial(a)b+a\partial(b)$ for all $a,b \in K$.} Building on ideas from Dahn and G\"oring (1987), van den Dries, Macintyre and Marker (2001) showed that $\mathbb{T}$ may be constructed as the union of an inductively defined chain of Hahn fields, and Aschenbrenner, van den Dries and van der Hoeven (forthcoming) subsequently showed that $\mathbb{T}$, so construed, has a natural embedding $\iota$ of ordered differential fields into $\langle \mathbf{No}, \partial _{\rm{BM}} \rangle$. Making use of Berarducci and Mantova's analysis of the image of $\iota$ (forthcoming, Theorem 4.11), Elliot Kaplan has shown that the image of $\iota$ is an initial subtree of $\mathbf{No}$.\footnote{Kaplan's result is announced in (Aschenbrenner, van den Dries and van der Hoeven 2019).}

Thus, as the material in \S 4 - \S9 reveals, there are trails beginning with the work of Stolz, Veronese, Levi-Civita and du Bois-Reymond, extended further by Hahn, Hausdorff, Artin and Schreier, and further still by Robinson, Conway and a host of other nonstandard analysts and surrealists, collectively inspired by interrelated considerations regarding the arithmetic continuum and an array of generalizations thereof, that enjoy a harmonious integration in the   absolute $s$-hierarchical arithmetic continuum of surreal numbers.\footnote{Unfortunately, for lack of space, we have not been able to take into account G. H. Hardy's aforementioned development of the ideas of du Bois-Reymond, the subsequent development of the theory of Hardy fields, and the far-reaching theory of \emph{H-fields} of Aschenbrenner, van den Dries, and van der Hoeven (2017).}

\section{hjelmslev's nilpotent infinitesimalist continuum}

Heretofore we have been concerned primarily with infinitesimalist conceptions of continua that are either based upon, or can be modeled by, non-Archimedean ordered fields. However, not all infinitesimalist theories of continua have such a structure. In this and the subsequent section we will touch on the principal ones that do not.

The system of \emph{dual numbers}, commonly denoted $\mathbb{R}[\epsilon]$, is a commutative ring with identity consisting of $$\{r_0+r_1\epsilon: r_0,r_1 \in \mathbb{R} \}$$ with sums and products defined termwise, it being understood that $\epsilon \neq 0$ and $\epsilon^2 =0$. Thus, for any two dual numbers $a_0+a_1\epsilon$ and $b_0+b_1\epsilon$ $$(a_0+a_1\epsilon)+(b_0+b_1\epsilon)=(a_0+b_0) +(a_1+b_1)\epsilon$$ 
\vspace{-.2 in}
$$(a_0+a_1\epsilon)(b_0+b_1\epsilon)=a_0b_0 +(a_0b_1+a_1b_0)\epsilon.$$

\smallskip
\noindent
$\mathbb{R}[\epsilon]$ is not a field since its nonzero elements of the form $0 + a_1\epsilon$ are not invertible, and all such elements are proper zero-divisors--nonzero elements $x$ such that for some nonzero $y$, $xy=0$--which even precludes $\mathbb{R}[\epsilon]$ from being an integral domain. $\mathbb{R}[\epsilon]$ so defined is isomorphic to the quotient $\mathbb{R}[X]/(X^2)$ of the polynomial ring $\mathbb{R}[X]$ by the ideal $(X^2)$ generated by the polynomial $X^2$, which is also frequently referred to as the system of dual numbers and is likewise denoted  $\mathbb{R}[\epsilon]$.

The element $\epsilon$ in $\mathbb{R}[\epsilon]$ is a \emph{nilpotent}, an element $x$ such that $x^n =0$ for some positive integer $n$. The least $n$ for which $x^n =0$ is the \emph{index of nilpotency} of $x$. The index of nilpotency of $\epsilon$ in $\mathbb{R}[\epsilon]$ is 2. Nonzero nilpotent elements of index $n$ have many applications in algebra, one of them being a convenient way of representing quantities up to \emph{infinitesimal order} $n$. When nilpotents are thus interpreted, they are referred to as \emph{nilpotent infinitesimals}. The number systems we have hitherto considered have no nonzero nilpotent elements.

$\mathbb{R}[\epsilon]$ was introduced by William Clifford (1873) and has been widely employed by mathematicians, physicists and engineers ever since, including W. Study (1903), J. 
Gr\"unwald  (1906), and C. Segre (1911), whose work initiated the study of geometrical spaces whose coordinates are $n$-tuples of elements of a ring that is neither a field nor a division ring more generally. However, it was Johannes Hjelmslev (1923) who first employed $\mathbb{R}[\epsilon]$ to develop a \emph{nilpotent-based infinitesimalist theory of the continuum}. Unfortunately, like the aforementioned works of Veronese, Levi-Civita, Hahn, du Bois-Reymond and Hausdorff, his work is not as well known among historians and philosophers of mathematics as it deserves to be.

Hjelmslev regarded classical geometry as a crude approximation of the empirical world. In particular, inspired by the view of the pre-Socratic philosopher Protagoras, as recounted in Aristotle's \emph{Metaphysics} (III: 2, 998a), he held that the axiom that two straight lines always share at most one point is incompatible with perceptual experience, as is the assertion that a circle and a line tangent to it meet at a single point (1923, pp. 1-2). This led him to devise a finitary ``geometry of reality" (1916) or a ``natural geometry" (1923), as he later called it,  whose subject is the lines and circles of perception constructed with real rulers and compasses.  He further devised (1923, pp. 12-13) an abstraction of the latter, coordinated by $\mathbb{R}[\epsilon]$, that is a prototype of what are today called \emph{affine Hjelmslev geometries} (e.g. Lorimer 1985, pp.  87-88). In these geometries, which are coordinated by \emph{affine Hjelmslev rings}, whereas for each pair of distinct points there is a line joining them, it need not be unique. Indeed, a pair of distinct points may lie on a pair of distinct lines; when this happens the points are said to be \emph{neighboring points} and the lines are said to be \emph{neighboring lines}, the two notions of neighbor being equivalence relations. \emph{Remote points}, by contrast, are joined by a unique line and \emph{remote lines} that intersect intersect in a unique point. In these geometries a circle and a line tangent to it intersect in a nondegenerate ``infinitesimal segment" of neighboring points. Moreover, when neighboring points and neighboring lines of these geometries are equated (i.e. collected together into equivalence classes), they give rise to their more familiar affine geometric counterparts in which two points determine a line and intersecting lines as well as circles and lines tangent to them intersect in a unique point. Hjelmslev's prototype, so conceived, reduces to standard Euclidean geometry over the reals.

Hjelmslev's plane over $\mathbb{R}[\epsilon]$, henceforth $H$, consists of all ordered pairs $(A,B) \in \mathbb{R}[\epsilon] \times \mathbb{R}[\epsilon]$. As usual, a straight line of $H$ is defined by a first-degree equation $$A x + By + C = 0,$$ 
where now, however, $A, B, C  \in \mathbb{R}[\epsilon]$ and $(A,B) \notin J\times J$, $J$ being $\{r\epsilon: r \in \mathbb{R}\}$. $\mathbb{R}[\epsilon]$ admits a relational expansion to a non-Archimedean (totally) ordered ring, where the order $<'$ is defined lexicographically by the condition: $a_0+a_1\epsilon <' b_0+b_1\epsilon$ if $a_0 < b_0$ or  $a_0 = b_0$ and $a_1 < b_1$. In virtue of the just-said ordering, the points on a line of $H$ are themselves totally ordered. 

In $H$, the circle defined by the equation $$Ax^2 +By^2 =1$$ intersects the line defined by the equation $$y=1$$ in a segment consisting of all points $(r\epsilon, 1)$ where $r \in \mathbb{R}$. In virtue of the order on $\mathbb{R}[\epsilon]$, the points of the segment of intersection are isomorphic to $\mathbb{R}$ considered as an ordered set.

Hjelmslev continued to develop these ideas in an influential six-part work under the rubric  \emph{general congruence theory} (1929-1949). In the third installment (1942, p. 48), he identified for each positive integer $n$ a ring containing nilpotents having index of nilpotency $n+1$ over which one can define geometries analogous to $H$.  These are the familiar quotient rings $\mathbb{R}[X]/(X^{n+1})$, the system $\mathbb{R}[X]/(X^2 )$ of dual numbers being the smallest, consisting of elements of the form $$r_0+...+r_{n}\epsilon^{n}$$
\noindent
with $r_0,...,r_{n} \in \mathbb{R}$ and sums and products defined termwise, it being understood that $\epsilon \neq 0$ and $\epsilon^{n+1} =0$.

Inspired by these ideas, the theory of Hjelmslev geometries and corresponding Hjelmslev rings was introduced by W. Klingenberg (1954, 1954a, 1955), and through subsequent work of Klingenberg (1956), Veldkamp (1981, 1985, 1987, 1995)  and others these geometries and corresponding rings have undergone substantial generalization and variation, yielding, for example, geometries in which not all neighboring points are connected by a line. In this family of neighbor-based geometries, Hjelmslev's idea that lines, lines and circles, and so on may intersect in infinitesimal segments characterized using nilpotent infinitesimals plays a central role.\footnote{For English presentations of Klingenberg's axiomatization, see (Dembrowski 1968, Appendix 7.2) and (Lorimer 1985), the latter of which contains detailed discussions of Hjelmslev planes and rings. For an extension of Klingenberg's axiomatization to higher dimensional spaces, see (Kreuzer 1987, 1988). For an overview with references of Hjelmslev geometries, Hjelmslev rings and various generalizations of both, see (Veldkamp 1995). For F. Bachmann's  generalization of Hjelmslev planes and their incorporation into his geometry of reflections, see (Bachmann 1971, 1973, 1973a, 1989). In Bachmann's treatment, a pair of points may either have more than one line or no line joining them, both conceptions he traces to Hjelmslev. Finally, for finite Hjelmslev planes, see (Drake and Jungnickel 1985).} These ideas, which emerged from Hjelmslev's nilpotent-based infinitesimalist conception of the geometrical continuum, have found application in other areas of mathematics as well, including  algebra, algebraic geometry, differential geometry and, most recently, in a number of nilpotent-infinitesimalist approaches to those portions of differential geometry concerned with \emph{smooth} (i.e. infinitely differentiable) \emph{manifolds} and various generalizations thereof. It is to the latter that we now turn.	

\section{infinitesimalist approaches to differential geometry of smooth manifolds and their underlying continua}
Abraham Robinson begins the \emph{Elementary Differential Geometry} section of his book on nonstandard analysis by observing that: 
\begin{quote}
Quite early in the history of the Differential Calculus, infinitesimals were employed in the development of the theory of curves, and long after the classical 
$\epsilon,\delta$-method had displaced the naive use of infinitesimals in Analysis they survived in Differential Geometry and Physics. Even now there are many classical results in Differential Geometry which have never been established in any other way, the assumption being that somehow the rigorous but less intuitive, $\epsilon,\delta$-method would lead to the same result. So far as one can see without a complete check this assumption is usually correct.  (1966, 1974, p. 83). 
\end {quote}
He then goes on to illustrate a point he had already made in earlier works, namely, using nonstandard analytic techniques 

\begin{quote}
we may now justify the use of infinitesimals in all these problems [from differential geometry] directly. (1961, p. 437; also see 1963, p. 265) 
\end{quote}

It is therefore not surprising that influenced in part by Robinson's placement of analysis on a logically sound infinitesimal foundation, logically sound infinitesimal foundations for portions of differential geometry would likewise be sought. At present there are at least four such approaches, at varying stages of development, all of which are concerned with smooth manifolds or generalizations thereof: \emph{Synthetic Differential Geometry} (SDG), introduced by F. W. Lawvere in 1967 (see Note 41)  and developed by Anders Kock (1981, 2010) and others; \emph{Infinitesimal Differential Geometry} (IDG) developed by Paolo Giordano (2001, 2010, 2011; 2016 with Wu); \emph{Topological Differential Calculus (over General Base Fields and Rings)} (TDC), introduced by Wolfgang Bertram (2008)\footnote{The appellation ``Topological Differential Calculus," which was requested by Bertram, comes from the title of his (2011).}; and \emph{Nonstandard Differential Geometry} (NSDG), based on nonstandard analysis, the most recent version being that of Tahl Nowik and Mikhail Katz (2015)\footnote{For earlier treatments of portions of differential geometry based on nonstandard analysis, see (Stroyan and Luxemburg 1976; Stroyan 1977; Lutz and Goze 1981; Almeida, Neves and Stroyan 2014).}.  In simplest terms, differential geometry of smooth manifolds is concerned with the application of calculus to spaces that locally behave like some Euclidean space $\mathbb{R}^n$. This includes Riemannian geometry, the branch of differential geometry that studies smooth manifolds with a Riemannian metric.

Whereas NSDG makes sole use invertible infinitesimals, SDG, IDG and TDC are developed using nilpotent infinitesimals. Moreover, unlike SDG, IDG and NSDG, which provide infinitesimal-based approaches to classical differential geometry of smooth manifolds (albeit sometimes within enriched frameworks), TDC seeks to provide a vast generalization of (the differential aspects of) classical differentiable geometry of smooth manifolds by providing a framework that is applicable to a wide range of underlying ``model[s] of the continuum" (Bertram 2008, p. 1) including the ordered field of reals, systems of hyperreals, the ring of dual numbers and ``good" \emph{topological fields} and \emph{topological commutative unital rings} more generally.\footnote{By a \emph{topological ring (field)}, Bertram means a commutative unital ring (field) $\mathbb{K}$ together with a topology (assumed to be Hausdorff) whose ring operations $+:\mathbb{K} \times \mathbb{K} \rightarrow \mathbb{K}$ and $\cdot:\mathbb{K} \times \mathbb{K} \rightarrow \mathbb{K}$ are continuous and for which the set $\mathbb{K}^\times$ of invertible elements is open in $\mathbb{K}$ and the inversion map $i:\mathbb{K}^\times\rightarrow \mathbb{K}$ is continuous. Bertram says $\mathbb{K}$ is \emph{good} if in addition $\mathbb{K}^\times$ is dense.}

Since the underlying continua in NSDG are hyperreal number systems, which we have already touched on, and of the remaining three theories SDG is at present the most developed, we will direct the bulk of our brief remarks to it, its prehistory, and some of the important differences that exist between it and the other theories, the latter of which deserve more attention than space will allow.

At roughly the time Klingenberg was developing the theory of Hjelmslev rings, the system of dual numbers was undergoing a distinct though overlapping  generalization and incorporated into differential geometry by Andr\'e Weil (1953) to treat Charles Ehresmann's (1951) theory of \emph{jets} and \emph{prolongations}, which are concerned with entities that generalize the notion of a tangent vector. In addition to Ehresmann's now classic work, Weil was motivated by the desire to ``return to the methods of Fermat in the first-order infinitesimal calculus" (1953, p. 111), methods that were introduced by Fermat to treat tangent constructions making use of nilpotent infinitesimals. Soon thereafter, Alexander Grothendieck (with the assistance of Jean Dieudonn\'e) (1960, 1971), making use of the system of dual numbers (and generalizations thereof), developed a nilpotent-infinitesimal-based algebraic treatment of differential calculus adequate for the needs of algebraic geometry, which remains a prominent fixture of the theory till this day.\footnote{In classical algebraic geometry the rings that were studied are integral. In Serre's subsequent treatment zero-divisors were permitted but there were no nonzero nilpotents. Grothendieck considered arbitrary commutative rings with identity, thereby permitting nilpotent infinitesimal elements. Differential properties of functions are of critical importance in algebraic geometry and this is what motivated Grothendieck's nilpotent-infinitesimal-based algebraic approach to differential calculus for the subject. For an introduction to these ideas, see (Shafarevich 2005, \S7: \emph{The Algebraic View of Infinitesimal Notions}) and  (Perren 2008, ch. 5 and Appendix B).

Nilpotent infinitesimals already appear on the second page of mathematics of the first installment of Grothendieck's \emph{\'El\'ements de g\'eom\'etrie alg\'ebrique I} (1960, p. 12), and they are found throughout the subsequent installments published from 1961 through 1967. In the revised edition, a paragraph was added explicitly stating their importance for the theory (1971, p. 11). However, while it was Grothendieck who made the rigorous use of nilpotent-infinitesimalist ideas a mainstay of algebraic geometry, an informal use of infinitesimals in algebraic geometry is already found in important work of Enriques of the 1930s. For a fascinating mathematico-historical discussion of how Grothendieck's nilpotent-infinitesimalist ideas can be employed to lend precision to Enriques's ideas, see (Mumford 2011).}  In Grothendieck's approach to algebraic geometry, as in Hjelmslev's treatment of the geometric continuum, curves meet their tangents at nondegenerate infinitesimal segments or arcs.

In May of 1967, utilizing an extension of Grothendieck's just-said ideas and incorporating them into a topos-theoretic setting, F. W. Lawvere proposed a nilpotent-infinitesimal-based approach to differential geometry of smooth manifolds\footnote{This took place as part of a series of lectures--summarized in (Lawvere 1979)--given at the University of Chicago. Among the central goals of Lawvere's lectures was to employ category theory in a substantial manner in physical theorizing. Given the role played by differential geometry therein, the connection was natural. For remarks on the critical role of category theory, in general, and topos theory, in particular, in the historical development and formulation of SDG, see (Kock 1981; McLarty 1990, \S 5; Bunge, Gago, San Luis 2018).} in which Weil's generalization of the dual number system--\emph{Weil algebras}, as they are now called--would come to play a prominent role.\footnote{For a general introduction to Weil algebras, see (Moerdijk and Reyes 1991, ch. I). For the present, we simply note that for each positive integer $n$, the Hjelmslev ring $\mathbb{R}[X]/(X^{n+1})$ is a simple example of a Weil algebra.} Unlike Robinson, who was stimulated by Leibniz's idea that the properties of infinitesimals should reflect the properties of the reals, Lawvere's ideas more closely mirror the heuristic ideas of geometers (like Cartan, Hjelmslev, Klingenberg, Lie, Weil and Grothendieck) who envision a vector tangent to a surface at a point as a tiny arc of a curve having the vector tangent to it. Building on Lawvere's ideas, G. Wraith, G. E. Reyes, E. J. Dubuc, and A. Kock developed the basic ingredients of the theory, which were given a systematic treatment by Kock (1981) under the now familiar rubric ``Synthetic Differential Geometry".

SDG is developed in an axiomatic framework in which a line is modeled by a commutative ring $\Re$ with unity containing a subset $D=\left \{d \in \Re: d^{2}=0 \right \}$ of nilpotent infinitesimals that satisfies the 
	
\medskip	
\noindent
\emph{Kock-Lawvere axiom}: For every mapping $f:D\rightarrow \Re$, there is a unique $b \in \Re$, such that for all $d \in D$, $f(d)=f(0)+d \cdot b$.

\medskip
\noindent
Geometrically speaking, this axiom asserts that the graph of every function $f:D\rightarrow \Re$ is a piece of the unique straight line through $(0,f(0))$  with slope $b$. It is a consequence of this assumption that in SDG a tangent vector to a curve $C$ at a point $p$ is a nondegenerate infinitesimal line segment around $p$  coincident with $C$.

For motivational purposes, one can employ the Kock-Lawvere axiom as stated above, where $\Re=\mathbb{R}[\epsilon]$ (without the lexicographic total ordering). For the development of SDG, however, the Kock-Lawvere axiom is replaced by a more sophisticated axiom scheme, due to Kock (1981, p. 64), which is formulated for an arbitrary Weil algebra. The importance of Weil algebras for SDG was first recognized by G. E. Reyes and G. Wraith (1978), and subsequently used to great effect by E. J. Dubuc, who introduced and developed the critical idea of a \emph{well-adapted model} of SDG and established the existence of such models making use of Weil algebras (1979, 1981, 1981a). Well-adapted models are of fundamental importance for SDG because they are precisely the models of SDG that make it possible to recover theorems of classical (i.e. standard $\epsilon, \delta$) differential geometry. That is, it is the existence of well-adapted models that makes the axioms of SDG relevant to classical mathematics. 

Using the Kock-Lawvere axiom, one may define the \emph{derivative} $f'(x)$ of \emph{any} function $f:\Re \rightarrow \Re$ at $x \in \Re$ to be $b$. More specifically, using the Kock-Lawvere axiom one may derive the following version of 

\medskip	
\noindent
\emph{Taylor's formula}: For any function $f:\Re \rightarrow \Re$ and any $x \in \Re$, $$f(x+d)=f(x)+d \cdot f'(x) \;\;\;\; \forall d \in D. $$

At first sight this appears to be impossible since it implies that all functions are smooth, which of course is not the case. Indeed, by a classical result of topology due to Banach (1931), even if we limit ourselves to the space of continuous functions, the set of functions having a derivative at a point is \emph{meager} (negligible in a precise sense). However, as authors of works on SDG hasten to add, the seeming catastrophic consequence of the Kock-Lawvere axiom does not materialize since the logic employed in SDG is \emph{intuitionistic} as opposed to classical, and the impossibility cannot be established using \emph{intuitionistic logic}, the latter of which results from classical logic by omitting the law of excluded middle.\footnote{The idea of introducing the derivative in the manner described above has roots that predate SDG. Neder (1941), for example, employed this approach (albeit, of course, not for all functions) in his attempt to develop a ``Leibnizian differential calculus with actual infinitesimal magnitudes (differentials) of the first order"  based on the system of dual numbers ``that is free from contradiction" (p. 251). Indeed, according Neder: ``Optimists may hope that one day the considerations presented here will be used to give a new (and by no means disparaged) justification for differential calculus" (p. 251). Moreover, according to Neder (p. 252), the basic idea was essentially in (Peterson 1898) and (Predella 1912). What Neder does not mention is that J. Petersen is J. Hjelmslev (see \S 10). Thus, Hjelmslev not only anticipated the idea of a nilpotent-infinitesimalist conception of a continuum based on the system of dual numbers, he also envisioned therein the rudiments of a corresponding theory of differentiation. The latter idea was further developed in (Gesztelyi 1958) and in the aforementioned work of Grothendeick. These ideas are also central to the theory of ``automatic differentiation" over the system of dual numbers, which has a vast literature whose roots go back to the 1950s.}

Neither the Kock-Lawvere axiom nor Kock's generalization thereof implies results about basic integration and to obtain such results one requires further axioms. In basic treatments, one usually introduces the following axiom due to Kock and Reyes (1981):

\medskip	
\noindent
\emph{Integration axiom}: For every function $f:[0,1]\rightarrow \Re$, there is a unique function $g:[0,1]\rightarrow \Re$, such that $g' = f$ and $g(0)=0$. 

\medskip
\noindent
In harmony with standard notation, for each $x$ in $[0,1]$, the value of the hypothesized function $g$ is denoted $$\int_{0}^{x}f(t)dt.$$ It is worth noting, however, that the integral in SDG has a different character than it does in ordinary calculus, where it is defined as the limit of an infinite sum, not as an antiderivative. Thus far, however, no one knows how to treat infinite sums or limits in SDG.\footnote{Remarks of Kock (25 Sept 2017, p. 23) suggest that he does not preclude such a possibility down the road, albeit perhaps with modifications of SDG.}

Unlike the Kock-Lawvere axiom or Kock's generalization thereof, the Integration axiom presupposes an order $\leq$ on the ring $\Re$ to make sense of the notion of the interval $[0,1]$. In particular, the order is assumed to be a \emph{preorder} (a reflexive and transitive relation), and on the basis of this one defines $$[a,b]=\{x\in \Re:a\leq x \leq b\}.$$ 
\noindent
It is further assumed that $\leq$ satisfies conditions which ensures its compatibility with the unitary ring structure of $\Re$,\footnote{$\forall xyz(x\leq y \rightarrow (x+z \leq y+z))$; $\forall xy((0\leq x \wedge 0\leq y) \rightarrow 0\leq xy )$; $0\leq1 \land \neg(1\leq0)$.} and in addition it is assumed that 
$$(*)\;\;\;\;\;\;\forall d \in D (0\leq d \land d\leq 0).$$

It is important to note that $\leq$ is not a total ordering, or even a partial ordering, since if it were assumed that $\leq$ is antisymmetric (i.e. $\forall xy((x\leq y \land y\leq x) \rightarrow x=y)$), it would follow from $(*)$ that $D=\{0\}$. But this is not the case; in fact, it is a consequence of the Kock-Lawvere axiom that $D\neq\{0\}$. On the other hand, since the underlying logic of SDG is intuitionistic, in SDG one cannot prove $\exists d \in D(d\neq 0)$ from $D\neq\{0\}$, as one might expect, and would be the case if the underlying logic were classical.  In fact, since (by definition) $a<b$ if and only if $a\leq b \land a\neq b$, if one could prove $\exists d \in D(d\neq 0)$, it would follow from $(*)$ that $\exists d \in D(d< 0 \land 0<d)$, which is impossible (even in intuitionistic logic). 

While the just-said axioms of SDG permit the development of swaths of basic calculus (Belair 1981; Belair and Reyes 1982; Bell 2008), they go just so far in the development of differential geometry of smooth manifolds. Whereas they provide the basis for developing many aspects of the theory concerned with infinitesimal behavior (Kock 2010), they provide little basis for the treatment of local finite behavior and are entirely silent on issues of a global nature. For example, they do not imply the existence and uniqueness of solutions of ordinary differential equations, permit the treatment integration of vector fields, or provide the means for the representability of germs of a smooth nature. As a result, over the years a variety of additional axioms have been proposed to make the study of SDG closer to what is ordinarily considered to be the analysis of differentiable functions of a real variable, early examples being those suggested by L. Belair (1981) and C. McLarty (1983).  However, the most far reaching work in this direction was initiated by Marta Bunge and Dubuc (1987), carried further in (Bunge and Gago 1988; San Luis 1999), and recently integrated and extended under the appellation \emph{Synthetic Differential Topology} (Bunge, Gago, San Luis 2018). In this work (SDT), not only are the three just-mentioned issues addressed, but additional topics such as Morse theory of singularities and the theory of stable germs of smooth mappings are covered including a proof of Mather's theorem.

	As in Hjelmslev geometries, the notion of neighborhood and corresponding incidence structure is of fundamental importance in SDG (Kock 2003, 2010, 25 Sept 2017). For example, it is a consequence of the Kock-Lawvere axiom that in SDG, unlike in Euclidean geometry, there are pairs of points in the plane that are \emph{not} connected by a unique straight line. However, unlike in Hjelmslev geometries, in SDG there are pairs of points in a plane that are not connected by any line at all. This arises in part from the fact that whereas the nilpotent infinitesimals in Hjelmslev geometries have ``a quantitative (linear ordered) character," those employed in SDG do not (Kock 2003, pp. 226-228). In this respect, the lines of SDG more closely resemble some of the lines of geometric spaces with neighbor relations investigated by Veldkamp (e.g. Veldkamp 1985). In these spaces, unlike those of Hjelmslev, the neighbor relation is \emph{not} transitive. This permits a hierarchy of neighborhoods--first order, second order and so on for each positive integer $n$--as is found in the algebraic geometry of Grothendieck, which lent inspiration to SDG.\footnote{For further discussion of the differences in the underlying continua of Hjelmslev's geometry and SDG, see (Kock 2003); for a detailed discussion of the significance of the presence or absence of the transitivity of the neighbor relations, see (Veldkamp 1985, 1995);  and for a penetrating discussion of neighborhoods in the context of SDG, see (Kock 25 Sept 2017).}  
	
	A space $X$ in SDG is said to be \emph{indecomposable} if there are no disjoint nonempty subsets $U$ and $V$ of $X$ such that $U \cup V =X$ (Bell  2001, 2005, 2008, 2009). There are models of SDG in which a classical space $\mathbb{R}^n$  has a counterpart $X$ that is indecomposable if $X$ is connected. John Bell takes this to imply that Òthe connected continua of SDG are true continua in something like the Anaxagoran senseÓ (1995, p. 56, 2009). In this respect, they are also reminiscent of the unsplittable continuum of Brouwer; however, the similarity is not perfect and varies depending on the axioms adopted for SDG (Bell 2001, 2009)\footnote{According to Bell (2001, p. 22, 2005, p. 302), the critical axiom is: $\forall xy (\neg x=y \;{\rm implies} \;x<y  \; {\rm or} \; y<x)$.   Though not in (Kock 1981), this assertion is a consequence of axioms assumed in (Moerdijk and Reyes 1991) and (Bunge, Gago, San Luis 2018).} 

	Another respect in which SDG is similar to Brouwer's theory is the failure of the intermediate value theorem in its underlying theory of analysis. In fact, in SDG, unlike in Brouwer's system, the theorem even fails for some polynomials (Moerdijk and Reyes 1991, pp. 317-318), a failure that runs contrary to the thinking of Leibniz and Euler let alone Bolzano, Cauchy, and Weierstrass. Accordingly, while SDG may in time provide a viable treatment of the smooth aspects of differential geometry, its underlying analysis may not be as well suited to provide a natural alternative for classical analysis, at least not if it hopes to mirror the latter's most central ideas regarding continuity.
	
How then should SDG/SDT be viewed vis-\'a-vis classical differential geometry/topology and classical analysis? Not surprisingly, the answer one gets depends on who one asks. In a recent statement, Kock expresses his view thus:  

\begin{quote}
I prefer not to think of SDG as a monolithic global theory, but as a method to
be used locally, in situations where it provides insight and simplification of a
notion, of a construction, or of an argument. The assumptions, or axioms that
are needed, may be taken from the valuable treasure chest of real analysis. (Kock 25 Sept 2017, p. 24)
\end{quote}

These remarks suggest that Kock does not see SDG so much as a theory intended to replace classical differential geometry, but rather as a methodological accompaniment that should be part of the classical differential geometer's toolkit (to borrow a term from Keisler) much as classical differential geometry is part of the synthetic differential geometer's toolkit, though the roles they would serve would be different. 

In a private communication, Marta Bunge voices a somewhat similar sentiment, writing:
\begin{quote}
\vspace{-.0375 in}
I agree with the view expressed by Anders in that SDG (or SDT) is not meant to replace classical differential geometry (or topology) but to be inspired by it and enhance it. Due to the rich structure of a topos and to the presence of infinitesimals, classical differential geometry (or topology) can be studied in a context that permits better conceptual simplification and unification in ways that often differ from those of its classical sources. I have stated this point of view in my paper ``Toposes in Logic and Logic in Toposes" (1984) intended for philosophers and in the Introduction to SDT. However, I do view SDG (or SDT) as a theory and not just as a method. In it some developments are possible which are totally non classical some of which originating in the work of Jacques Penon (1981, 1985).\footnote{Bunge is referring to the fact that in addition to the nilpotent infinitesimals of SDG, SDT employs what she calls ``logical infinitesimals" (in particular  $\Delta = \neg \neg \{0\}$) and also ``logical (or Penon) opens" neither of which could be discussed in a classical setting.} It is as a theory that and therefore one can consider models of it, preferably well-adapted in that, by means of such models, themselves often inspired by classical differential geometry (or topology), one can recover classical results. The well adapted part is I think important and Archimedeanness of the ring of line type is inherent in it. (M. Bunge, Private Communication, November 26, 2017)\footnote{As Bunge alludes to above, among the assumptions adopted in SDT is what contributors to SDG typically call the \emph {Archimedean condition}, the assertion $$(A)\;\;\;\;\;\forall x \in \Re \exists n\in \mathbb{N} (x<n).$$ Perhaps it would be more appropriate to say $(A)$ requires $\Re$ to be \emph{finitely bounded} or \emph{weakly Archimedean}, since there is a vast array of non-Archimedean ordered rings that satisfy this weak form of the Archimedean condition for rings, including the aforementioned ones of Hjelmslev. In fact, if $F$ is a non-Archimedean ordered field, then the non-Archimedaen subring of $F$ consisting of all the finite and infinitesimal members of $F$ satisfies $(A)$.}   
\end{quote}

At present, in addition to five axioms for SDG, which are taken to include the Archimedean axiom, SDT is based on two basic axioms, which address the  representability of germs of smooth mappings and the existence and uniqueness of solutions of ordinary differential equations, and four specialized postulates that are employed to treat different parts of the theory (Bunge, Gago, San Luis 2018, \S 2.1, \S 6.1-\S 6.2).\footnote{As the last line of Bunge's remark suggests, the assumption of the Archimedean axiom in SDT is motivated by the desire for a well-adapted model. Indeed, the main theorem of (Bunge and Dubuc 1986), developed in collaboration with Andr\'e Joyal, ``explains in part why the axiom is important if one wants to develop a theory (or imply a method) which does not go against classical mathematics but which incorporates it...." (M. Bunge, Private Communication, June 26, 2018)}  Bunge anticipates that additional axioms may be needed to further extend SDG/SDT and she has expressed confidence that the \emph{Dubuc topos} (Dubuc 1981), which is the only known well-adapted model for SDT (Bunge, Gago, San Luis 2018, pp. 5, 182 and ch. 12), will be found to model them.\footnote{Private communication.}

A reason that Kock offers for viewing SDG as a method rather than as a theory is because he believes that 
\begin{quote}
full fledged analysis in axiomatic terms, incorporating SDG, quickly becomes
overloaded with axioms, and is better developed as a descriptive theory,
describing what actually holds in specific models.... (Kock 25 Sept 2017, p. 24)
\end{quote}

Despite the methodological difference between Kock and Bunge, it is worth noting the difference between the attitudes expressed by these two leading architects of SDG/SDT and those of some of its champions, such as John Bell, who on occasion do not portray real analysis as a treasure chest from which to draw upon and seek partnership with, but rather as an edifice with paradoxical consequences of the axiom of choice and unwarranted applications of the law of excluded middle that need ``jettisoning", at least if one wishes to obtain ``a faithful account of the truly continuous" (Bell 2008, p. 5, 2005, pp. 294-297, 2009).\footnote{Hellman (2006) and Shapiro (2014) offer philosophical discussions of SDG and its underlying analysis viewed as constructivist replacements for standard analysis and portions of standard differential geometry. We suspect that philosophical assessments of SDG/SDT and its underlying analysis from the more nuanced perspectives of Kock and Bunge, as expressed in the above quotations, also would be worthy of attention. For remarks on SDG from an author who ``invites you to leave the overcrowded Cantor's paradise" (p. 496), see (Bauer 2017), and for further details on Bell's view, see his contribution to the present volume.} 
	
As was noted above, it is the Kock-Lawvere axiom that underlies the incompatibility of SDG with classical logic (e.g. Kock 1981, pp. 2-5); Lavendhomme 1996, pp. 2-5). In particular, it is the contention that the axiom applies to \emph{every} mapping $f:D \rightarrow \Re$ that runs afoul of the law of excluded middle. Accordingly, since the idea that a tangent vector to a curve intersects the curve in a nondegenerate infinitesimal line segment is entirely compatible with classical logic as is the employment of nilpotent infinitesimals for the representation of such intersecting segments, it is natural to inquire if it is possible to develop a variation of SDG in which the underlying logic is classical rather than intuitionistic. It is the idea of providing just such a variation of SDG that motivates Giordano's IDG. 

At the heart of the original formulation of IDG is the idea of replacing the Kock-Lawvere axiom by a theorem asserting: for each smooth function $f:\mathbb{R} \rightarrow \mathbb{R}$  (that is pathwise Lipschitzian) there is a function ${^\bullet}f:{^\bullet}\mathbb{R} \rightarrow {^\bullet}\mathbb{R}$ extending $f$, and a unique $b \in \mathbb{R}$ such that $$\forall d \in D: {^\bullet}f(d)={^\bullet}f(0)+d\cdot b,$$
\noindent
where  $^{\bullet}\mathbb{R}$ is Giordano's ring of \emph{Fermat reals} that extends $\mathbb{R}$, $D=\left \{d \in {^\bullet}\mathbb{R}: d^{2}=0 \right \}$ and the functions ${^\bullet}f:{^\bullet}\mathbb{R} \rightarrow {^\bullet}\mathbb{R}$ are analogs for $^{\bullet}\mathbb{R}$ of the standard smooth functions $f: \mathbb{R} \rightarrow \mathbb{R}$ (Giordano 2001).\footnote{The universe of $^{\bullet}\mathbb{R}$ is introduced in three steps. One starts with the class $\mathbb{R}_{o}[t]$ of ``little-oh polynomials", i.e. functions $x:\mathbb{R}_{\geq 0}\rightarrow \mathbb{R}$ that can be written as $$x(t)=r+\sum_{i=1}^{k}a_{i}\cdot t^{a_{i}}+o(t)\;\;\;{\rm {as}} \;\;t\rightarrow 0^{+},$$  where the coefficients, powers and $r$ are reals, $k \in \mathbb{N}$ and $o$ is a variant of Landau's little-$o$ notation. Next, one defines equivalence between $x,y \in \mathbb{R}_{o}[t]$ by the condition: $x \sim y$ if and only if $x(t)=y(t)+o(t)$ as $t \rightarrow 0^+$. Finally, $^{\bullet}\mathbb{R}$ is taken to be the  quotient set of $\mathbb{R}_{o}[t]$ with respect to $\sim$.

For the order and ring operations on $^{\bullet}\mathbb{R}$, as well as detailed discussions of its properties, see  (Giordano 2010a; Giordano 2011a; Giordano and Kunzinger 2013).} 

Thus, in addition to the tuples of reals contained in $f$, $^{\bullet}f$ contains tuples of nonstandard Fermat reals. Each Fermat real number $x$ differs from a unique real number $^{\rm{o}}x$ by an infinitesimal amount, and can be written in a unique form $$x={^{\rm{o}}x+\sum_{i=1}^{k}\alpha_{i}\cdot {\rm{d}}t_{a_{i}}},$$ where $^{\rm{o}}x$ is the \emph{standard part }of $x$ (the unique real number closest to $x$), $k \in \mathbb{N}$, the $\alpha_{i}$s are nonzero real numbers and the ${\rm{d}}t_{a_{i}}$s constitue a descending sequence of nilpotent infinitesimals (in the sense that ${\rm{d}}t_{a_{n}}$ is infinitesimal relative to ${\rm{d}}t_{a_{m}}$ whenever $1 \leq m<n \leq k$), it being understood that $x= {^{\rm{o}}x}$ if $k=0$ (e.g. Giordano 2010, 2011 \S11). Moreover, for each $n>1$, there are nilpotent infinitesimals in ${^\bullet}\mathbb{R}$ having index of nilpotency $n$.

However, despite their overlapping motivations, there are notable differences between Giordano's theorem and the Kock-Lawvere axiom. For example, in addition to specifying that the slope $b$ is a real number rather than a member of $\Re$, the theorem refers to all standard smooth functions (that are pathwise Lipschitzian) as opposed to \emph{all} functions $f:D \rightarrow \Re$. Moreover, unlike $\Re$ and the Weil algebras more generally employed in SDG, $^{\bullet}\mathbb{R}$ is \emph{totally ordered}. Accordingly, since the set $D$ in Giordano's theorem contains nilpotent infinitesimals that are provably nonzero, the totally ordered commutative unitary ring $^{\bullet}\mathbb{R}$ is necessarily non-Archimedean despite the fact that it satisfies the finitely bounded version of the Archimedean condition employed in SDG. Furthermore, whereas the product of any two members of $D$ is equal to $0$ in IDG (Giordano 2010: Theorem 24, p. 172), this is not the case in SDG, the latter being a consequence of the Kock-Lawvere axiom (cf. Kock 1981, pp. 6, 15). In this respect, the nilpotent infinitesimals employed in IDG more closely resemble those employed in Hjelmslev geometries.\footnote{In his (2010), Giordano suggests that due to the just-stated difference his nilpotent infinitesimals are more intuitively satisfying than those employed in SDG. In particular, he suggests that an intuitive picture of infinitesimal segments of lengths $h$ and $k$ where $h^{2}=k^{2} =0$ and $h \cdot k \neq 0$ is not possible. However, he has since recanted and now more modestly maintains that in the end ``it is a matter of taste about what approaches are felt as beautiful, manageable and in accordance with our philosophical approach to mathematics" (Giordano and Wu 2016, p. 898). For historical motivation as well as remarks regarding the impact this difference has on SDG, see (Kock 2003, p. 228). Also see (Hellman 2006, \S 4) for philosophical remarks bearing on the difference.}

 Like SDG, IDG is developed in a category-theoretic framework. However, unlike SDG, IDG as developed above does not lead to a Cartesian closed category of spaces.\footnote{In SDG, one is assured of working in a Cartesian closed category since every topos is a Cartesian closed category. For an overview of the virtues of Cartesian closed categories and some of the roles they play in a variety of mathematical and physical theories, see (Giordano 2011a).} To address this, Giordano has replaced the just-described treatment based on pathwise Lipschitzian functions with one based on \emph{quasi-standard smooth functions}. In their most basic form, these are functions $g:{^\bullet}\mathbb{R} \rightarrow {^\bullet}\mathbb{R}$ that locally can be written as $g(x)={^\bullet}\varphi (x,p)$, where the $p \in {^\bullet}\mathbb{R}^n$ are nonstandard fixed parameters.\footnote{The class of quasi-standard smooth functions is a ${^\bullet}\mathbb{R}$ analog of the class of \emph{quasi-standard functions} employed by Robinson in the early days of nonstandard analysis (Robinson 1961, pp. 437-438, 1963/1965, pp. 266-268; Giordano 2011, p. 864; Giordano and Wu 2013, p. 890) to deal with the theory of distributions. It was subsequently replaced by the more versatile notion of an \emph{internal} function (Robinson 1966, 1974, p. 42).} The derivative of such a function is defined by appealing to what Giordano calls the \emph{Fermat-Reyes Theorem}, which asserts: for each quasi-standard smooth function $g$ there is precisely one quasi-standard smooth function $r:{^\bullet}\mathbb{R}^{2} \rightarrow {^\bullet}\mathbb{R}$ such that $$ \forall (x,h) \in {^\bullet}\mathbb{R}^{2}:g(x+h)=g(x)+h\cdot r(x,h).$$

\noindent
This permits one to define $g'(x)=r(x,0)$ for all $x \in {^\bullet}\mathbb{R}$. Integration for quasi-standard smooth functions is developed by proving the existence and uniqueness of primitives. 

On the basis of this (or rather a treatment based on a more general class of quasi-standard smooth functions $^{\bullet}f:S \rightarrow T$ where $S$ and $T$ are subsets of ${^\bullet}\mathbb{R}^n$ and ${^\bullet}\mathbb{R}^m$, respectively), Giordano and Wu have constructed a calculus that is at least as strong as its counterpart in SDG (Giordano 2011; Giordano and Wu 2016, forthcoming), and which in some respects more closely resembles the standard calculus of smooth functions. For example, in their system, the intermediate value theorem holds. The quasi-standard smooth functions are also arrows of the Cartesian closed category of \emph{Fermat Spaces}, which are essentially ${^\bullet}\mathbb{R}$-enriched extensions of diffeological spaces\footnote{Since the latter decades of the twentieth century the needs of theoretical physics have repeatedly challenged the limitations of classical differential geometry of smooth manifolds. Diffeological spaces (which are defined on $\mathbb{R}$) are generalizations of smooth manifolds that have been motivated by those needs. Diffeology, which is a substantial extension of differential geometry, treats such matters. See (Iglesias-Zemmour 2013) for a comprehensive discussion of the history, foundations and physical motivation of the subject.}, and the category of Fermat spaces together with its underlying calculus of quasi-standard smooth functions serves as the basis for the forthcoming development of IDG referred to in (Giordano and Wu 2016, pp. 889-890). 

Echoing the attitudes expressed by Kock and Bunge regarding SDG, Giordano does not envision IDG and its underlying calculus as providing replacements for their classical $\epsilon, \delta $-counterparts, but rather as companions that provide new insights into the structure of the smooth world while lending precision to some of the nilpotent infinitesimalist techniques that have been employed on occasion by analysts, differential geometers, physicists and engineers since the time of Fermat. 

As was alluded to above, Bertram's TDC seeks to provide a generalization of the differential aspects of classical differentiable geometry of smooth manifolds by providing a nilpotent-infinitesimalist framework that is applicable to a wide range of models of the continuum. Underlying TDC is a general nilpotent-infinitesimalist approach to differential calculus developed earlier by Bertram, Gl\"ocker and Neeb (2004).  Since integral calculus does not appear to admit such a unified theory, Bertram maintains that integration theory for the various underlying continua needs to be developed in localized settings, adding assumptions of an analytic or a topological nature, according to the strength of the results one desires. Unlike SDG and IDG, TDC is developed in a set-theoretic (rather than a category-theoretic) framework, and like IDG the varying sorts of logical apparatus that are part and parcel of SDG (topos theory) and NSDG (model theory) are absent from TDC (Bertram 2008, p. 160). Of the three above mentioned nilpotent-infinitesimalist approaches to differential geometry, TDC is the most purely algebraic.

Among the motivations underlying TDC are considerations of simplicity. As Bertram puts it:  

\begin{quote}
The suggestion to use dual numbers in differential geometry is not new--one of the earliest steps in this direction was by A. Weil [1953]; one of the most recent is [IDG due to Giordano]. However, most of the proposed constructions [including SDG] are so complicated that one is discouraged to reiterate them. But this is what makes the dual number formalism [of TDC] so useful. (Bertram 2008, p. 3) 
\end{quote}

For the sake of space, we will limit our overview of TDC to a few snapshots of the theory presented by Bertram himself, beginning with the following motivational observation, where $\mathbb{K}$ is understood to be a good topological ring in the sense defined above (Note 39). 
 
\begin{quote}
We define \emph{manifolds} and \emph{tangent bundles} in the classical way using charts and atlases...; then, intuitively, we may think of the tangent space $T_{p}M$ of [a manifold] $M$ at [point] $p$ as a ``(first order) infinitesimal neighborhood" of $p$. This idea may be formalized by writing, with respect to some fixed chart of $M$, a tangent vector at the point $p$ in the form $p+\epsilon\nu$, where $\epsilon$ is a ``very small" quantity (say, Planck's constant), thus expressing that the tangent vector is ``infinitesimally small" compared to elements of $M$ (in a chart). The property of being ``very small" can mathematically be expressed by requiring that $\epsilon^2$ is zero, compared with all ``space quantities". This suggests that, if $M$ is a manifold modelled on a  $\mathbb{K}$-module $V,$ then the tangent bundle [$TM$] should be modelled on the space $V\oplus \epsilon V=:V\times V$, which is considered a module over the  [good topological] ring $\mathbb{K}[\epsilon]=\mathbb{K}\oplus \epsilon  \mathbb{K}$ of \emph{dual numbers} over $\mathbb{K}$ (it is constructed from $\mathbb{K}$ in a similar way as the complex numbers are constructed from $\mathbb{R}$, replacing the condition $i^2 =-1$ by $\epsilon^2 =0$).  All this would be really meaningful if $TM$ were a manifold not only over $\mathbb{K}$, but also over the extended ring $\mathbb{K}[\epsilon].$ (2008, p. 2).
\end{quote}

Bertram proves that this is indeed the case. More specifically, he shows that: (i) if $M$ is a smooth manifold over $\mathbb{K}$, then $TM$ is, in a natural way, a manifold over $\mathbb{K}[\epsilon]$; and (ii) if $f : M \rightarrow N$ is smooth over  $\mathbb{K}$, then $Tf : TM \rightarrow TN$ is smooth over $\mathbb{K}[\epsilon]$, i.e. the familiar tangent functor $T$ of differential geometry can be regarded as a functor of scalar extension from $\mathbb{K}$ to $\mathbb{K}[\epsilon]$ in the category of smooth manifolds (2008: Theorem 6.2, p. 36).

\begin{quote}
Hence [says Bertram] one may use the ``dual number unit" $\epsilon$ when dealing with tangent bundles with the same right as the imaginary unit $i$ when dealing with complex manifolds. The proof...is conceptual and allows [one] to understand why dual numbers naturally appear in this context....(2008, p. 2)
\end{quote}

Indeed, according to Bertram:

\begin{quote}
in a way this result contains a new justification of ÒinfinitesimalsÓ as proposed
by A. Weil in [1953], having the advantage not to use the heavy logical or model-theoretic machinery
introduced in synthetic differential geometry.... (Bertram 2006, p. 96)
\end{quote}
\noindent
Moreover, and quite importantly, by iterations of the above result one gets much more. 
\begin{quote}
[F]or instance, if $TM$ is the scalar extension of $M$ by $\mathbb{K}[\epsilon_1]$, then the double tangent bundle $T^{2}M:=T(TM)$ is simply a scalar extension of $M$ by the ring ... $$\mathbb{K}[\epsilon_1][\epsilon_2]\cong \mathbb{K} \oplus\epsilon_1\mathbb{K} \oplus\epsilon_2\mathbb{K}\oplus\epsilon_1\epsilon_2\mathbb{K},$$ and so on for all higher tangent bundles $T^{k}M$. As a matter of fact, most of the important notions of differential geometry deal, in one way or another, with the second order tangent bundle $T^{2}M$ (e.g. Lie bracket, exterior derivative, connections) or with $T^{3}M$ (e.g. curvature). Therefore second and third order differential geometry really is the central part of all differential geometry and finding a good notation concerning second and third order tangent bundles becomes a necessity. Most textbook, if at all $TTM$ is considered, use a component notation in order to describe objects related to this bundle. In this situation, the use of \emph{different} symbols $\epsilon_1, \epsilon_{2},...$ for the infinitesimal units of the various scalar extensions is a great notational progress, combining algebraic rigour and transparency. It becomes clear that many structural features of $T^{k}M$ are simple consequences of corresponding structural features of the rings ... $\mathbb{K}[\epsilon_1, ..., \epsilon_{k}]$...[where $\epsilon_{i}^{2} =0$ and $\epsilon_{i}\epsilon_{j}=\epsilon_{j}\epsilon_{i}$ for all $1\leq i,j, \leq k$]. (Bertram, 2008, p. 3)

\end{quote}

It is worth emphasizing that while Bertram believes nilpotent infinitesimals do have an important role to play in smooth differential geometry, they need not necessarily be viewed as elements of the underlying continuum itself. As Bertram remarks: 

\begin{quote}
Some readers may wish to avoid the use of rings which are not fields, or even
to stay in the context of the real base field. In principle, all our results that do
not directly involve dual numbers can be proved in a ``purely real" way, i.e., by
interpreting $\epsilon$ just as a formal (and very useful !) label....But in the end, just
as the ``imaginary unit" $i$ got its well-deserved place in mathematics, so will the
``infinitesimal unit" $\epsilon$. (Bertram, 2008, p. 3)
\end{quote}

While we have focused our attention in this section on the nilpotent-infinitesimalist approaches to smooth differential geometry,  we would be remiss if we failed to at least mention the promise of NSDG. Formally speaking, in SDG tangent vectors are treated as infinitesimal displacements and vector fields are regarded as infinitesimal transformations. The aforementioned paper of Nowik and Katz presents a similar approach from the standpoint of nonstandard analysis.  Though still in the early stages of its development, the theory of Nowik and Katz is already promising enough that in his review of their paper for \emph{Mathematical Reviews}, H. Nishimura, a longstanding contributor to the development and application of SDG, wrote: ``A serious comparison between synthetic differential geometry and nonstandard differential geometry might presumably be intriguing" (MR3457545). Building on Nishimura's observation,  we draw this section to a close by suggesting that a serious comparison of all four infinitesimalist approaches to smooth differential geometry touched upon in this section along with an examination of their corresponding underlying treatments of smooth continua would be very intriguing indeed.

\section{invertible and nilpotent infinitesimals: afterthoughts}	

As was earlier mentioned, Robinson expressed the opinion that using the techniques of nonstandard analysis one could justify the residual uses of infinitesimals in the writings of post $\delta,\epsilon$ differential geometers and physicists. Since some of those uses were of nilpotent infinitesimals his words are in part misleading unless (as we suspect) he was suggesting that using the techniques of nonstandard analysis one can prove the classical results they established using infinitesimals, regardless of their nature.\footnote{While outside the province of this paper, we note that there is an ongoing debate concerning whether, and if so, in what respects, nonstandard analysis is a useful lens to evaluate and interpret various infinitesimalist techniques employed by Leibniz, Euler, Cauchy and other pre-$\epsilon,\delta$ analysts. The debate grew out of some early remarks of Robinson (1966, 1974, ch. x, 1967) and Lakatos (1978). For an introduction to the debate along with many references, see (Bos 1974, 2004; Laugwitz 1987, 1989; Fraser 2015; Borovik and Katz 2012; Katz and Sherry 2013; Blaszczyk, Katz and Sherry 2013; Bair \emph{et al} 2017).} By the same token, proponents of SDG often say their use of nilpotent infinitesimals lends precision to the techniques employed by classical differential geometers such as Sophus Lie and \'Elie Cartan. However, as I. Moerdijk and G. Reyes emphasize in the Preface to their monograph on models of SDG (1991, p. v): 
\begin{quote}
\emph{two} kinds of infinitesimals were used by geometers like S. Lie and E. Cartan, namely invertible infinitesimals and nilpotent ones. 
\end{quote}
\noindent
This, together with the occasional use of infinite elements by physicists, led Moerdijk and Reyes (1991, pp. 241-243, 285-286) to introduce models of the basic axioms of SDG having both sorts of infinitesimals as well as infinite elements and to issue the following proclamation (1991, p. 239):

\begin{quote}
[There are] at least two different kinds of infinitesimals that have appeared in the literature. On the one hand, there are the \emph{nilpotent}, [used] in handling ``infinitesimal" structures like jets, prolongations and connections....On the other hand, there are the \emph{invertible} infinitesimals, which together with infinitely large integers, are used to analyze such notions as limits and convergence along the lines of non-Standard Analysis, as exemplified by Robinson's book. Thus, these types of infinitesimals serve different (and complimentary purposes), and both should appear in a theory of infinitesimals worth its salt.  
\end{quote}

Thus far, however, models having both types of infinitesimals as well as infinite elements have been employed sparingly in SDG--originally by Moerdijk and Reyes in connection with the theory of distributions and the Dirac $\delta$ function (1991, pp. 320-337) and later in connection with differential equations (Kennison 1999). Moreover, their use runs contrary to the stated Archimedean nature of SDT, and according to Kock, ``the[ir] value...for synthetic reasoning in geometry is not evident" (1993, p. 354). Invertible infinitesimals along with their multiplicative inverses can also be added to IDG (Giordano 2001, p. 77), though thus far they have not. And, while it has been suggested that nilpotent infinitesimals could be added to nonstandard analysis (Giordano and Wu 2016, p. 896), that would require substantial modifications of the theory which would no longer give rise to an elementary extension of the reals and it is by no means clear the gains would outweigh the losses. Collaterally, as Robinson appears to have maintained, and further development of the work of Nowik, Katz and others may eventually show, it may be possible to develop an entirely adequate infinitesimalist approach to the smooth aspects of differential geometry without nilpotent infinitesimals at all. 

On the other hand, like the models of Moerdijk and Reyes, systems containing both types of infinitesimals as well as infinite numbers fall within the framework of TDC (Bertram 2008, p. viii); and so, at least from the standpoint of the purely differential aspects of smooth differential geometry as characterized by TDC the transition to such a system is seamless. Moreover, such systems have found application in algebra and geometry since Zemmer (1953, p. 177; also see, Fuchs 1963, pp. 108-109) and Klingenberg (1954a) pioneered their use well over half a century ago. 

The question whether a theory of infinitesimals worth its salt requires both types of infinitesimals, as Moerdijk and Reyes maintain, is a fascinating one, as is the corresponding question regarding an infinitesimalist theory of continua built thereon. While support for Moerdijk and Reyes's contention has thus far been muted, whether with time this will change, only time will tell.

\section{concluding remarks} Since the discovery that a diagonal of a square is incommensurable with its sides, the question of the possibility of bridging the gap between the domains of discreteness and of continuity, or between arithmetic and geometry, has been a central problem in the foundations of mathematics. Cantor and Dedekind of course believed they had bridged the gap with the creation of their arithmetico-set-theoretic continuum of real numbers, and it remains a central tenet of standard mathematical philosophy that indeed they had. Nevertheless, Cantor was overly optimistic when he suggested that his theory of the continuum, unlike that of the ancients, had ``been thought out ... with the clarity and completeness ... required to exclude the possibility of different opinions among [its] posterity" (Cantor 1883/1996, p. 903). After all, whereas Cantor and Dedekind had succeeded in replacing the vague ancient conception with a clear and precise arithmetico-set-theoretic conception that proved itself adequate for the needs of analysis, differential geometry and the empirical sciences of their day, they could neither free their theory of its logical, theoretical, and philosophical presuppositions, nor preclude the possibility that other adequate conceptual schemes, each self-consistent, could be devised offering alternative visions of the continuum.

However, it was critiques of the former and/or the realization of the logical possibility of the latter that has given rise to a host of non-Cantor-Dedekindean conceptions, including the recent infinitesimalist conceptions canvassed above. In the latter cases, with the exception of the surreals, the architects were motivated by the belief, or at least the hope, that their respective theories are, or with time would be, adequate for the needs of analysis or portions of differential geometry, and the empirical sciences served thereby. Like their late nineteenth- and early twentieth-century non-Archimedean geometric forerunners, nonstandard analysis and the infinitesimalist approaches to portions of differential geometry have drawn attention to the possibility of physical continua whose logical cogency, let alone physical possibility, had long been in doubt. Whether empirical science will require such a theory, as some, like Fenstad (1987, 1988), already contend, and others, like Veronese (1909/1994, p. 180) and Hahn (1933/1980 p. 100), would not rule out, only time will tell. On the other hand, while showing no sign of displacing the Cantor-Dedekind theory and the collateral theories based thereon, the infinitesimalist approaches have performed, and continue to perform, important logical and philosophical service. Nonstandard analysis has also had substantial success in shedding important light on, and establishing significant new results in, various areas of analysis, theoretical physics and economics, and SDG/SDT has provided important insights into the relation between algebraic geometry, differential geometry and smooth spaces more generally. Along with their theoretical accomplishments, the practitioners of nonstandard analysis and SDG/SDT maintain that their approaches have heuristic and intuitive advantages over their standard counterparts. The authors of IDG and TDC likewise maintain such advantages for their approaches, but without, they add, either the need for model-theoretic apparatus or the abandonment of classical logic. Not eschewing model-theoretic apparatus, NSDG seeks an extension of the range of application of nonstandard analysis, which, in principle, could provide an infinitesimalist approach to all of analysis and the smooth portions of differential geometry, something that is beyond the scope of SDG/SDT, IDG and TDC. 

Of course, whether nonstandard analysis, SDG/SDT or any of the other just-said infinitesimalist approaches, together with its corresponding conception of the continuum, will eventually displace its standard counterpart or even become a widely employed implement in the toolkits of mathematicians remains to be seen. Also remaining to be seen is whether the s-hierarchical ordered field of surreal numbers will come to be widely regarded as an absolute arithmetic continuum modulo NBG. However, regardless of how these questions are ultimately answered, one cannot help but wonder whether or not F. W. Lawvere was prophetic when he recently maintained: 
``Contrary to common opinion, the question ``what is the continuum?" does not have a final answer, the immortal work of Dedekind notwithstanding" (2011, p. 249).

\section*{references}

\begin{footnotesize}

\smallskip
\noindent
Alling, N. (1987): \emph{Foundations of Analysis Over Surreal Number Fields}, North-Holland Publishing Company, Amsterdam.

\smallskip
\noindent
Alling, N. and Ehrlich, P. (1986): \emph{An alternative construction of Conway's surreal numbers}, Comptes Rendus Math\'ematiques de l'Acad\'emie des Sciences, La Soci\'et\'e Royale du Canada VIII, pp. 241-246.

\smallskip
\noindent
Alling, N. and Ehrlich, P. (1987): Sections 4.02 and 4.03 of (Alling 1987).

\smallskip
\noindent
Almeida, R., Neves, V. and Stroyan, K. (2014): \emph{Infinitesimal differential geometry: cusps and envelopes}, Differential Geometry--Dynamical Systems 16, pp. 1-13.

\smallskip
\noindent
Anderson, R. M. (1976): \emph{A non-standard representation for Brownian motion and It\^o integration},
Israel Journal of Mathematics 25, pp. 15-46. 

\smallskip
\noindent
Artin, E. and Schreier, O. (1926): \emph{Algebraische Konstruktion reeller K\"orper}, Abhandlungen aus dem Mathematischen Seminar der Universit\"at Hamburg 5, pp. 85-99. Translated as \emph{Algebraic construction of real fields}, in \textbf{\emph{Exposition by Emil Artin: A Selection}}, edited by  M. R. Rosen, American Mathematical Society, London Mathematical Society, Providence, R.I., 2007, pp. 273-283.

\smallskip
\noindent
Aschenbrenner, M., van den Dries, L. and van der Hoeven, J. (2017): \emph{Asymptotic Differential Algebra and Model Theory of Transseries}, Annals of Mathematics Studies, 195, Princeton University Press, Princeton, NJ.

\smallskip
\noindent
Aschenbrenner, M., van den Dries, L. and van der Hoeven, J. (2019): \emph{Numbers, germs and transseries}, in \textbf{\emph{Proceedings of the International Congress of Mathematicians, Rio De Janeiro, 2018, Volume 2}}, edited by B. Sirakov, P. N. de Souza and M. Viana, World Scientific Publishing Company, NJ, pp. 19-43.  

\smallskip
\noindent
Aschenbrenner, M., van den Dries, L. and van der Hoeven, J. (forthcoming): \emph{Surreal numbers as a universal $H$-field}, Journal of the European Mathematical Society.

\smallskip
\noindent
Bachmann, F. (1971): \emph{Hjelmslev planes}, in  \textbf{\emph{Atti del Convegno di Geometria Combinatoria e sue Applicazioni (University of Perugia, Perugia, 1970)}}, Istituto di Matematica, University of Perugia, Perugia, Italy, pp. 43-56.   

\smallskip
\noindent
Bachmann, F. (1973): \emph{Geometry of reflections}, in  \textbf{\emph{Atti del Convegno internazionale sul tema Storia, pedagogia e filosofia della scienza. A celebrazione del centenario della nascita di Federigo Enriques (Pisa, Bologna e Roma, 1971)}}, Accademia nazionale dei Lincei, Roma, pp. 101-108.

\smallskip
\noindent
Bachmann, F. (1973a): \emph{Aufbau der Geometrie aus dem Spiegelungsbegriff, Zweite erg\"anzte Auflage}, 
Die Grundlehren der mathematischen Wissenschaften, Band 96, Springer-Verlag, Berlin-New York.

\smallskip
\noindent
Bachmann, F. (1989): \emph{Ebene Spiegelungsgeometrie. Eine Vorlesung \"uber Hjelmslev-Gruppen.  With a foreword by M. G\"otzky  and H. Wolff}, Bibliographisches Institut, Mannheim, 1989.

\smallskip
\noindent
Baer, R. (1927): \emph{\"Uber nicht-archimedisch geordnete K\"orper}, Sitzungsberichte der Heidelberger Akademie
der Wissenschaften 8, pp. 3-13.

\smallskip
\noindent
Baer, R. (1929): \emph{Zur Topologie der Gruppen}, Journal f\"ur die Reine und Angewandte Mathematik 160, pp. 208-226.

\smallskip
\noindent
Baer, R. (1970): \emph{Dichte Archimedizit\"at und Starrheit geordneter K\"orper}, Mathematische Annalen
188, pp. 165-205.

\smallskip
\noindent
Bair, J., B\l aszczyk, P., Ely, R., Henry, V., Kanovei, V., Katz, K., Katz, M., Kutateladze, S., McGaffey, T., Reeder, P., Schaps, D., Sherry, D. and Shnider, S.
(2017): \emph{Interpreting the infinitesimal mathematics of Leibniz and Euler}, 
Journal of General Philosophy of Science 48, pp. 195-238.

\smallskip
\noindent
Banach, S. (1931): \emph{\"Uber die Baire'sche Kategorie gewisser Funktionenmengen}, Studia Mathematica 3, pp. 174--179.

\smallskip
\noindent
Bauer, A. (2017): \emph{Five stages of accepting constructive mathematics}, Bulletin of the American Mathematical Society (N. S.) 54, pp. 481--498. 

\smallskip
\noindent
Belair, L. (1981): \emph{Calcul infinit\'esimal en g\'eom\'etrie diff\'erentielle synth\'etique}, Masters Thesis, Universiti\'e de Montr\'eal.

\smallskip
\noindent
Belair, L. and Reyes, G. E. (1982):
\emph{Calcul infinit\'esimal en g\'eom\'etrie diff\'erentielle synth\'etique}, Sidney Category Seminar Reports.

\smallskip
\noindent
Bell, J. L. (1995): \emph{Infinitesimals and the continuum}, The Mathematical Intelligencer 17 (2), pp. 55-57.

\smallskip
\noindent
Bell, J. L. (2001): \emph{The Continuum in smooth infinitesimal analysis}, in \textbf{\emph{Reuniting the Antipodes-Constructive and Nonstandard Views of the Continuum}}, edited by P. Schuster, U. Berger and H. Osswald, Dordrecht, Holland.

\smallskip
\noindent
Bell, J. L. (2005): \emph{The Continuous and Infinitesimal in Mathematics and Philosophy}, Polimetrica, International Scientific Publisher, Monza-Milano.

\smallskip
\noindent
Bell, J. L. (2008): \emph{A Primer of Infinitesimal Analysis, 2nd Edition}, Cambridge University Press, Cambridge.

\smallskip
\noindent
Bell, J. L. (2009): \emph{Cohesiveness}, Intellectica 51, pp. 1-24.

\smallskip
\noindent
Benci, V.,  Forti, M.  and Di Nasso, M. (2006): \emph{The eightfold path to nonstandard analysis}, in
\textbf{\emph{Nonstandard methods and applications in mathematics}}, Lecture Notes
in Logic  25, Association of Symbolic Logic, La Jolla, CA., pp. 3-44.

\smallskip
\noindent
Berarducci, A. and Mantova, V. (2018): \emph{Surreal numbers, derivations and transseries}, Journal of the European Mathematical Society 20, pp. 339-390.

\smallskip
\noindent
Berarducci, A. and Mantova, V. (forthcoming): \emph{Transseries as germs of surreal functions}, Transactions of the American Mathematical Society.

\smallskip
\noindent
Bertram, W. (2006): \emph{Differential geometry over general base fields and rings}, in \textbf{\emph{Modern Trends in Geometry and Topology. Proceedings of the 7th International Workshop on Differential Geometry and its Applications held in Deva, September 5-11, 2005}}, edited by D. Andrica, P. A. Blaga and S. Moroianu, Cluj University Press, Cluj-Napoca, pp. 95-101.

\smallskip
\noindent
Bertram, W. (2008):  \emph{Differential Geometry, Lie Groups and Symmetric Spaces over General Base Fields and Rings}, Memoirs of the American Mathematical Society  192, No. 900.  

\smallskip
\noindent
Bertram, W. (2011):  \emph{Calcul diff\'erentiel topologique \'el\'ementaire}, Calvage and Mounet, Paris.

\smallskip
\noindent
Bertram, W., Gl\"ockner, H. and Neeb, K-H. (2004):
\emph{Differential calculus over general base fields and rings}, 
Expositiones Mathematicae  22, pp. 213--282. 

\smallskip
\noindent
Bettazzi, R. (1890): \emph{Teoria Delle Grandezze}, Pisa.

\smallskip
\noindent
B\l aszczyk, P., Katz, M. and Sherry, D. (2013): \emph{Ten misconceptions from the history of analysis and their debunking}, 
Foundations of  Science 18, pp. 43-74. 

\smallskip
\noindent
Borovik, A. and Katz, M. (2012): \emph{Who gave you the Cauchy-Weierstrass tale? The dual history of rigorous calculus}, 
Foundations of Science 17, pp. 245-276.

\smallskip
\noindent
Borovik, A., Jin, R. and Katz, M. (2012): \emph{An integer construction of infinitesimals: toward a theory of Eudoxus hyperreals}, Notre Dame Journal of Formal Logic 53, pp. 557-570.

\smallskip
\noindent
Bos, H. (1974); \emph{Differentials, higher-order differentials and the derivative in the Leibnizian calculus}, Archive for History of Exact Sciences 14, pp. 1-90. 

\smallskip
\noindent
Bos, H. (2004); \emph{Philosophical challenges from history of mathematics. New trends in the history and philosophy of mathematics}, University of Southern Denmark Studies in Philosophy 19, pp. 51-66.

\smallskip
\noindent
Bunge, M. (1984): \emph{Toposes in logic and logic in toposes}, Topoi 3, pp. 13-22.

\smallskip
\noindent
Bunge, M. and Dubuc, E. J. (1986): \emph{Archimedian local $C^\infty$ rings and models of SDG}, Cahiers de Topologie et G\'eom\'etrie Diff\'erentielle Cat\'egoriques  27, pp. 3-22.

\smallskip
\noindent
Bunge, M. and Dubuc, E. J. (1987): \emph{Local concepts in synthetic differential geometry and germ representability}, in \textbf{\emph{Mathematical logic and theoretical computer science (College Park, Md., 1984--1985)}}, edited by D. Kueker, E. Lopez-Escobar and C. Smith, Lecture Notes in Pure and Applied Mathematics 106, Dekker, New York, pp. 93-159. 

\smallskip
\noindent
Bunge, M. and Gago, F. (1988):
\emph{Synthetic aspects of $C^{\infty}$-mappings. II. Mather's theorem for infinitesimally represented germs},
Journal of Pure and Applied Algebra 55, pp. 213-250.

\smallskip
\noindent
Bunge, M., Gago, F. and San Luis, A. M. (2018): \emph{Synthetic Differential Topology}, Cambridge University Press, Cambridge.

\begin{sloppypar}
\smallskip
\noindent
Cantor, G. (1872): \emph{\"Uber die Ausdehnung eines Satzes der Theorie der trigonometrischen Reihen}, Mathematische Annalen 5, pp. 122-132. Reprinted in \textbf{\emph{Georg Cantor Gesammelte Abhandlungen mathematischen und philosophischen Inhalts}}, edited by E. Zermelo, J. Springer, Berlin, 1932, pp. 92-102.  
\end{sloppypar}

\smallskip
\noindent
Cantor, G. (1883/1996): \emph{Grundlagen einer allgemeinen Mannigfaltigkeitslehre. Ein mathematisch-philosophischer Versuch in der Lehre des Unendlichen}. Translated by W. B. Ewald as \emph{Foundations of a general theory of manifolds: a mathematico-philosophical investigation into the theory of the infinite}, in \textbf{\emph{From Kant to Hilbert: A Source Book in the Foundations of Mathematics, Volume II}}, edited by W. B. Ewald,  Clarendon Press, Oxford, 1996, pp. 878-920.

\smallskip
\noindent
Cantor, G. (1895): \emph{Beitr\"age zur Begr\"undung der transfiniten Mengenlehre, Part I}, Mathematische Annalen 46, pp. 481-512. English translation in \textbf{\emph{Contributions to the Founding of Transfinite Numbers, Translated, and Provided
with an Introduction and Notes by Philip E. B. Jourdain}}, 1915. Reprinted by Dover Publications,
Inc., New York, 1955, pp. 85-136.

\smallskip
\noindent
Chang, C. C. and  Keisler, H. J. (1990): \emph{Model Theory, Third Edition}, North-Holland Publishing Company, Amsterdam, New York. Reprinted Dover Publications, Inc., Mineola, New York, 2012.

\smallskip
\noindent
Clifford, W. (1873): \emph{Preliminary sketch of bi-quaternions}, Proceedings of the London Mathematical Society 4, pp. 381--95.

\smallskip
\noindent
Cohen, L. W. and Goffman, C. (1949): \emph{The topology of ordered abelian groups}, Transactions of the American Mathematical Society 67, pp. 310-319.

\smallskip
\noindent
Connes, A. (1994): \emph{Noncommutative Geometry}, Academic Press, London.

\smallskip
\noindent
Connes, A. (1998): \emph{Noncommutative differential geometry and the structure of space-time}, in \textbf{\emph{The Geometric Universe: Science, Geometry, and the Work of Roger Penrose}}, edited by S. A. Huggett, L. J. Mason, K. P. Tod, S. T. Tsou and N. M. J. Woodhouse, Oxford University Press, Oxford, pp. 48-80.	

\smallskip
\noindent
Connes, A. (2006): \emph{On the foundations of noncommutative geometry}, in \textbf{\emph{The Unity of Mathematics, In Honor of the Ninetieth Birthday of I. M. Gelfand}}, edited by  P. Etingof, V. S. Retakh and I. M. Singer, Birkh\"auser, Boston, pp. 173-204.

\smallskip
\noindent
 Conway, J. H. (1976): \emph{On Numbers and Games}, Academic Press, London. \emph{Second Edition}, A K Peters, Ltd., Natick, Massachusetts (2001).
 
\smallskip
\noindent
Costin, O. (2009): \emph{Asymptotics and Borel Summability}, Chapmann \& Hall, New York.

 \smallskip
\noindent
Costin, O. and Ehrlich, P. (2016): \emph{Integration on the surreals: a conjecture of Conway, Kruskal and Norton}, in \textbf{\emph{Mini-workshop: surreal numbers, surreal analysis, Hahn fields and derivations. Abstracts from the mini-workshop held December 18-23, 2016. Organized by A. Berarducci, P. Ehrlich and S. Kuhlmann}}, \emph{Oberwolfach Reports}  13, No. 4, Gerhard Huisken, Editor in Chief, European Mathematical Society Publishing House, Z\"urich, Switzerland, pp. 3364-3365.

 \smallskip
\noindent
Costin, O., Ehrlich, P. and Friedman, H. (24 Aug 2015): \emph{Integration on the surreals: a conjecture of Conway, Kruskal and Norton}, preprint, arXiv:1334466.

 \smallskip
\noindent
Dahn, B. G\"oring and P. (1987): \emph{Notes on exponential-logarithmic terms}, Fundamenta Mathematica 127, pp. 45-50.

\smallskip
\noindent
Dedekind, R. (1872): \emph{Stetigkeit und irrationale Zahlen}, Braunschweig, Vieweg. Translated by W. W. Beman (with corrections by W. B. Ewald) as \emph{Continuity and Irrational Numbers}, in \textbf{\emph{From Kant to Hilbert: A Source Book in the Foundations of Mathematics, Volume II}}, edited by W. B. Ewald, Oxford: Clarendon Press, 1996.

\smallskip
\noindent
Dembowski, P. (1968): \emph{Finite Geometries}, Springer-Verlag, Berlin, New York.

\smallskip
\noindent
Drake, D. A. and Jungnickel, D. (1985): \emph{Finite Hejemslev planes and Klingenberg epimorphisms}, in \textbf{\emph{Rings and Geometry}}, edited by R. Kaya, P. Plaumann and K. Strambach, D. Reidel Publishing Company, Dordrecht, pp. 153-231.

\smallskip
\noindent
du Bois-Reymond, P. (1870-71): \emph{Sur la grandeur relative des infinis des fonctions}, Annali di matematica pura ed applicata 4, pp. 338-353.

\smallskip
\noindent
du Bois-Reymond, P. (1875): \emph{Ueber asymptotische Werthe, infinit\"are Approximationen und infinit\"are Aufl\"osung von Gleichungen}, Mathematische Annalen 8, pp. 363-414; (\emph{Nachtr\"age zur Abhandlung: ueber asymptotische Werthe etc.}), pp. 574-576).

\smallskip
\noindent
du Bois-Reymond, P. (1877): \emph{Ueber die Paradoxen des Infinit\"arcalc\"uls}, Mathematische Annalen 11, pp. 149-167.

\smallskip
\noindent
du Bois-Reymond, P. (1882): \emph{Die allgemeine Functionentheorie I}, Verlag der H. Laupp'schen Buchhandlung, T\"ubingen.

\smallskip
\noindent
Dubuc, E. J. (1979): \emph{Sur les mod\`eles de la g\'eom\'etrie diff\'erentielle synth\'etique}, Cahiers de Topologie et G\'eom\'etrie Diff\'erentielle Cat\'egoriques 20, pp. 231-279.

\smallskip
\noindent
Dubuc, E. J. (1981): \emph{$C^{\infty}$-schemes}, American Journal of Mathematics 103, pp. 683--690.

\smallskip
\noindent
Dubuc, E. J. (1981a): \emph{Open covers and infinitary operations in $C^{\infty}$-rings}, Cahiers de Topologie et G\'eom\'etrie Diff\'erentielle Cat\'egoriques 22, pp. 287-300.

\smallskip
\noindent
\'Ecalle, J. (1981-1985): \emph{Fonctions Resurgentes I-III}, Publications Math\'ematiques d'Orsay, 81, Universit\'e de Paris-Sud, Depart\'ement de Math\'ematique, Orsay.

\smallskip
\noindent
 \'Ecalle, J. (1992): \emph{Introduction aux fonctions analysables et preuve constructive de la conjecture de Dulac},  Actualit\'es Math\'ematiques, Hermann, Paris.

\smallskip
\noindent
 \'Ecalle, J. (1993): \emph{Six Lectures on Transseries, Analysable Functions and the Constructive Proof of Dulac's Conjecture, Bifurcations and Periodic Orbits of Vector Fields}, NATO Advanced Science Institutes Series C: Mathematical and Physical Sciences 408, Kluwer, Dordrecht.

\smallskip
\noindent
Edgar, G.A. (2010): \emph{Transseries for beginners}, Real Analysis Exchange 35, pp. 253-310.

\smallskip
\noindent
Ehresmann, C. (1951): \emph{Les prolongements d'une vari\'et\'e diff\'erentiable. I-III}, Comptes Rendus Math\'ematique, Acad\'emie des Sciences, Paris  233, pp. 598--600, 777-779, 1081-1083.

\smallskip
\noindent
Ehrlich, P. (1987): \emph{The absolute arithmetic and geometric continua}, in \textbf{\emph{PSA 1986, Volume 2}}, edited by A. Fine and P. Machamer, Philosophy of Science Association, Lansing, MI, pp. 237-247.

\smallskip
\noindent
Ehrlich, P. (1988): \emph{An alternative construction of Conway's ordered field No}, Algebra Universalis 25, pp. 7-16; errata, ibid. 25 (1988), p. 233.

\smallskip
\noindent
Ehrlich, P. (1989): \emph{Absolutely saturated models}, Fundamenta Mathematica 133, pp. 39-46.

\smallskip
\noindent
Ehrlich, P. (1992): \emph{Universally extending arithmetic continua}, in  \textbf{\emph{Le Labyrinthe du Continu, Colloque de Cerisy}}, edited by H. Sinaceur and J.-M. Salanskis, Springer-Verlag France, Paris, pp. 168-177.

\smallskip
\noindent
Ehrlich, P. (1994): \emph{All number great and small}, in \textbf{\emph{Real Numbers, Generalizations of the Reals, and Theories of Continua}}, edited by P. Ehrlich, Kluwer Academic Publishers, Dordrecht, pp. 239-258.

\smallskip
\noindent
Ehrlich, P. (1995): 
\emph{Hahn's \"Uber die nichtarchimedischen Gr\"ossensysteme  and the origins of the modern theory of magnitudes and numbers to measure them}, in \textbf{\emph{From Dedekind to G\"odel: Essays on the Development of the Foundations of Mathematics}}, edited by J. Hintikka, Kluwer Academic Publishers, Dordrecht, pp. 165-213.

\smallskip
\noindent
Ehrlich, P. (1997): \emph{Dedekind cuts of Archimedean complete ordered abelian groups}, Algebra Universalis  37, pp. 223-234.

\smallskip
\noindent
Ehrlich, P. (1997a): 
\emph{From completeness to Archimedean completeness: an essay in the foundations of Euclidean geometry}, in \textbf{\emph{A Symposium on David Hilbert}}, edited by A. Tauber and A. Kanamori, Synthese 110, pp. 57-76.

\smallskip
\noindent
Ehrlich, P. (2001): \emph{Number systems with simplicity hierarchies: a generalization of Conway's theory of surreal numbers}, The Journal of Symbolic Logic 66, pp. 1231-1258.

\smallskip
\noindent
Ehrlich, P. (2002): \emph{Surreal numbers: an alternative construction (Abstract)}, The Bulletin of Symbolic Logic 8, p. 448.

\smallskip
\noindent
Ehrlich, P. (2005): \emph{Continuity}, in \textbf{\emph{Encyclopedia of Philosophy, Second Edition, Volume 2}}, D. M. Borchart, Editor in Chief, Macmillan Reference USA, pp. 489-518.

\smallskip
\noindent
Ehrlich, P. (2006):
\emph{The Rise non-Archimedean mathematics and the roots of a misconception I: the emergence of non-Archimedean systems of magnitudes}, Archive for History of Exact Sciences 60, pp. 1-121.

\smallskip
\noindent
Ehrlich, P. (2007): \emph{Review of} 
 \textbf{\emph{The Continuous and the Infinitesimal in Mathematics and Philosophy}} \emph{by John L. Bell}, The Bulletin of Symbolic Logic 13, pp. 361-363.

\smallskip
\noindent
Ehrlich, P. (2010): \emph{The absolute arithmetic continuum and its Peircean counterpart}, in \textbf{\emph{New Essays on Peirce's Mathematical Philosophy}}, edited by M. Moore, Open Court Press, pp. 235-282.

\smallskip
\noindent
Ehrlich, P. (2011): \emph{Conway names, the simplicity hierarchy and the surreal number tree}, Journal of Logic and Analysis 3, no. 1, pp. 1-26.

\smallskip
\noindent
Ehrlich, P. (2012): \emph{The absolute arithmetic continuum and the unification of all numbers great and small}, The Bulletin of Symbolic Logic 18, pp. 1-45.

\smallskip
\noindent
Ehrlich, P. and Kaplan, E. (2018): \emph{Number systems with simplicity hierarchies: a generalization of Conway's theory of surreal numbers II}, The Journal of Symbolic Logic 83, pp. 617-633.

\smallskip
\noindent
Euler, L. (1778): \emph{De infinities infinitis gradibus tam infinite magnorum quam infinite parvorum}, Acta Academiae Scientiarum Petroplitanae (II) I, 1780, pp. 102-118. Reprinted in \textbf{\emph{Opera Omnia (I) XV}}, pp. 298-313.

\smallskip
\noindent
Fenstad, J. E. (1987): \emph{The discrete and the continuous in mathematics and the natural sciences}, in \textbf{\emph{L'Infinito Nella Scienza}}, edited by G. T. di Francia, Istituto Della Enciclopedia Italiana, Fondata Di G. Treccani, Roma.

\smallskip
\noindent
Fenstad, J. E. (1988):  \emph{Infinities in mathematics and the natural sciences},Ó in \textbf{\emph{Methods and Applications of Mathematical Logic}}, Contemporary Mathematics Volume 69, edited by W. A. Carnielli and L. P. de Alcantara, American Mathematical Society, Providence, RI.

\smallskip
\noindent
Fisher, G. (1981): \emph{The infinite and infinitesimal quantities of
du Bois-Reymond and their reception}, Archive for History of Exact Sciences 24, pp. 101-164.

\smallskip
\noindent
Fletcher, P., Hrbacek, K., Kanovei, V., Katz, M. G., Lobry, C. and Sanders, S. (2017): \emph{Approaches to analysis with infinitesimals following Robinson, Nelson and others}, Real Analysis Exchange 42, pp. 193-251.

\smallskip
\noindent
Fornasiero, A. (2016): \emph{Initial embeddings in the surreals}, in \textbf{\emph{Mini-workshop: surreal numbers, surreal analysis, Hahn fields and derivations. Abstracts from the mini-workshop held December 18-23, 2016. Organized by A. Berarducci, P. Ehrlich and S. Kuhlmann}}, \emph{Oberwolfach Reports} 13, No. 4, Gerhard Huisken, Editor in Chief, European Mathematical Society Publishing House, Z\"urich, pp. 3325-3326.

\smallskip
\noindent
Fraser, C. (2015): \emph{Nonstandard analysis, infinitesimals, and the history of calculus}, in \textbf{\emph{A Delicate Balance: Global Perspectives on Innovation and Tradition in the History of Mathematics: A Festschrift in Honor of Joseph W. Dauben}}, edited by D. Rowe and W-S. Horng, Birkh\"auser/Springer, Cham, pp. 25-49. 

\smallskip
\noindent
Fuchs, L. (1963): \emph{Partialy Ordered Algebraic Systems}, Pergamon Press, Oxford. Reprinted by Dover Publications, Inc., New York, 2011.

\smallskip
\noindent
Gamboa, J. M. (1987): \emph{ Some new results on ordered fields}, Journal of Algebra 110, pp. 1-12.

\smallskip
\noindent
Gesztelyi, E. (1958): \emph{Eine neue Begr\"undung der Differentialrechnung}, Matematikai Lapok 9, pp. 91-114.

\smallskip
\noindent
Giordano, P. (2001): \emph{Nilpotent infinitesimals and synthetic differential geometry in classical logic}, in \textbf{\emph{Reuniting the Antipodes--Constructive and Nonstandard Views of the Continuum}}, edited by P. Schuster, U. Berger, and H. Osswald, Kluwer Academic Publishers, The Netherlands,  pp. 75-92.

\smallskip
\noindent
Giordano, P. (2010): \emph{Infinitesimals without logic}, Russian Journal of Mathematical Physics 17, pp. 159-191.

\smallskip
\noindent
Giordano, P. (2010a): \emph{The ring of Fermat reals}, Advances in Mathematics 225, pp. 2050-2075.

\smallskip
\noindent
Giordano, P. (2011): \emph{Fermat-Reyes method in the ring of Fermat reals}, Advances in Mathematics 228, pp. 862-893.

\smallskip
\noindent
Giordano, P. (2011a): \emph{Infinite dimensional spaces and Cartesian closedness}, Journal of Mathematical Physics, Analysis and Geometry 7, pp. 225--284.

\smallskip
\noindent
Giordano, P. and Wu, E. (2016): \emph{Calculus in the ring of Fermat reals, Part I: integral calculus}, Advances in Mathematics 289, pp. 888-927.

\smallskip
\noindent
Giordano, P. and Wu, E. (forthcoming): \emph{Calculus in the ring of Fermat reals, Part II: differential calculus}.

\smallskip
\noindent
Giordano, P. and Kunzinger, M. (2013): \emph{Topological and algebraic structures on the ring of Fermat reals}, Israel Journal of Mathematics 193, pp. 459-505.

\smallskip
\noindent
Goldblatt, R. (1998): \emph{Lectures on the Hyperreals: An Introduction to Nonstandard Analysis}, Springer, New York.

\smallskip
\noindent
Gonshor, H. (1986): \emph{An Introduction to the Theory of Surreal Numbers}, Cambridge University Press, Cambridge. 

\smallskip
\noindent
Grothendieck, A. (in collaboration with J. A. Dieudonn\'e) (1960): \emph{\'El\'ements de G\'eom\'etrie Alg\'ebrique I. Le Langage des Sch\'emas},
Institut des Hautes \'Etudes Scientifiques, Publications Math\'ematiques  No. 4. 

\smallskip
\noindent
Grothendieck, A. and Dieudonn\'e, J. A. (1971): \emph{\'El\'ements de G\'eom\'etrie Alg\'ebrique I}, 
Grundlehren der Mathematischen Wissenschaften 166, Springer-Verlag, Berlin.

\smallskip
\noindent
Gr\"unwald, J. (1906): \emph{\"Uber duale Zahlen und ihre Anwendung in der Geometrie}, Monatshefte f\"ur Mathematik 17, pp. 8-136.

\smallskip
\noindent
Hahn, H. (1907): \emph{\"Uber die nichtarchimedischen Gr\"ossensysteme}, Sitzungsberichte der Kaiserlichen  Akademie der Wissenschaften, Wien, Mathematisch-Naturwissenschaftliche Klasse 116 (Abteilung  IIa), pp. 601-655.

\smallskip
\noindent
Hahn, H. (1933/1980): \emph{The crisis in intuition}, in \textbf{\emph{Empiricism, Logic and Mathematics: Philosophical Papers by Hans Hahn}}, edited by B. F. McGuinness, D. Reidel Publishing Company, Dordrecht, Holland, 1980, pp. 73-102. First published in \textbf{\emph{Krise und Neuaufbau in den exakten Wissenschaften}}, F\"unf Wiener Vortr\"age, Leipzig and Vienna, 1933.

\smallskip
\noindent
Hardy, G. H. (1910): \emph{Orders of Infinity, The ÒInfinit\"arcalc\"ulÓ of Paul Du Bois-Reymond}, Cambridge University Press, Cambridge. Reprinted  Hafner, New York, 1971.

\smallskip
\noindent
Hardy, G. H. (1912): \emph{Properties of logarithmico-exponential functions}, Proceedings of the
London Mathematical Society, Second Series 10, pp. 54-90.

\smallskip
\noindent
Hauschild, K. (1966): \emph{\"Uber die Konstruktionon von Erweiterungsk\"orpern zu nichtarchimedisch angeordneten K\"orpern mit Hilfe von H\"olderschen Schnitten}, Wissenschaftliche Zeitschrift der Humboldt-Universit\"at zu Berlin, Mathematisch-Naturwissenschaftliche Reihe 15, pp. 685-688.

\smallskip
\noindent
Hausdorff, F. (1907): \emph{Untersuchungen \"uber Ordungtypen IV, V}, Berichte \"uber die Verhandlungen der K\"oniglich S\"achsischen Gesellschaft der Wissenschaften zu Leipzig, Matematisch-Physische Klasse 59, pp. 84-159. Translated as \emph{Investigations into order types IV, V}, in (Plotkin 2005), pp. 113-171.

\smallskip
\noindent
Hausdorff, F. (1909): \emph{Die Graduierung nach dem Endverlauf}, Abhandlungen der K\"oniglich S\"achsischen Gesellschaft der Wissenschaften zu Leipzig, Matematisch-Physische Klasse 31, pp. 295-335. Translated as \emph{Graduation by final behavior} in (Plotkin 2005), pp. 271-301.

\smallskip
\noindent
Hellman, G. (2006): \emph{Mathematical pluralism: the case of smooth infinitesimal analysis}, Journal of Philosophical Logic 35, pp. 621-651.

\smallskip
\noindent
Henson, C. W. (1974): \emph{The isomorphism property in nonstandard analysis and its use in the theory of Banach space}, The Journal of Symbolic Logic 39, pp. 717-771.

\smallskip
\noindent
Henson, C. W. (1979): \emph{Unbounded Loeb measures}, Proceedings of the American Mathematical Society 74, pp. 143-150.

\smallskip
\noindent
Henson, C. W. and Keisler, H. J. (1986): \emph{On the strength of nonstandard analysis}, The Journal of Symbolic Logic 51, pp. 377-386.

\smallskip
\noindent
Hewitt, E. (1948): \emph{Rings of real-valued continuous functions. I}, Transactions of the American Mathematical Society 64, pp. 44-99.

\smallskip
\noindent
Hilbert, D. (1899): \emph{Grundlagen der Geometrie}, Teubner, Leipzig.

\smallskip
\noindent
Hilbert, D. (1871): \emph{Foundations of Geometry 10th edition} of (Hilbert 1899), translated by L. Unger, Open Court, LaSalle, IL. 

\smallskip
\noindent
Hjelmslev, J. (1916): \emph{Die Geometrie der Wirklichkeit}, Acta Mathematica 40, pp. 3-66.

\smallskip
\noindent
Hjelmslev, J. (1923): \emph{Die nat\"urliche Geometrie},  Abhandlungen aus dem Mathematischen Seminar der Universit\"at Hamburg  2, pp. 1-36. 

\smallskip
\noindent
Hjelmslev, J. (1929--1949): \emph{Einleitung in die allgemeine Kongruenzlehre I,II,III,IV,V, VI}, Det Kongelige Danske Videnskabernes Selskab Matematisk-Fysiske Meddelelser 8, Nr. 11 (1929); 10, Nr. 1 (1929); 19, Nr. 12 (1942); 22, Nr. 6 (1945);  22, Nr. 13 (1945); 25, Nr. 10 (1949).

\smallskip
\noindent
H\"older, O. (1901): \emph{Die Quantit\"at und die Lehre vom Mass},Ó Berichte \"uber die Verhandlungen der
k\"oniglich s\"achsischen Gesellschaft der Wissenschaften zu Leipzig, Mathematisch - Physische Classe 53, pp. 1-64.

\smallskip
\noindent
Hurd, A. E. and Loeb, P. (1985): \emph{An Introduction to Nonstandard Real Analysis}, Academic Press, Inc., Orlando.

\smallskip
\noindent
Iglesias-Zemmour, P. (2013): \emph{Diffeology}, Mathematical Surveys and Monographs, Volume 185, American Mathematical Society, Providence RI.

\smallskip
\noindent
Jin, R. (1992): \emph{The isomorphism property versus the special model axiom}, 
The Journal of Symbolic Logic 57, pp. 975-987. 

\smallskip
\noindent
Jin, R. (1997): \emph{Better nonstandard universes with applications}, in \textbf{\emph{Nonstandard analysis (Edinburgh, 1996)}}, 
NATO Advanced Science Institutes Series C: Mathematical and Physical Sciences 493, Kluwer Academic Publishers, Dordrecht, pp. 183-208.

\smallskip
\noindent
Jin, R. and Keisler, H. J. (1993): \emph{Game sentences and ultrapowers}, Annals of Pure and Applied Logic 60, pp. 261-274.

\smallskip
\noindent
J\'onsson, B. (1960): \emph{Homogeneous universal relational structures}, Mathematica Scandinavica 8, pp. 137-142.

\smallskip
\noindent
Kamo, S: (1981): \emph{Nonstandard real number systems with regular gaps}, Tsukuba Journal of Mathematics 5, pp. 21-24.

\smallskip
\noindent
Kamo, S. (1981a): \emph{Nonstandard natural number systems and nonstandard models}, The Journal of Symbolic Logic 46, pp. 365-376.

\smallskip
\noindent
Kanovei, V. and Reeken, M. (2004): \emph{Nonstandard Analysis, Axiomatically}, Springer
Monographs in Mathematics, Springer, Berlin.

\smallskip
\noindent
Kanovei, V. and Shelah, S. (2004): \emph{A definable nonstandard model of the reals}, The Journal of Symbolic Logic 69, pp. 159-164.

\smallskip
\noindent
Kaplansky, I. (1942): \emph{Maximal fields with valuations}, Duke Mathematical Journal 9, pp. 303-321.

\smallskip
\noindent
Katz, M. and Sherry, D. (2013): \emph{Leibniz's infinitesimals: their fictionality, their modern implementations, and their foes from Berkeley to Russell and beyond}, Erkenntnis 78, pp. 571-625.

\smallskip
\noindent
Keisler, H. J. (1963): \emph{Limit ultrpowers}, Transactions of the American Mathematical Society 107, pp. 383-408.

\smallskip
\noindent
Keisler, H. J. (1964): \emph{Ultraproducts and saturated models}, Indagationes Mathematicae 26, pp. 178-186.

\smallskip
\noindent
Keisler, H. J. (1974): \emph{Monotone complete fields}, in \textbf{\emph{Victoria Symposium on Nonstandard Analysis (University of Victoria, Victoria, B.C., 1972)}}, Lecture Notes in Mathematics, Vol. 369, Springer, Berlin, pp. 113-115. 

\smallskip
\noindent
Keisler, H. J. (1976): \emph{Foundations of Infinitesimal Calculus}, Prindle, Weber \& Schmidt, Incorporated, 20 Newbury Street, Boston. Second Edition (2011), available at www.math.wisc.edu/~keisler/foundations.pdf.

\smallskip
\noindent
Keisler, H. J. (1984): \emph{ An infinitesimal Approach to Stochastic Analysis}, Memoirs of the American Mathematical Society, No. 297.

\smallskip
\noindent
Keisler, H. J. (1994): \emph{The hyperreal line}, in \textbf{\emph{Real Numbers, Generalizations of Reals, and Theories of Continua}}, edited by
P. Ehrlich, Kluwer Academic Publishers, Dordrecht, pp. 207-237.

\smallskip
\noindent
Keisler, H. J. and Schmerl, J. (1991): \emph{Making the hyperreal line both saturated and complete}, The Journal of Symbolic Logic 56, pp. 1016-1025. 

\smallskip
\noindent
Kennison, J. (1999): \emph{Synthetic solution manifolds for differential equations}, Journal of Pure and Applied Algebra 143, pp. 255-274.

\smallskip
\noindent
Klingenberg, W. (1954): \emph{Projektive und affine Ebenen mit Nachbarelementen}, 
Mathematische Zeitschrift 60, pp. 388-406. 

\smallskip
\noindent
Klingenberg, W. (1954a): \emph{Euklidische Ebenen mit Nachbarelementen}, Mathematische Zeitschrift 61, pp. 1-25.

\smallskip
\noindent
Klingenberg, W. (1955): \emph{Desarguessche Ebenen mit Nachbarelementen}, Abhandlungen aus dem Mathematischen Seminar der Universit\"at Hamburg 20, pp. 97-111.

\smallskip
\noindent
Klingenberg, W. (1956): \emph{Projektive Geometrien mit Homomorphismus}, Mathematische Annalen 132, pp. 180-200.

\smallskip
\noindent
Kock, A. (1981): \emph{Synthetic Differential Geometry}, London Mathematical Society Lecture Note Series, 333. Cambridge University Press, Cambridge. Second Edition, 2006.

\smallskip
\noindent
Kock, A. (1993): \emph{Review of} \textbf{\emph{Models for Smooth Infinitesimal Analysis}} \emph{by I. Moerdijk and G. Reyes}, The Journal of Symbolic Logic 58, pp. 354-355.

\smallskip
\noindent
Kock, A. (2003): \emph{Differential calculus and nilpotent real numbers}, The Bulletin of Symbolic Logic 9, pp. 225-230.

\smallskip
\noindent
Kock, A. (2010): \emph{Synthetic Geometry of Manifolds}, Cambridge University Press, Cambridge.

\smallskip
\noindent
Kock, A. (25 Sept 2017): \emph{New methods for old spaces: synthetic differential geometry}, preprint, arXiv:1610.00286v3. To appear in \textbf{\emph{New Spaces in Mathematics and Physics}}, edited by M. Anel and G. Catren.

\smallskip
\noindent
Kock, A. and Reyes, G. E. (1981): \emph{Models for synthetic integration theory}, Mathematica Scandinavica 48, pp. 145-152.

\smallskip
\noindent
Krull, W. (1932): \emph{Allgemeine Bewertungstheorie}, Journal f\"ur die Reine und Angewandte Mathematik 167, pp. 160-196.

\smallskip
\noindent
Kreisel, G. (1969): \emph{Axiomatizations of nonstandard analysis that are conservative extensions of formal systems for classical analysis}, in \textbf{\emph{Applications of Model Theory to Algebra, Analysis, and Probability}}, edited by W. A. J. Luxemburg, Holt, Rinehart and Winston, New York, pp. 93-106.

\smallskip
\noindent
Kreuzer, A. (1987): \emph{Hjelmslev-R\"aume}, Results in Mathematics 12, pp. 148-156. 

\smallskip
\noindent
Kreuzer, A. (1988): \emph{ Projective Hjelmslev-R\"aume}, Dissertation TU, Munich.

\smallskip
\noindent
Lakatos, I. (1978): \emph{Cauchy and the continuum: the significance of nonstandard analysis for the history and philosophy of mathematics}, in \textbf{\emph{Imre Lakatos, Philosophical Papers, Volume 2. Mathematics, Science and Epistemology}}, edited and with an introduction by J. Worrall and G. Currie, Cambridge University Press, Cambridge-New York. Reprinted in Mathematical Intelligencer 1 (1978), no. 3, pp 151-161.

\smallskip
\noindent
Laugwitz, D. (1983): \emph{$\Omega$-calculus as a generalization of field extension--an alternative approach to nonstandard analysis}, in \textbf{\emph{Nonstandard Analysis--Recent Developments}}, edited by A. E. Hurd, Springer, Berlin-New York, pp. 120-133.  

\smallskip
\noindent
Laugwitz, D. (1987): \emph{Infinitely small quantities in Cauchy's textbooks}, Historia Mathematica 14, pp. 258-274.

\smallskip
\noindent
Laugwitz, D. (1989):
\emph{Definite values of
infinite sums: aspects of the foundations of infinitesimal analysis
around 1820}, Archive for History of Exact Sciences 39, pp. 195-245.

\smallskip
\noindent
Laugwitz, D. (1992): \emph{Leibniz' principle and omega calculus}, in \textbf{\emph{Le Continu Mathematique, Colloque de Cerisy}},  edited by H. Sinaceur and J-M. Salanskis, Springer-Verlag France, Paris, pp. 144-154. 

\smallskip
\noindent
Laugwitz, D. (2001): \emph{Kurt Schmieden's approach to infinitesimals. An eye-opener to the historiography of analysis}, in \textbf{\emph{Reuniting the Antipodes-Constructive and Nonstandard Views of the Continuum}}, edited by P. Schuster, U. Berger and H. Osswald, Dordrecht, Holland, pp. 127-142.

\smallskip
\noindent
Lavendhomme, R. (1996): \emph{Basic Concepts of Synthetic Differential Geometry}, Kluwer Academic Publishers, Dordrecht.

\smallskip
\noindent
Lawvere, F. W. (1979):
\emph{Categorical dynamics}, in \textbf{\emph{Topos Theoretic Methods in Geometry}}, edited by A. Kock, Various Publications Series, Aarhus University, Aarhus, No. 30, pp. 1-28.

\smallskip
\noindent
Lawvere, F. W. (2011): \emph{Euler's continuum functorially vindicated}, in \textbf{\emph{Logic, Mathematics, Philosophy: Vintage Enthusiasms: Essays in Honor of John L. Bell}}, edited by D. DeVidi, M. Hallet and P. Clark, The Western Ontario Series in Philosophy of Science 75, Springer Science+Business Media, pp. 249-254.

\smallskip
\noindent
Levi-Civita, T. (1892-93): \emph{Sugli infiniti ed infinitesimi attuali quali elementi analitici}. Reprinted in \textbf{\emph{Tullio Levi-Civita, Opere Matematiche, Memorie e Note, Volume primo 1893-1900}},  Nicola Zanichelli, Bologna, 1954.

\smallskip
\noindent
Levi-Civita, T. (1898): \emph{Sui Numeri Transfiniti}. Reprinted in \textbf{\emph{Tullio Levi-Civita, Opere Matematiche, Memorie e Note, Volume primo 1893-1900}}, Nicola Zanichelli, Bologna, 1954.

\smallskip
\noindent
Lindstr\o m, T. (1988): \emph{ An invitation to nonstandard analysis}, in \textbf{\emph{Nonstandard Analysis and its Applications (Hull, 1986)}}, edited by N. Cutland, Cambridge University Press, Cambridge, pp. 1-105.

\smallskip
\noindent
Loeb, P. (1975): \emph{Conversion from nonstandard to standard measure spaces and applications in probability theory},
Transactions of the American Mathematical Society 211, pp. 113-122.

\smallskip
\noindent
Loeb, P. (2000): \emph{An Introduction to Nonstandard Analysis (With an Appendix by Horst Osswald)}, in  \textbf{\emph{ Nonstandard Analysis for the Working Mathematician}}, edited by P. Loeb and M. Wolff, Kluwer Acadademic Publihsers, Dordrecht, pp. 1-95. 

\smallskip
\noindent
Kyle,

Since I will be going over the exam tomorrow, I am willing to give you a make up tomorrow in my office (231S Lindley Hall) at 11:50.

Please let me know if you will be there.

PE
: \emph{A topological characterization of Hjelmslev's classical geometries}, in \textbf{\emph{Rings and Geometry}}, edited by R. Kaya, P. Plaumann and K. Strambach, D. Reidel Publishing Company, Dordrecht, pp. 81-151.

\smallskip
\noindent
\L o\'s, J. (1955): \emph{Quelques remarques, th\'or\`emes, et probl\`emes sur les classes d\'efinissables d'alg\`ebres}, in  \textbf{\emph{Mathematical Interpretations of Formal Systems}}, North Holland, Amsterdam, pp. 98-113.

\smallskip
\noindent
Lutz, R. and Goze, M. (1981): \emph{Nonstandard Analysis, a Practical Guide
with Applications}, Springer-Verlag, Berlin-New York.

\smallskip
\noindent
Luxemburg, W. A. J. (1962): \emph{Non-Standard Analysis, Lectures on A. Robinson's Theory of Infinitesimal and Infinitely Large Numbers}, mimeographed notes, California Institute of Technology Bookstore, Pasadena, CA. Revised 1964.

\smallskip
\noindent
Luxemburg, W. A. J. (1969): \emph{A general theory of monads}, in \textbf{\emph{Applications of Model Theory to Algebra, Analysis, and Probability}}, edited by W. A. J. Luxemburg, Holt, Rinehart and Winston, New York, pp. 18-86.

\smallskip
\noindent
McLarty, C. (1983): \emph{Local, and some global, results in synthetic differential geometry}, in \textbf{\emph{Category theoretic methods in geometry}}, edited by A. Kock, Aarhus University, pp. 226-256.

\smallskip
\noindent
McLarty, C. (1990): \emph{The uses and abuses of the history of topos theory}, The British Journal for the Philosophy of Science 41, pp. 351-375.

\smallskip
\noindent
Moerdijk, I. and Reyes, G. E. (1991): \emph{Models for Smooth Infinitesimal Analysis}, Springer-Verlag, New York.

\smallskip
\noindent
Morley, M. and Vaught, R. (1962): \emph{Homogeneous universal models}, Mathematica Scandinavica 11, pp. 37-57.

\smallskip
\noindent
Mumford, D. (2011): \emph{Intuition and rigor and Enriques's quest}, Notices of the American Mathematical Society 58, pp. 250-260.

\smallskip
\noindent
Neder, L. (1941): \emph{Modell einer Differentialrechnung mit aktual unendlich kleinen Gr\"ossen erster Ordnung}, 
Mathematische Annalen 118, pp. 251-262.

\smallskip
\noindent
Neder, L. (1943): \emph{Modell einer Leibnizischen Differentialrechnung mit aktual unendlich kleinen Gr\"ossen s\"amtlicher Ordnungen}, 
Mathematische Annalen 118, pp. 718-732.

\smallskip
\noindent
Neumann, B. H. (1949): \emph{On ordered division rings}, Transactions of the American Mathematical Society 66, pp. 202-252.

\smallskip
\noindent
 Nowik, T. and Katz, M. (2015): \emph{Differential geometry via infinitesimal displacements}, Journal of Logic and Analysis 7, pp. 1-44.

\smallskip
\noindent
Ozawa, M. (1995): \emph{Scott incomplete Boolean ultrapowers of the real line}, The Journal of Symbolic Logic 60, pp. 160-171.

\smallskip
\noindent
Peirce, C. S. (circa 1897): \emph{Multitude and continuity}, in \textbf{\emph{The New Elements of Mathematics, Volume III/1. Mathematical Miscellanea}}, edited by C. Eisele, Moulton Publishers, The Hague-Paris; Humanities Press, Atlantic Highlands, N.J. 1976, pp. 82-100.

\smallskip
\noindent
 Peirce, C. S. (1898): \emph{The logic of relatives}, in \textbf{\emph{Reasoning and the Logic of Things: The Cambridge Conferences Lectures of 1898}}, edited by K. L. Ketner with an Introduction by K. L. Ketner and H. Putnam, Harvard University Press, Cambridge, MA., 1992, pp. 146-164.

\smallskip
\noindent
Peirce, C. S. (1900): \emph{Infinitesimals}, Science 2, pp. 430-33. Reprinted in \textbf{\emph{Collected Papers of Charles Sanders Peirce, Volume III}}, edited by C. Hartshone and P. Weiss, Harvard University Press, Cambridge, MA., 1935, pp. 360-365.

\smallskip
\noindent
Penon, P. (1981): \emph{Infinit\'esimaux et intuitionnisme}, Cahiers de Topologie et G\'eom\'etrie Diff\'erentielle Cat\'egoriques 22, pp. 67-72.

\smallskip
\noindent
Penon, J. (1985): \emph{De l'infinit\'esimal au local}, Th\`ese de Doctorat d'\'Etat, Universit\'e de Paris VII. Reprinted in \emph{Diagrammes}, Paris VII, 1985. 

\smallskip
\noindent
Perren, D. (2008): \emph{Algebraic Geometry: An Introduction}, Springer-Verlag London Limited, London. 

\smallskip
\noindent
Peterson, J. (= Hjelmslev, J.) (1898): \emph{Nouveau principe pour \'etudes de g\'eom\'etrie des droites}, Oversigt over det K. Danske Vidensk. Selskabs Forhandl, pp. 283-344.

\smallskip
\noindent
Pincherle, S. (1884): \emph{Alcune osservazioni sugli ordini dÕinfinito delle funzioni}, Memorie della
R. Accademia delle scienze di Bologna 5 (Serie 4), pp. 739-750.

\smallskip
\noindent
Plotkin, J. M. (2005): \textbf{\emph{Hausdorff on Ordered Sets}}, edited and translated by J. M. Plotkin, American Mathematical Society, London Mathematical Society, Providence, RI.

\smallskip
\noindent
Poincar\'e, H. (1893/1952): \emph{Le continu math\'ematique}, Revue de m\'etaphysique et de morale 1, pp. 26-34. Reprinted with modifications as Chapter II of \textbf{\emph{Science et l'hypoth\'ese}}, Flammarion, Paris,
1902; English translation: \textbf{\emph{Science and Hypothesis}}, Dover Publications, Inc., New York, 1952.

\smallskip
\noindent
Predella, P. (1912): \emph{Saggio di Geometria non-Archimedeia (Nota II)}, Giornale di matematiche di Battaglini 50 [(3), 3], pp. 161-171.

\smallskip
\noindent
Reyes, G. E. and Wraith, G. C. (1978): \emph{A note on tangent bundles in a category with a ring object}, Mathematica Scandinavica 42, pp. 56-63.

\smallskip
\noindent
Robinson, A. (1961): \emph{Non-standard analysis}, Indagationes Mathematicae 23, pp. 432-440. Reprinted in \textbf{\emph{Abraham Robinson Selected Papers, Volume 2: Nonstandard Analysis and Philosophy}}, edited by J. A. W. Luxemburg and S. K\"orner, Yale University Press, New Haven and London, 1979, pp. 3-11.

\smallskip
\noindent
Robinson, A. (1963): \emph{Introduction to Model Theory and to the Metamathematics of Algebra}, North-Holland Publishing Company, Amsterdam. Second Edition, 1965.

\smallskip
\noindent
Robinson, A. (1966): \emph{Non-standard Analysis}, North-Holland Publishing Company, Amsterdam. 

\smallskip
\noindent
Robinson, A. (1967): \emph{The Metaphysics of the calculus}, in \textbf{\emph{Problems in the Philosophy of Mathematics}}, edited by I. Lakatos, North-Holland Publishing Company, Amsterdam, pp. 28-46. Reprinted in \textbf{\emph{Abraham Robinson Selected Papers, Volume 2: Nonstandard Analysis and Philosophy}}, edited by J. A. W. Luxemburg and S. K\"orner, Yale University Press, New Haven and London, 1979, pp. 537-549.

\smallskip
\noindent
Robinson, A. (1973): \emph{Standard and nonstandard number systems}, The Brouwer memorial lecture 1973, Leiden, April 26, 1973. Nieuw Archief voor Wiskunde (3) 21, pp. 115-133. Reprinted in \textbf{\emph{Abraham Robinson Selected Papers, Volume 2: Nonstandard Analysis and Philosophy}}, edited by J.A.W. Luxemburg and S. K\"orner, Yale University Press, New Haven and London,  1979, pp. 426-444.

\smallskip
\noindent
Robinson, A. (1974): \emph{Non-standard Analysis, Second Edition}, North-Holland Publishing Company, Amsterdam.

\smallskip
\noindent
Ross, D. (1990): \emph{The special model axiom in nonstandard analysis}, The Journal of Symbolic Logic 55, pp. 1233-1242.

\smallskip
\noindent
Rubinstein-Salzedo, S. and Swaminathan, A. (2014):  \emph{Analysis on surreal numbers}, Journal of Logic and Analysis 6, Paper 5, pp 1-39. 

\smallskip
\noindent
Sanders, S. (2018): \emph{To be or not to be constructive, that is not the question}, Indagationes Mathematicae (New Series) 29, pp. 313-381.

\smallskip
\noindent
Schmieden, C. and Laugwitz, D. (1958): \emph{Eine Erweiterung der Infinitesimalrechnung}, Mathematische Zeitschrift 69, pp.1-39.

\smallskip
\noindent
Scott, D. (1969): \emph{Boolean models and nonstandard Analysis}, in \textbf{\emph{Applications of Model Theory to Algebra, Analysis, and Probability}}, edited by W. A. J. Luxemburg, Holt, Rinehart and Winston, New York, pp. 87-92.

\smallskip
\noindent
Scott, D. (1969a): \emph{On completing ordered fields}, in \textbf{\emph{Applications of Model Theory to Algebra, Analysis, and Probability}}, edited by W. A. J. Luxemburg, Holt, Rinehart and Winston, New York, pp. 274-278.

\smallskip
\noindent
Segre, C. (1911): \emph{Le geometrie proiettive nei campi di numeri duali}, Atti della R. Academia della Scienze di Torino 47, pp. 114-133, 164-185.

\smallskip
\noindent
Shafarevich, I. R. (2005): \emph{Basic Notions of Algebra}, Springer-Verlag, Berlin. Reprinting of the 1990 edition.

\smallskip
\noindent
Shapiro, S. (2014): \emph{Varieties of Logic}, Oxford University Press, Oxford.

\smallskip
\noindent
Shulman, M. (2013): \emph{The surreal numbers}, in \textbf{\emph{Homotopy Type Theory: Univalent Foundations of Mathematics}}, The Univalent Foundations Program, \url{https://homotopytypetheory.org/book}, Institute for Advanced Study, Princeton, N. J., Section 11.6, pp. 401-411.

\smallskip
\noindent
Siegel, A. N. (2013): \emph{Combinatorial Game Theory}, American Mathematical Society, Providence, RI.

\smallskip
\noindent
Skolem, T. (1934): \emph{\"Uber die Nichtcharakterisierbarkeit der Zahlenreihe mittels endlich oder abz\"ahlbar unendlich vieler Aussagen mit ausschliesslich Zahlenvariablen}, Fundamenta Mathematica 23, pp. 150-161. Reprinted in \textbf{\emph{Selected Works in Logic by T. H. Skolem}}, edited by J. E. Fenstad, Universitetsforlaget, Oslo, 1970, pp. 355-366.

\smallskip
\noindent
Smullyan, R. and Fitting, M. (2010): \emph{Set theory and the Continuum Problem}, Dover Publications, New York.

\smallskip
\noindent 
Stolz, O. (1879): \emph{Ueber die Grenzwerthe der Quotienten}, Mathematische Annalen 
14, pp. 231-240.

\smallskip
\noindent
Stolz, O. (1883): \emph{Zur Geometrie der Alten, insbesondere \"uber ein Axiom des Archimedes}, Mathematische Annalen 22, pp. 504-519.

\smallskip
\noindent
Stolz, O. (1884): \emph{Die unendlich kleinen Gr\"ossen}, Berichte des Naturwissenschaftlich-Medizinischen Vereines in Innsbruck 14, pp. 21-43.

\smallskip
\noindent
Stolz, O. (1885): \emph{Vorlesungen \"uber Allgemeine Arithmetik, Erster Theil: Allgemeines und Arithmetik der Reelen Zahlen}, Teubner, Leipzig.

\smallskip
\noindent
Stroyan, K. (1977): \emph{Infinitesimal analysis of curves and surfaces}, in \textbf{\emph{Handbook
of Mathematical Logic}}, edited by J. Barwise, North-Holland Publishing Company, Amsterdam, pp. 197-231.

\smallskip
\noindent
Stroyan, K. and Luxemburg, W. A. J. (1976): \emph{Introduction to the Theory
of Infinitesimals}, Academic Press, New York.

\smallskip
\noindent
Study, E. (1903): \emph{Geometrie der Dynamen}, Teubner, Leipzig.

\smallskip
\noindent
Tarski, A. (1939): \emph{The completeness of elementary algebra and geometry}, in  \textbf{\emph{Centre National de la Recherche Scientifique, Institut Blaise Pascal}}, Paris, 1967, pp. 289-346. Reprinted in \textbf{\emph{The Collected Papers of Alfred Tarski, Vol. IV}}, edited by S. R. Givant and R. N. McKenzie, Birkh\"auser, 1986.

 \smallskip
\noindent
Tarski, A. (1948): \emph{A Decision Method for Elementary Algebra and Geometry}, Rand Corporation, Santa Monica, CA. Second Edition (prepared for publication with the assistance of J. C. C. McKinsey) University of California Press, Berkeley and Los Angeles, 1951. 

\smallskip
\noindent
Tarski, A. (1959): \emph{What is elementary geometry?} in \textbf{\emph{The Axiomatic Method, with Special Reference to Geometry and Physics}}, edited by L. Henkin, P. Suppes  and A. Tarski, North-Holland Publishing Company, Amsterdam, pp. 16-29.

\smallskip
\noindent
Tao, T. (2008): \emph{Structure and Randomness. Pages from Year One of a Mathematical
Blog}, American Mathematical Society, Providence, RI. 

\smallskip
\noindent
Tao, T. (2013): \emph{Compactness and Contradiction}, American Mathematical Society, Providence, RI. 

\smallskip
\noindent
Tao, T. (2014): \emph{Hilbert's Fifth Problem and Related Topics}, Graduate Studies in Mathematics, American Mathematical Society, Providence, RI. 
153. 

\smallskip
\noindent
Thomae, J. (1870): \emph{Abriss einer Theorie der complexen Functionen und der Thetafunctionen einer
Ver\"anderlichen}, Verlag von Louis Nebert, Halle.

\smallskip
\noindent
van den Dries, L. and Ehrlich, P. (2001): \emph{Fields of surreal numbers and exponentiation}, Fundamenta Mathematica 167, pp. 173-188; erratum, ibid. 168 (2001), pp. 295-297. 

\smallskip
\noindent
van den Dries, L. and Ehrlich, P. (forthcoming): \emph{Homogeneous universal $H$-fields}, Proceedings of the American Mathematical Society; preprint, arXiv:1807.08861. 

\smallskip
\noindent
van den Dries, L.. Macintyre, A. and Marker, D. (2001): \emph{Logarithmic-exponential series}, Annals of Pure and Applied Logic 111, pp. 61-113. 

\smallskip
\noindent
van der Hoeven, J. (2006): \emph{Transseries and Real Differential Algebra}, Springer-Verlag, Berlin.

\smallskip
\noindent
Veldkamp, F. D. (1981): \emph{Projective planes over rings of stable rank $2$}, Geometriae Dedicata 11, pp. 285-308. 

\smallskip
\noindent
Veldkamp, F. D. (1985): \emph{Projective ring planes and their homomorphisms}, in \textbf{\emph{Rings and Geometry}}, edited by R. Kaya, P. Plaumann and K. Strambach, D. Reidel Publishing Company, Dordrecht, pp. 289-349.

\smallskip
\noindent
Veldkamp, F. D. (1987): \emph{Projective Barbilian spaces}, Results in Mathematics 12, pp. 222-240, 434-449.

\smallskip
\noindent
Veldkamp, F. D. (1995): \emph{Geometry over rings},  in \textbf{\emph{Handbook of Incidence Geometry}}, edited by F. Buekenhout, North-Holland, Amstredam, pp. 1033-1084.

\smallskip
\noindent
Veronese, G. (1889): \emph{Il continuo rettilineo e l'assioma V di Archimede}, Memorie della Reale
Accademia dei Lincei, Atti della Classe di scienze naturali, fisiche e matematiche (4), 6,
pp. 603Ð-624.
 
 \smallskip
\noindent
Veronese, G. (1891): \emph{Fondamenti di geometria a pi dimensioni e a pi specie di unit rettilinee esposti in forma elementare}, Padova: Tipografia del Seminario.

\smallskip
\noindent
Veronese, G. (1894): \emph{Grundz\"uge der Geometrie von mehreren Dimensionen und mehreren Arten gradliniger Einheiten in elementarer Form entwickelt. Mit Genehmigung des Verfassers nach einer neuen Bearbeitung des Originals \"ubersetzt von Adolf Schepp}, Teubner, Leipzig.

\smallskip
\noindent
Veronese, G. (1909/1994): \emph{La geometria non-Archimedea}, translated by M. Marion as \emph{On non-Archimedean geometry}, in  \textbf{\emph{Real Numbers, Generalizations of the Reals, and Theories of Continua}}, edited by P. Ehrlich, Kluwer Academic Publishers, Dordrecht, Netherlands, 1994. First published in \textbf{\emph{Atti del IV Congresso Internazionale dei Matematici (Roma 6-11 Aprile 1908), Volume I}}, 1909, pp. 197-208.

\smallskip
\noindent
Warner, S. (1965): \emph{Modern Algebra, Volumes I and II}, Prentice Hall, Inc., Englewood Cliffs, N J. Reprinted by Dover Publishing Company, New York, 1990.

\smallskip
\noindent
Weil, A. (1953): \emph{Th\'eorie des points proches sur les vari\'et\'es diff\'erentiables}, in \textbf{\emph{G\'eom\'etrie diff\'erentielle}}, Colloques Internationaux du Centre National de la Recherche Scientifique, Strasbourg, CNRS, Paris, pp. 111-117.

\smallskip
\noindent
Zakon, E. (1969): \emph{Remarks on the nonstandard real axis}, in \textbf{\emph{Applications of Model Theory to Algebra, Analysis, and Probability}}, edited by W. A. J. Luxemburg, Holt, Rinehart and Winston, New York, pp. 195-227.

\smallskip
\noindent
Zemmer, J. L. (1953): \emph{Ordered algebras that contain divisors of zero}, Duke Mathematical Journal 20, pp. 177-183.

\end{footnotesize}	

\end{document}